\documentclass[authoryear,preprint,11pt]{elsarticle}

\usepackage{amsthm,amsmath,amssymb,eurosym}
\usepackage{color}
\usepackage{multirow}
\usepackage{longtable}
\usepackage{arydshln}
\usepackage{geometry}
\usepackage{enumitem}
\usepackage{graphicx}
\usepackage{booktabs}
\usepackage[table,xcdraw]{xcolor}
\usepackage{tabularx}
\usepackage[para]{threeparttable}
\usepackage{adjustbox}
\usepackage{float}
\usepackage{caption}
\usepackage{subcaption}
\usepackage[utf8]{inputenc}
\usepackage{hyperref}
\usepackage{comment}

\usepackage{tikz}
\usepackage{pgfplots}
\usetikzlibrary{patterns}
\usetikzlibrary{arrows}
\usetikzlibrary[arrows,snakes,shapes,shapes.misc,backgrounds,trees,shadows,positioning,fit,petri,patterns]
\usetikzlibrary[matrix,decorations]
\usetikzlibrary[calc,through,intersections]
\usetikzlibrary[calendar,mindmap,scopes]
\usetikzlibrary{decorations.markings,%
    shapes.arrows,chains,matrix,positioning,decorations.pathmorphing,
    shapes,backgrounds,decorations.text}
\usetikzlibrary{decorations,calc}
\usetikzlibrary{decorations,calc}
\usepackage{tikz-3dplot}
\usepackage{tikz-qtree}
\usepackage{pict2e,picture}

\newtheorem{remark}{Remark}

\journal{European Journal of Operational Research}

\usepackage{lineno, blindtext}

\begin{document}

\begin{frontmatter}
    \title{A multiple criteria approach for bulding a pandemic impact assessment composite indicator: The case of {\sc{Covid-19}} in Portugal}
    \author[cegist]{Jos\'e Rui {\sc  Figueira}\corref{cor1}}\ead{figueira@tecnico.ulisboa.pt}
    \author[cemat]{Henrique M. {\sc Oliveira}}
    \author[cqe]{Ana Paula {\sc Serro}}
    \author[idmec]{Rogério {\sc  Colaço}}
    \author[gcom]{\\Filipe {\sc  Froes}}
    \author[gcom]{Carlos {\sc  Robalo Cordeiro}}
    \author[gcom]{António {\sc  Diniz}}
    \author[gcom]{Miguel {\sc  Guimarães}}
    \address[cegist]{CEGIST, Instituto Superior T\'{e}cnico,  Universidade de Lisboa, Portugal}
    \address[cemat]{CAMGSD, Instituto Superior T\'{e}cnico,  Universidade de Lisboa, Portugal}
    \address[cqe]{CQE, Instituto Superior T\'{e}cnico,  Universidade de Lisboa, Portugal}
    \address[idmec]{IDMEC, Instituto Superior T\'{e}cnico,  Universidade de Lisboa, Portugal}
    \address[gcom]{GCOM, Ordem dos Médicos, Portugal}
    \cortext[cor1]{Corresponding author at: CEGIST, Instituto Superior T\'{e}cnico,  Av. Rovisco Pais, 1, 1049-001 Lisboa, Portugal.}
    \begin{abstract}
    \noindent  The {\sc{Covid-19}} pandemic has caused major damage and disruption to social, economic, and health systems (among others). In addition, it has posed unprecedented challenges for public health and policy/decision makers who have been responsible for designing and implementing measures to mitigate its strong negative impact. The Portuguese health authorities have used decision analysis-like techniques to assess the impact of the pandemic and implemented measures for individual counties, regions, or across the whole country. These decision tools have been subject to some criticism and many stakeholders asked for novel approaches, in particular those which took into account the dynamic changes in the pandemic's behaviour, for example, as a result of new virus variants or vaccines. A multidisciplinary team formed by researchers from the {\sc{Covid-19}} Committee of Instituto Superior Técnico at Universidade de Lisboa (CCIST analyst team) and physicians from the Crisis Office of the Portuguese Medical Association (GCOM expert team) joined forces and worked together to create a new tool to help politicians and decision-makers to fight the pandemic. This paper presents the main steps that led to the building of a pandemic impact assessment composite indicator applied to the specific case of {\sc{Covid-19}} in Portugal. A multiple criteria approach based on an additive multi-attribute value theory aggregation model was used to build the pandemic assessment composite indicator. The parameters of the additive model were devised based on a sociotechnical co-constructive interactive process between the CCIST and GCOM team members. The deck of cards method was the technical tool adopted to help in building the value functions and the assessment of the criteria weights. The final tool was presented at a press conference and had a strong impact in the Portuguese media and on the main health decision-making stakeholders in the country.
    \end{abstract}
    \vspace{0.25cm}
    \begin{keyword}
        Multiple criteria analysis  \sep  Composite indicator  \sep Multi-attribute value theory (MAVT) \sep Robustness and validation analyses \sep Deck of cards method \sep Multidisciplinary team work.
    \end{keyword}
\end{frontmatter}

\vfill\newpage

\tableofcontents

\vfill\newpage

\section{Introduction}\label{sec:introduction}
\noindent A pandemic causes major damage and disruptions to social, economic, and health systems (among others) and has major implications on the lives of populations throughout the world. Not only has it leads to very serious physical and mental health problems, but also to poverty and hunger. {\sc{Covid-19}}, the most recent pandemic poses  unprecedented challenges to public health and for policy/decision-makers, especially when designing and implementing measures to mitigate the pandemic's negative impacts.

Different countries and researchers around the world have presented tools for mitigating the impact of {\sc{Covid-19}}. Two types of the literature review can be presented: on the one hand an analysis of the tools used by other countries, and on the other hand, a review of the published literature in the field. One of the most widely adopted tools for assessing the impact of {\sc{Covid-19}} has been the use of chromatic systems, especially the use of territorial unit (counties, districts, municipalities, departments, provinces, regions, states, and the entire countries) risk maps. These chromatic systems assign a colour to each territorial unit, which represents the risk from the lowest to the highest level (i.e., from a lighter to a darker colour or from green to red hues) of that particular unit. The methodologies and factors taken into account to obtain the risk maps are somewhat different from case to case, as can be seen is the following examples:

\begin{enumerate}
    \item	In Spain\footnote{https://flowmaps.life.bsc.es/flowboard/board\_what\_is\_risk}, a risk flow-map is produced based on a score computed from three criteria (indicators): the daily and the cumulative incidences, and the population mobility partners. The score is used to estimate the number of cases that can be exported or imported between pairs of regions.
    \item In Italy\footnote{https://www.salute.gov.it/portale/nuovocoronavirus/dettaglioContenutiNuovoCoronavirus.jsp?lingua \=english\&id=5367\&area\=nuovoCoronavirus\&menu\=vuoto}, the risk classification of each region is based on ordinances issued by the Ministry of Health. The risk is determined from the impact of several factors and a probability is associated with each impact level. A risk matrix (RM) is formed with different impact levels.
    \item In France\footnote{https://sante.journaldesfemmes.fr/fiches-maladies/2667643-carte-covid-france -europe-voyage-monde-pays-zone-rouge-epidemie-contamination/}, a {\sc{Covid-19}} map was designed according to the administrative divisions (department) to model the progression of the pandemic based on the seven day incidence of {\sc{Covid-19}} \textit{per} $100\;000$ inhabitants.
    \item In Germany\footnote{https://www.rki.de/DE/Content/InfAZ/N/Neuartiges\_Coronavirus/Risikobewertung\_Grundlage.html}, the Robert Koch Institute’s coloured maps are only based on the number of {\sc{Covid-19}} cases and cumulative incidence (per $100\;000$ inhabitants) reported by each county/federal state.
    \item In the United Kingdom (UK)\footnote{https://coronavirus.data.gov.uk/details/interactive-map/cases}, the map shows the seven-day case rate \textit{per} $100\;000$ inhabitants, and like German, it is a single criterion-based assessment tool.
    \item In the United States of America (USA)\footnote{USA (North Carolina): https://covid19.ncdhhs.gov/dashboard}, in particular in North Carolina, the county map considers the number of cases \textit{per} $100\;000$  residents.
    \item In Canada\footnote{https://health-infobase.canada.ca/covid-19/}, the map is similar to number 6: the provinces and territories are coloured according to the number of cases over the past seven days.
    \item In Brazil\footnote{https://coronavirus.es.gov.br/mapa-de-gestao-de-risco}, the risk map is based on several factors including the number of active cases, number of tests, and lethality and it is then aggregated into two dimensions (threats and vulnerabilities); this results in a RM with several levels being built.
\end{enumerate}

In summary, the tools presented by Brazil and Italy are closest to that developed by the {\sc{Covid-19}} Crises Office of the Portuguese Medical Association (PMA), i.e., a first step for the construction of a multi-criteria decision aiding/analysis (MCDA) tool.

Although not officially adopted, a significant number of other decision support tools have recently been developed for pandemic mitigation policy making purposes. For instance, \cite{Haghighat2021} proposed a combined multilayer perception neural network and Markov chain approach for predicting the number of future patients and deaths in the Bushehr province, Iran. \cite{CatalaEtAl2021} proposed three risk indicators to estimate the status of the pandemic and applied them to the evolution of different European countries. These indicators quantify both the propagation and the number of estimated cases. \cite{NelkenEtAl2020} conducted a review of the different {\sc{Covid-19}} indicators proposed in the literature and explore the social role of these indicators in the pandemic from different perspectives. \cite{HaleEtAl2021} presented databases and composite indicators analysing the effect of policy responses on the spread of {\sc{Covid-19}} cases and deaths and on economic and social welfare. In their study, the composite indicators are simple formulas, which aggregate several partial indicators of a both qualitative and quantitative nature. The conversion of the qualitative nature of the scale levels by assigning numbers is very questionable. \cite{PangEtAl2021} carried out a study on risk environmental assessment in the Hubei province of China and put forward a composite indicator of the disaster loss for {\sc{Covid-19}} transmission. This indicator is based on five environmental perspectives and $38$ partial indicators. Statistical and component analysis methods were used to analyse and build the indicators.

Some indicators are also related to risk, vulnerability or other impact concepts maps.  \cite{NeyensEtAl2020} proposed a statistical-based method to assess the risk map of each Belgium municipality by making use of spatial data on {\sc{Covid-19}} gathered from a large online survey. This study enables predictions to be made on the incidence of the disease and establish a comparison (analysing the proportion of heterogeneity) with respect to the number of confirmed cases. \cite{LiEtAl2021} presented a risk analysis of the {\sc{Covid-19}} infection (modelled using the classic impact $X$ probability formula) of the different regions of China, from the Wuhan region to the other $31$ regions. The authors use the high-speed rail network to assess and predict the regional risk of infection of each region. \cite{DlaminiEtAl2020} proposed several risk assessment indicators for identifying the risk areas in Eswatini, Iran. The risk overall indicators use socio-economic and demographic partial indicators. \cite{SarkarCh2021} presented a socio-environmental vulnerability indicator of the potential risk of community spread of {\sc{Covid-19}}. The overall composite indicator was built from the four most influent socio-economic and environmental partial indicators selected through principal component analysis. It was then applied to assess the vulnerability risk of each district of India. \cite{GhimireEtAl2021} also proposed indicators for {\sc{Covid-19}} risk assessment with geo-visualisation map tools applied to Nepal. The composite indicator results from a weighted-sum which takes into account a positive case score, a quarantined people score, a community exposure score, and a population density score.

Among the papers reviewed, three were of particular interest: two of a multi-criteria nature and one of a single criterion nature, but applied to the Portuguese scenario:

\begin{enumerate}
    \item \cite{SangiorgioPa2020} presented a very interesting composite indicator for the prediction of the risk of contagion risk in urban districts of the Apulia region in Italy. The indicator considers relevant socio-economic data from three perspectives: hazard, vulnerability, and exposure. Each of them comprises several dimensions that are normalised and weighted. The composite indicator is based on a factorial formula, which is calibrated through an optimisation procedure.
    \item \cite{ShadeedAl2021} introduced a multi-criteria index for estimating vulnerability. This is based on the analytical hierarchy process method and was applied in the Governorates of Palestine. The criteria considered are the following: population, population density, elderly population, accommodation and food service activities, school students, chronic diseases, hospital beds, health insurance, and pharmacy.
    \item \cite{AzevedoEtAl2020} presented an interesting indicator for infection risk assessment applied to each municipality on the Portuguese mainland. The indicator is based on the daily number of infected people and uses a direct block sequential simulation.
\end{enumerate}

Despite the relevance of these proposal for policymaking, they do not assess, however, the impact of {\sc{Covid-19}} in terms of activity and severity. Building a model of the {\sc{Covid-19}} impact using a composite indicator is a conceptual activity which can provide a means of observing the evolution of the pandemic; it may also be an important tool for policymaking.

In Portugal a team from the {\sc{Covid-19}} Crisis Office of PMA (GCOM experts team) and a team from the {\sc{Covid-19}} Committee of Instituto Superior T\'ecnico (CCIST analysts team), joined forces after an initial period in which they acted separately to help in fight against {\sc{Covid-19}} and mitigate its negative impact on people's lives.

A ``Risk Matrix''(RM) tool (Figure \ref{fig: OM_HA_Matrix} in the Appendix of this paper, subsequently altered to accommodate some \textit{ad hoc rules}) has been used by the Portuguese health authorities to help in the pandemic decision-making process. This tool attributes a colour coded risk status to each county. It has, however, been subject to some criticism, mainly due to it being incomplete and unable to provide an adequate idea of the pandemic's evolution in the country. RM, in this context, has a meaning different from the well known decision aiding tool of the same name in the field of Decision Analysis. The term ``matrix'' is also unrelated to the mathematical concept. Here, RM is a two-dimensional (2D) referential accompanied by a visual chromatic system (from light green to dark red), where the abscissa axis represents the rawdata on the transmission rate ($R(t)$) and the ordinate axis the average incidence of new positive cases over the past seven days \textit{per} $100\; 000$ inhabitants. In addition, two cut-off lines (one horizontal and one vertical) are used as criteria for separating the referential in four regions: the southwest region with the lowest risk impact; the northeast region with the highest risk impact; and, the other two regions (northwest and southeast), the regions with intermediate risk impact. The different territorial units (the counties and regions) were coloured according to this system and some measures were assigned to each colour. It is important here, in summary, to highlight that the concept of RM and criticism of this system in a different context (i.e., in Decision Analysis, as mentioned previously) can be seen in \cite{Cox2008}.

The RM main drawbacks can be succinctly presented as follows:
\begin{enumerate}
  \item Despite the usefulness and advantages of the visual chromatic system for communication purposes, it suffers from a major pitfall, which renders very difficult to see the evolution of the pandemic over the time-line (plotting each daily situation in the referential and linking all the successive points by a line leads to a very confusing evolution curve, see Figure \ref{fig: RM_Evolution} in the Appendix).
  \item Despite the importance of the dimensions used in the referential (incidence and transmission), they are only part of the problem. Both are related to the activity of the pandemic, but more dimensions should be considered, especially those related to the severity of the pandemic.
  \item The impact of risk does not seem to be appropriately modelled since moving from an $R(t) = 0.1$ to an $R(t) = 0.2$ has the same risk impact as moving from an $R(t) = 0.9$ to an $R(t) = 1.0$, which does not represent how the population, in general, feels about the impact of the transmission (the same reasoning can be applied to incidence).
  \item There is no differentiation between the contributions of incidence and transmission to the overall risk impact; they count equally (this can be acceptable, but it is not always the case).
\end{enumerate}

The CCIST team had previously presented an improved RM containing more cut-off lines than the one used by the Portuguese health authorities. It enabled a closer analysis of the situation to be made, but still suffered from the same drawbacks as the original RM. In parallel the GCOM team also proposed an improved RM (see Figure \ref{fig: RM_OM} in the Appendix) different from the one used by CCIST team. To overcome some of the drawbacks of the original RM, the GCOM team recommended the use of a 2D referential system but one that considered several indicators in both the abscissa axis and the ordinate axis. This was a first important step towards an MCDA-based indicator. The two sets of dimensions (called ahead pillars) on this new RM are the ``activity'' and the ``severity'' of the pandemic. Unfortunately, the way the activity and the severity indicators were considered was questionable, and although this revised model led to a more complete and finer analysis of the problem, it also had some drawbacks (1), (3), and (4). This GCOM proposal was made public in the first week of June 2021. At the beginning of July 2021, the two teams (CCIST and GCOM), began working together to propose the composite indicator presented in this paper. This new proposal had a strong impact in Portugal, especially in the \textit{media} and among health policy and decision-makers. At the end of July 2021, the RM used by the Portuguese health authorities was changed to include some \textit{ad hoc rules}, based on the severity aspects of the pandemic, such as the one proposed in our pandemic assessment composite indicator (PACI).

What was missing in the proposed RM approaches? In short, a more adequate system was needed to characterise the pandemic impact and to recommend the most suitable measures to mitigate its impact. Therefore, the \textit{main decision problem} we faced was how to build a state indicator of the pandemic's impact for a given territorial unit (country, region, county, etc), with the purpose of assigning mitigating measures and/or recommendations for each state (the least to the most restrictive ones). This paper does not, however, present the measures for each state since this is a matter for the Portuguese health authorities and varies with time. An \textit{essential observation} is that it is essential to follow the recent evolution of the pandemic's impact for better planning when a given territory unit moves from a given state to another one. In addition, it is extremely important to assess the impact of the \textit{Portuguese vaccination plan}. Each territorial unit is assessed on a daily basis, taking into account a set of criteria (also called, in our case, indicators or dimensions) grouped in two perspectives or pillars: the activity and the severity of the pandemic. This problem statement can be viewed as belonging to the field of MCDA. For more details, the reader can consult \cite{BeltonSt2002} and \cite{Roy1996}. The problem is known in the literature as an ordinal classification (or sorting) MCDA problem \citep{DoumposZo2002,ZopounidisDo2002}. There are several ways of building a state composite (aggregation) indicator \citep[see][for a recent survey in this topic]{ElGibariEtAl2019}. The main MCDA approaches for designing composite indicators are as follows:

      \begin{enumerate}
        \item \textit{Scoring-based approaches}, as for example, multi-attribute utility/value theory (MAUT/MA\-VT) aggregation models \citep[e.g.,][]{Dyer2016,KeeneyRa1993,TsoukiasFi2006}, analytical hierarchy process \citep[e.g.,][]{Saaty2016}, fuzzy sets techniques \citep[e.g.,][]{DuboisPe2016}, and fuzzy measure based aggregation functions \citep{GrabischLa2016}.
        \item \textit{Outranking-based approaches}, as for example, {\sc{Electre}} methods \citep{FigueiraEtAl2016}, {\sc{Promethee}} methods, \citep{BransDe2016}, and other outranking techniques \citep{MartelBe2016}.
        \item \textit{Rule-based systems}, as for example, decision rule methods \citep{GrecoEtAl2016}, and  verbal decision analysis \citep{MoshkovichEtAl2016}.
      \end{enumerate}

According to the definition of our decision problem, outranking-based methods and rule-based systems are powerful MCDA techniques for ordinal classification and could be adequate tools for building a state composite indicator (ordinal scale). However, the need to analyse the evolution of the pandemic (see the fundamental observation stated before) requires the construction of a richer scale, of a cardinal nature. It is true that outranking-based methods and rule-based systems can be adapted to produce such a cardinal scale (see, for example, \citealt{FigueiraEtAl2021}), but this is a complex process, which is harder to explain to the main actors and the public in general. Consequently, the most adequate approach for dealing with our problem was a scoring-based approach. Since the model needed to be simple enough for interaction and communication with the experts and the public in general, without losing sight of the reality it sought to represent, our study focused solely on MAVT methods. More complex scoring-based methods were discarded. After a more in depth analysis, we finally decided to keep the simple additive MAVT approach. It was suitable for modelling our problem and rendered the communication easy. At this point, another question arose: how should we built the additive model? There were two possible answers:

    \begin{enumerate}
      \item Through a constructive learning approach (machine learning like approaches), such as UTA type methods \citep{SiskosEtAl2016}, or an adaptation of more sophisticated techniques as the GRIP method \citep{FigueiraEtAl2009} with representative functions.
      \item Through a co-constructive sociotechnical interactive process between analysts and policy/de\-ci\-sion-ma\-kers or experts using, for example, the classical MAVT method \citep{KeeneyRa1993}, the MACBETH method \citep{BanaEtAl2016}, or the deck of cards method \citep{Corrente2021}.
    \end{enumerate}

In every co-constructive sociotechnical process, the analyst must be familiar with the technicalities of the method. In addition, the policy/decision-makers or experts must understand the basic questions for assessing their judgements. The improved version of the deck of cards method by \cite{Corrente2021} was found to be an adequate tool. Its adequacy comes from some important aspects: time limitation to produce a meaningful indicator, easy to be understood by the experts, easy to communicate with the public in general, and easy to reproduce the calculations for a reader with an elementary background in mathematics.

In this study we apply MAVT theory through the improved deck of cards method (DCM) \citep{Corrente2021} to: the construction of a cardinal impact assessment composite indicator of the pandemic. The objectives were two-fold: on the one hand to observe the evolution of the pandemic and on the other hand to form a state ordinal indicator with measures and/or recommendations associated with each state to be applied in the case of {\sc{Covid-19}} in Portugal.

The paper is organised as follows. Section 2 introduces the basic mathematical concepts required throughout the paper. Section 3 presents the main models used to perform this study (criteria model, aggregation model, and graphical model). Section 4 displays lessons learn from practice, including the successful aspects, failures and improvements to the tool. Finally, Section 5, outlines the main conclusions and some avenues for future research.


\section{Concepts, definitions, and notation}\label{sec:concepts}
\noindent This section introduces the main concepts, definitions, and notation used along the paper. It comprises the criteria model basic data,  the MAVT additive model, and the chromatic classification system.

\subsection{Basic data}\label{sec:basic_data}
\noindent The basic data can be introduced as follows. Let, $T = \{t_1,\ldots,t_i,\ldots,t_m\}$, denote a set of actions or \textit{time periods} (in general, days) used for observing the pandemic state in a given territory unit (country, region, district, etc), and, $G = \{g_1,\ldots,g_j,\ldots,g_n\}$, denote the set of relevant \textit{criteria} (our problem dimensions or indicators) identified with the experts for assessing the actions or time periods. The \textit{performance} $g_j(t_i) = x_{jt_i} \in E_j$ represents the impact level of activity or severity  over the action or time period $t_i \in T$, according to criterion $g_j$, being $E_j$ the (continuous or discrete) scale of this criterion, for $j=1,\ldots,n$. We will assume, without  any loss of generality, that, for each criterion, the higher the performance level, the higher the impact on the pandemic. The set of criteria has been built according to certain desirable properties \citep[see][]{Keeney1992}.

\subsection{The multi-attribute value theory additive model}\label{sec:additive_model}
\noindent The proposed model is a conjoint analysis model (see, for example, \citealt{BouyssouPi2016}), more specifically an additive MAVT model. The origins of this type of models dates back to 1969, with the seminal work by H. Raiffa, only published in 2016, in \cite{TsoukiasFi2006}, with several comments from prominent researchers in the area. For more details about the additive model see \cite{KeeneyRa1993}.

Let $\succsim$ denote a comprehensive binary relation, over the actions in $T$, whose meaning is ``impacts at least as much as'',. Thus, an action $t^{\prime}$ is considered to impact at least as much as an action $t^{\prime\prime}$, denoted $t^{\prime} \succsim t^{\prime\prime}$, if and only if, the overall value of $t^{\prime}$, $v(t^{\prime})$ is greater than or equal to the overall value of $t^{\prime\prime}$, $v(t^{\prime\prime})$, i.e., $v(t^{\prime}) \geqslant v(t^{\prime\prime})$, where the overall value of each action is additively computed as follows:

\begin{equation}
\label{eq:add_model}
    v(t) = \sum_{j=1}^{n}w_jv_j(x_{jt}), \;\, \mbox{for all}\;\, t \in T
\end{equation}

\noindent in which $w_j$ is the weight of criterion $j$, for $j=1,\ldots,n$, (assuming that $\sum_{j=1}^{n}w_j =1$), and $v_j\big(x_{jt}\big)$ is the value of the performance $x_{jt}$ on criterion $g_j$, for all for $j=1,\ldots,n$.

The asymmetric part of the relation, $t^{\prime} \succ t^{\prime\prime}$, means that $t^{\prime}$ is considered to impact strictly more than  $t^{\prime\prime}$, while the symmetric of the relation, $t^{\prime} \sim t^{\prime\prime}$, means that $t^{\prime}$ is considered to impact equally as  $t^{\prime\prime}$. The three relations $\succsim$, $\succ$, and $\sim$ are transitive.

The construction of the value function, $v_j(x_{jt})$, for criterion $g_j$ and each action or time period $t \in T$, is done in such a way that its value increases with an increasing of the performances level of criterion $j$, $j=1,\ldots, n$ (this function is a non-decreasing monotonic function). Let $t^{\prime}$ and $t^{\prime\prime}$ denote two actions. The following conditions must be fulfilled:

\begin{enumerate}
    \item The strict inequality $v_j(x_{jt^{\prime}}) > v_j(x_{jt^{\prime\prime}})$ holds, if and only if, the impact of performance $x_{jt^{\prime}}$ is considered strictly higher than the impact of performance $x_{jt^{\prime\prime}}$, on criterion $g_j$ (it means that, $t^{\prime}$ impacts strictly more than $t^{\prime\prime}$), for $j=1,\ldots,n$.
    \item The equality $v_j(x_{jt^{\prime}}) = v_j(x_{jt^{\prime\prime}})$ holds, if and only if, the performance $x_{jt^{\prime}}$ impacts the same as the performance $x_{jt^{\prime\prime}}$, on criterion $g_j$, (it means that $t^{\prime}$ impacts equally as $t^{\prime\prime}$), for $j=1,\ldots,n$.
\end{enumerate}

In addition, the value functions are also used for modelling the impact of the performance differences. The higher the performance difference, the higher the strength of the value function impact. Let $t^{\prime}$, $t^{\prime\prime}$, $t^{\prime\prime\prime}$, and $t^{\prime\prime\prime\prime}$ denote four actions. The following conditions must be fulfilled:

\begin{enumerate}
    \item The strict inequality $v_j(x_{jt^{\prime}})-v_j(x_{jt^{\prime\prime}}) > v_j(x_{jt^{\prime\prime\prime}})-v_j(x_{jt^{\prime\prime\prime\prime}})$ holds, if and only if, the strength of the impact of $x_{jt^{\prime}}$ over $x_{jt^{\prime\prime}}$ is strictly higher than the strength of impact of $x_{jt^{\prime\prime\prime}}$ over $x_{jt^{\prime\prime\prime\prime}}$, on criterion $g_j$, , for $j=1,\ldots,n$.
    \item  The equality $v_j(x_{jt^{\prime}})-v_j(x_{jt^{\prime\prime}}) = v_j(x_{jt^{\prime\prime\prime}})-v_j(x_{jt^{\prime\prime\prime\prime}})$ holds, if and only if, the strength of impact of $x_{jt^{\prime}}$ over $x_{jt^{\prime\prime}}$ is the same to the strength of impact of $x_{jt^{\prime\prime\prime}}$ over $x_{jt^{\prime\prime\prime\prime}}$, on criterion $g_j$, for $j=1,\ldots,n$.
\end{enumerate}

In the construction of the value functions and the criteria weights we assume that the axioms of transitivity and independence hold (see \citealt{KeeneyRa1993}).

\subsection{Chromatic ordinal classification model}\label{sec:classification_model}
\noindent The chromatic ordinal classification model is an ordinal scale with categories and colours associated with them. Let $C = \{C_1,\ldots,C_r,\ldots,C_s\}$ denote a set of totally ordered (and pre-defined) categories, from the best $C_1$ (the lowest pandemic state impact), to the worst $C_s$ (the highest pandemic state impact): $C_1 \succ \cdots \succ C_r \succ \cdots \succ C_s$, where $\succ$ means ``impacts strictly more than''. The categories are used to define a set of states, as follows:

\begin{itemize}[label={--}]
    \item $C_1$ (green): Baseline state.
    \item $C_2$ (light green): Residual state.
    \item $C_3$ (yellow): Alarm state.
    \item $C_4$ (orange): Alert state.
    \item $C_5$ (red): Critical state.
    \item $C_6$ (dark red): Break state.
\end{itemize}

There are four fundamental states, from $C_2$ to $C_5$, with a particular set of associated measures/recommendations. It is worthy of note that, the colours assigned to each state change smoothly when reaching the boundaries of the neighbouring states and that they move quickly when passing from one state to the next in the upper part of the scale, as for example from $C_4$ to $C_5$, than when moving from a state to the next in the lower part of the scale, as for example from $C_2$ to $C_3$. It also goes quickly from top down, i.e., in a descending way. This can be done through the way the value functions are modelled and/or the choice of the values for setting the cut-off lines with the possible definition of thresholds (see subsection \ref{sec:ordinalscale}), for the justification.


\section{Modeling aspects}\label{sec:modeling}
\noindent This section provides the details of the three fundamental models used in our study: the criteria model, the aggregation model, and the graphical visualisation and communication model. The classification chromatic system and an illustrative example are also presented in this section.

\subsection{Criteria model}\label{sec:criteria}
\noindent A set of criteria built by the experts as the most relevant, taking into account the two main perspectives (called pillars) were used to characterise the pandemic, aimed to fulfil several desirable properties as stated in \cite{Keeney1992}: essential, controllable, complete, measurable, operational, decomposable, non-redundant, concise, and understandable. They were grouped as follows:

\begin{itemize}
    \item[A.] Pillar I ({\sc{ACT}}). \textit{Activity}.  This pillar was built to capture the main aspects of the {\sc{Covid-19}} registered or observed activity, i.e., the survival and development of the virus and its ability to still be active and cause infection in people in a given territorial unit. The following two {\sc{Covid-19}} activity criteria were considered to render this pillar operational.
    \begin{itemize}
    \item[1.] Criterion $g_1$ - \textit{Incidence} ({\sc{incid}}). The incidence  (see \citealt{Martcheva2015}) is the number of new {\sc{Covid-19}} positive cases presented daily, $N(\cdot)$, in the Official Health Reports. In most countries, the exact daily values vary periodically over each week. In Portugal, in particular, the evolution of new daily cases peaks markedly at day seven. Thus, to regularise the time series of the incidence, we consider the seven-day moving average and use this variable in our computations:
        \begin{equation}\label{eq:incidence}
        g_1(t) = x_{1t} = \frac{\displaystyle \sum_{u=t-6}^{t}N(u)}{7}
        \end{equation}
        We could use the raw data directly, but that choice would introduce artificial weekly fluctuation due to weak reporting at weekends. A longer periodic average, as for example, by considering the last fourteen days would lead to a slow effect of the impact. This analysis led us to consider a seven-day average as the most adequate for this criterion.
    \item[2.] Criterion $g_2$ - \textit{Transmission} ({\sc{trans}}). The transmission is modelled here as the rate of change in the active cases computed from the raw data of the daily incidence  $N(\cdot)$ (with no moving averages). With the goal of regularising these time series criterion values and smoothing the weekly fluctuations, we calculated the geometric mean over the last seven days. Our criterion is defined by the expression below:
        \begin{equation}\label{eq:transmissibiliy}
        g_2(t) = x_{2t} = \left(\prod_{u=t-6}^{t}\frac{\displaystyle \sum_{v=u-6}^{u}N(v)}{\displaystyle \sum_{v=u-7}^{u-1}N(v)}\right)^{\frac{1}{7}}
        \end{equation}
        With this formula, we have the advantage of a quicker response to the changes in incidence with respect to $R(t)$ transmission rate, the usual reproduction number of an epidemic with time. Moreover, our model has the same meaning as the $R(t)$, for $t=1$ (see \citealt{Koch2020}).
    \end{itemize}
    \item[B.] Pillar II  ({\sc{SEV}}) - \textit{Severity}. This pillar was built to capture the severity of the effects of {\sc{Covid-19}} on the Portuguese people, in particular on the health system. The following three {\sc{Covid-19}} severity criteria were considered to render this pillar operational.
    \begin{itemize}
    \item[3.] Criterion $g_3$ - \textit{Lethality} ({\sc{letha}}). The lethality is modelled here by taking into account the ratio of deaths at a given time period $u$ over the number of new cases in the fourteen days prior. Then, by considering the accumulate number of cases $N(\cdot)$ and number of deaths, $O(\cdot)$, it can be calculated by using the formula:
        \[
            \ell(t) =100\times\frac{O(t)- O(t-1)}{N(t-14) - N(t-15)},
        \]
        We hypothesise that the average time to death after the communication of the case is fourteen days. With the goal of regularising this variable and smoothing the fluctuations, we calculated a moving average of the last fourteen days for this criterion. The lethality formula used is given as follows:
        \begin{equation}
            g_4(t) = x_{3t} =\frac{1}{14}\sum\limits_{u=t-13}^{t}\ell(u) .
        \end{equation}
        Another formula could be defined for modeling the lethality, but this one has been considered the most adequate by the experts in our case. Lethality could be modelled using a seven-day moving average formula, however the 14 day moving average was more adequate given that the evolution of lethality is gradual and slow and the 14 days average regularises statistical fluctuations of the observed data.
    \item[4.] Criterion $g_4$ - \textit{Number of patients admitted to wards} ({\sc{wards}}).
        This criterion considers the total number of {\sc{Covid-19}} patients admitted to wards without counting those admitted to the intensive care unit, $H(\cdot)$, which is raw data. The formula is thus a direct one:
        \begin{equation}\label{eq:nurserie}
        g_4(t) = x_{4t} = H(t)
        \end{equation}
    \item[5.] Criterion $g_5$ - \textit{Number of patients admitted to wards} ({\sc{icu}}).
    Similar to the previous criterion, it counts the number of  {\sc{Covid-19}} patients admitted to the intensive care units, $U(\cdot)$, which is also raw data. The formula is also a direct one:
    \begin{equation}\label{eq:UCI}
        g_5(t) = x_{5t} = U(t)
    \end{equation}
    \end{itemize}
\end{itemize}

All the raw data $N(\cdot)$, $O(\cdot)$, $H(\cdot)$, and $U(\cdot)$, are available at the Direcção-Geral da Saúde (DGS) web site ({\tt www.dgs.pt}).

  \begin{remark}{(Fragility Point 1)}  \textit{Imperfect knowledge of criteria set (see \citealt{RoyEtAl2014})} This imperfect knowledge is mainly due to the imprecision of the tools and the procedures used to determine the raw data needed for the computation of the performance levels of the three criteria (namely $N(\cdot)$, since $O(\cdot)$, and $U(\cdot)$ do not suffer from significant imprecision) and also due to the arbitrariness of the formulas chosen for the three criteria ($g_1(\cdot)$, $g_2(\cdot)$, and $g_3(\cdot)$). Other models could have been selected and justified. Whenever a fragility point (weakness or vulnerability) is identified, sensitivity analyses are needed to guarantee the validity of the model and confidence in the results. These sensitivity analyses will be presented in Section \ref{sec:robustness}.
  \end{remark}

\subsection{About the rationale behind the set of criteria}\label{sec:rational_criteria}
\noindent The building this set of criteria followed a logic based on two fundamental principles that are somehow linked: familiarity with the problem and intelligibility. The first stems essentially from the information that the population received \textit{via} the media. All Portuguese citizens have knowledge, relatively informed, and with a certain degree of in-depth knowledge, of the impact of this pandemic, through two major concepts or aspects:, the activity of the pandemic and its severity. In defining the activity, two more further key concepts are included, which are also well known to the Portuguese people: incidence and transmission. As for the concept of severity, the concepts of mortality, occupancy of bed occupancy in the wards and in the ICU are published daily in the \textit{media} and are equally similarly familiar concepts to all Portuguese people. In total, there are five concepts that the Portuguese people deal with every day. These concepts came naturally to them, and it would be very difficult not to consider them in a set of criteria to assess the impact of the pandemic in the country. That was also the opinion of that was shared by the experts who are were part of the team that produced the indicator. We think that the justification the rationale for the importance of the consideration the familiarity principle was quite clear.

As for the principle of intelligibility, this is more closely related to the models used to operationalise each of the five concepts, then grouped into two pillars. Each of these models needed to be easy to understand for the population, when presented to them. Let us look at the case of incidence and the associated model. Incidence is defined as the number of positive cases detected daily (note that it is not possible to detect all the cases in a population and this number is not accurate). The direct use of this number in our indicator would imply greater fluctuation in values, since it also fluctuates throughout the week. We therefore chose to model incidence using the average number of cases \textit{per} week, a model that is clear and easy to understand by the population, thus confirming the principle of intelligibility. The same philosophy was applied in the construction of the models of the other four concepts.

In summary, the rationale for the set of criteria was made through concepts, two of a more global scope, the pillars of disease activity and severity, and five of a more local scope, which are the basis for building operational models of the two global pillars. Taking into account additional criteria, such as the number of tests, or others, would make the model more complex and we found in practice that these concepts did not have much influence on the perception of the impact on the pandemic, such as the ones we chose for our model

\subsection{Parameters of the aggregation model and the chromatic system}\label{sec:agregation}
\noindent This section presents the technical aspects related to the construction of the parameters of the additive aggregation model (i.e, the value functions and the weights), as well as the chromatic classification system.

The construction of the value functions (interval scales) and the weights of criteria (ratio scales) was performed using a simplified version of the Pairwise Comparison Deck of Cards Method (here called PaCo-DCM), proposed by \cite{Corrente2021}. This simplified version did not consider making all a pairwise comparison and did not take into account imprecise information leading to possible inconsistent judgements..

The origin of the DCM in MCDA dates back to the eighties, a procedure proposed by \cite{Simos1989}.  This method was later revised by \cite{FigueiraRo2002} and used for determining the weights of criteria in outranking methods. In this revised version, \cite{FigueiraRo2002} mention the possibility of using the method and SRF software, proposed in the same paper, to build not only ratio scales, in general, but also to build interval scales. For another extension DCM and a review of applications, see  \cite{SiskosTo2015}. Regarding the interval scales, a first attempt to build them was proposed by  \cite{PictetBo2008}, while \cite{BotteroEtAl2018} improved the DCM method to build more general interval scales (based on the definition of at least two reference levels with a precise meaning for policy/decision-makers, users, or experts). \cite{BotteroEtAl2018} also created another extension of the construction of ratio scales for determining the capacities of the Choquet integral aggregation method. \cite{DinisEtAl2021} made use of a tradeoff procedure for determining the weights of criteria for the additive MAVT model. The method used for computing the weights of criteria in this paper is very similar to the latter. Another recent and interesting extension of the DCM with visualisation tools was proposed by \cite{TsotsolasEtAl2019}.

\subsubsection{Value functions (interval scales)}\label{sec:intervalscale}
\noindent  The construction of the value functions using PaCo-DCM requires the use of pairwise comparison tables. This idea was introduced into MCDA by \cite{Saaty2016} and later adapted and improved to accommodate qualitative judgements by \cite{BanaEtAl2016}.

In what follows, we will show, step-by-step, the details of the application of PaCo-DCM. In sociotechnical processes, it is always important to provide key aspects of the context that enable and facilitate the evolution of these processes. It is important to note that we had a very limited time period for producing a first version of our tool (only three days) and for presenting a first prototype with meaningful results (a further ten days). This was possible as we benefitted from the help of a mathematician on the CCIST team, who also had strong expertise in {\sc{Covid-19}}  and was well very acquainted with the experts on the GCOM team. In addition, the CCIST member has substantial experience programming with Wolfram Mathematica\footnote{https://writings.stephenwolfram.com/2019/04/version-12-launches-today-big-jump-for-wolfram-language-and-mathematica/ }, which was crucial meaning we obtained the results of our tests and graphical tools almost instantaneously. This was a fundamental aspect for the interaction with the members of the GCOM team, comprised exclusively of physicians familiar with the fundamentals of mathematics.

We will present the main interactions between the CCIST team and the GCOM team for building together the value function of the first criterion as a sociotechnical process, which took into account the experts’ judgements and the technicalities of the PaCo-DCM tool. We will also present the details of all the computations. Readers can easily follow how we built an interval scale with a simple, but adequate version of the PaCo-DCM.

\begin{enumerate}
  \item \textit{The basics of PaCo-DCM for gathering and assessing the expert's judgements}. The method was introduced to the experts in a simplistic form, by explaining the meaning of the DCM used to assess their judgements \textit{via} a short example. The experts were provided with two sets of cards:
      \begin{enumerate}
        \item A very small set of labelled cards with very familiar objects (e.g., a lemon, an apple, and a mango) and that can easily be placed in preferential order (first mango, then apples, then lemon), from the best to the worst (for the sake of simplicity assume they are totally ordered, i.e., there are no ties). All the experts agreed on the same ranking.
        \item A large enough set of blank cards. These blank cards are used to model the intensity or strength of preference between pairs of objects.
        \item Assume we have three objects $o_i$, $o_k$, and $o_j$. If the experts feel the strength of preference difference between $o_i$ and $o_k$  is stronger than the strength of preferences difference between $o_k$ and $o_j$, they place more blank cards in between $o_i$ and $o_k$ than in between $o_k$ and $o_j$ (these are thus judgements for building a thermometer like scale). The experts can place as many cards as they want in between two objects and they do not need to count them, just hold them in their hands. In fact, in our case we used (wooden balls instead of blank cards; this does not invalidate the application of the method and it is more suited to gaining judgements from the experts, because the wooden balls are easy to handle and have a better visualisation effect. Experts may always revise their judgements about the strength of preference and change the number of cards in between two objects.
        \item We then explained the experts that:
            \begin{itemize}[label={--}]
              \item No blank card in between two objects does not mean that the two objects have the same value, but that the difference is minimal (minimal here means equivalent to the value of the unit, a concept the experts would subsequently understand better).
              \item One blank card means that the difference of preference is twice the unit.
              \item Two blank cards means that the difference of preference is three times the unit and so on.
            \end{itemize}
        \item Finally, we explained to the experts that in our case, we were modelling the strength or intensity of the pandemic impact instead of preferences, but the concept of preference was extremely useful to render the experts familiar with the concept of strength of impact and with the method.
      \end{enumerate}
  \item \textit{(At least) two well-defined reference levels}. PaCo-DCM requires the definition of two reference levels for the construction of an interval scale. These two levels must have a precise meaning for the experts. One level is, in general, located in the lower part of the scale and the other in the upper part of the scale. This is similar to the method proposed in \cite{BanaEtAl2016} where, in general, ``neutral'' and ``good'' reference levels are needed to build a scale, with the assignment of the values $0$ and $100$, respectively. In PaCo-DCM, the values of the reference levels do not need to be set at $0$ and $100$. Any two values can be used in PaCo-DCM for building the interval scale. In the application of the model to the pandemic situation, the two reference levels built from the interaction with the experts were the following:
      \begin{itemize}[label={--}]
        \item \textit{Baseline level}: Incidence value equal to $0$. This means that no new cases have been registered over the last seven days. It does not mean the absence of a pandemic, but the fact that no new cases have been observed. The value of the baseline impact level was set at $v_1(0) = 0$, which is an arbitrary origin on the interval scale for the $0$ preference level in the first criterion.
        \item \textit{Critical level}: Incidence value equal to $1125$. The value of the critical level was first set at $1100$, but after the discussion of the subsequent step, we made a slight adjustment to $1125$ and decided to set $v_1(1125) = 100$, which represents the highest value before entering a critical state. Due to the number of public health physicians and contact tracing (tracking), after $900$ new cases there is a saturation of resources, and the experts considered $1125$ to be an adequate number to model the critical level.
      \end{itemize}
  \item \textit{Setting the number of value function breakpoints}. In this step, we defined alongside the experts the most adequate way of discretising, by levels, the performances of the incidence, taking into account the initial two reference levels built in the previous step, $0$ and $1100$. A first discussion led us to consider only values in between $0$ and $2000$, more than this value would lead to an emergency state, even $2000$ seemed to be a very large value. Thus, we thought to discretise the range $[0, 2000]$ into six breakpoints $0$, $450$, $900$, $1350$, $1800$, and $2000$, but an width of $450$ in between two consecutive levels was considered quite large. Finally, we decided to discretise the range in ten points, with an width of $225$ between two consecutive points. The following values were finally considered, with an adjustment in the last one to be consistent with the width of $225$, and by considering the critical level at $1125$ instead of $1100$:

        \[
            \textcolor{blue}{\{0\}} \;\, \{225\} \;\; \{450\} \;\, \{675\} \;\, \{900\} \;\,  \textcolor{blue}{\{1125\}} \;\, \{1350\} \;\, \{1575\} \;\, \{1800\} \;\, \{2025\}.
        \]
  \item \textit{Inserting blank cards}. In this step, the experts were invited to insert blank cards in between consecutive levels; it corresponds to fill the diagonal of Table \ref{table:criterion_g1}. This process was performed for the initial number of breakpoints (ten) and led to several adjustments resulting from the sociotechnical constructive interaction process between experts (GCOM team) and the analysts (CCIST team). Each change was accompanied by a figure (see, for example, Figure \ref{fig:function_v1}), which was an important visualisation tool for assisting the experts. The consistency tests ones described in the next step were also performed. After building all the value functions (Appendix) and testing them with past observations of the pandemic, the two last levels below were discarded since a break level was set before reaching level $1575$.
        \[
        \textcolor{blue}{\{0\}} [0] \{225\} [2] \{450\} [4] \{675\} [6] \{900\} [8]  \textcolor{blue}{\{1125\}} [10] \{1350\} [13] \{1575\} \textcolor{red}{[11] \{1800\} [8] \{2025\}}
        \]
        However, we can observe that the number of blank cards increase until level $1575$, then it decreases. It means that the shape of the value function will move from a convex to a concave shape (it would be similar to a continuous sigmoid function). From a certain point, more new cases have almost the same impact on the pandemic as fewer new cases. We are referring to a part of the function where the situation would be out of control.
  \item \textit{Testing a more sophisticated version of the method}. In the PaCo-DCM method, more cells of Table \ref{table:criterion_g1} can be filled from the judgements provided by the experts. The time limitation and the good understanding of the method by themselves reduced the number of interactions for checking consistency judgements. Table \ref{table:criterion_g1} contains thus all the possible comparisons with the eight levels kept (here we are not considering the two last ones). An example of a test for assessing the impact differences between the two non-consecutive levels $\{900\}$ and $\{225\}$ was performed as follows. We placed the set of six blank cards in between $\{900\}$ and $\{675\}$, and the set of seven cards in between $\{675\}$ and $\{225\}$ (the experts do not really need to know how many cards are in between these levels, but these two numbers came from the previous interaction and were provided by the experts). Then, we placed a set of sixteen cards in between $\{900\}$ and $\{225\}$ and we asked the experts to compare the three sets of cards, asking whether they felt comfortable with the third set of sixteen cards. If not, we started by removing the blank cards, one by one. We removed blank cards, until there was a set of thirteen blank cards. Then, we told the experts that thirteen cards is slightly inconsistent and showed them why. We finally, asked them whether they felt comfortable with a set of fourteen blank cards (six, in between $\{675\}$ and $\{900\}$, plus seven, in between $\{225\}$ and $\{675\}$, plus one), and they agreed. The other cells of the table can be filled by transitivity, i.e., by following the consistency condition presented in \cite{Corrente2021}. The impact difference between two non-consecutive cells is determined as follows:
      \begin{equation}
      \label{eq:2}
        e_{ij} = e_{ik} + e_{kj} + 1 \;\;\, \mbox{for all} \;\; i,k,j = 1,\ldots,t \quad \textrm{and} \quad i < k < j
      \end{equation}
      We can see that $e_{\{900\},\{225\}} = e_{\{900\},\{675\}} + e_{\{675\},\{225\}} + 1 = 6+7+1=14$.
    \begin{table}[htbh]
    \small
    \centering
    \begin{center}
    \begin{tabular}{c|cccccccc|}
        \cline{2-9}
        & \textcolor{red}{$\,\;\;\{0\}\,\;\;$} & $\;\{225\}\;$ & $\;\{450\}\;$ & $\;\{675\}\;$ & $\;\{900\}\;$ & \textcolor{red}{$\{1125\}$} & $\{1350\}$ & $\{1575\}$ \\
        \hline
        \multicolumn{1}{|l|}{$\{0\}$}    & \cellcolor{gray} & \textbf{0}   & 3 & 8 & 15 & 24 & 35 & 49 \\
        \multicolumn{1}{|l|}{$\{225\}$}  & & \cellcolor{gray} & \textbf{2}     & 7 & 14 & 23 & 34 & 48 \\
        \multicolumn{1}{|l|}{$\{450\}$}  & & & \cellcolor{gray} & \textbf{4}     & 11 & 20 & 31 & 45 \\
        \multicolumn{1}{|l|}{$\{675\}$}  & & & & \cellcolor{gray} & \textbf{6}     & 15 & 26 & 40 \\
        \multicolumn{1}{|l|}{$\{900\}$}  & & & & & \cellcolor{gray} & \textbf{8}     & 19 & 33 \\
        \multicolumn{1}{|l|}{$\{1125\}$} & & & & & &  \cellcolor{gray} & \textbf{10}   & 24 \\
        \multicolumn{1}{|l|}{$\{1350\}$}  & & & & & & & \cellcolor{gray} & \textbf{13} \\
        \multicolumn{1}{|l|}{$\{1350\}$}  & & & & & & & & \cellcolor{gray}             \\
        \hline
    \end{tabular}
    \caption{Pairwise comparison table for criterion $g_1$ (incidence)}
    \label{table:criterion_g1}
    \end{center}
    \end{table}
  \item \textit{Computations}. The computation of the values of the breakpoints was done as follows:
    \begin{itemize}[label={--}]
      \item The values of the reference levels: $v_1(0) = 0$ and $v_1(1125) = 100$.
      \item The number of units in between them, $h = (0+1)+(2+1)+(4+1)+ (6+1) + (8+1) = 25$. Remember that $0$ cards does not mean the same value, but that the difference is equal to the unit. Thus, we need to add one more to all the number of blank cards in between two levels.
      \item The value of the unit, $\alpha = \big(v_1(1125)-v_1(0)\big)/h = (100-0)/25 = 4$. The value of the unit is equal to four points, and now the experts were able to understand better the concepts of unit and value of the unit.
      \item The values of the breakpoints are now easy to determine: $v_1(225) = 0 + 4\times 1 = 4$, $v_1(450) = 0 + 4 \times 4 = 16$, and so on, for the remaining:  $v_1(675) = 36$, $v_1(900) = 64$, $v_1(1125) = 100$, $v_1(1350) = 140$, $v_1(1575) = 200$.
    \end{itemize}
  \item \textit{The shape of the value function}. After assessing the values of the breakpoints, we can then draw a piecewise linear function as in Figure \ref{fig:function_v1}. Any value within each linear piece can then be obtained by linear interpolation.


    \begin{center}
    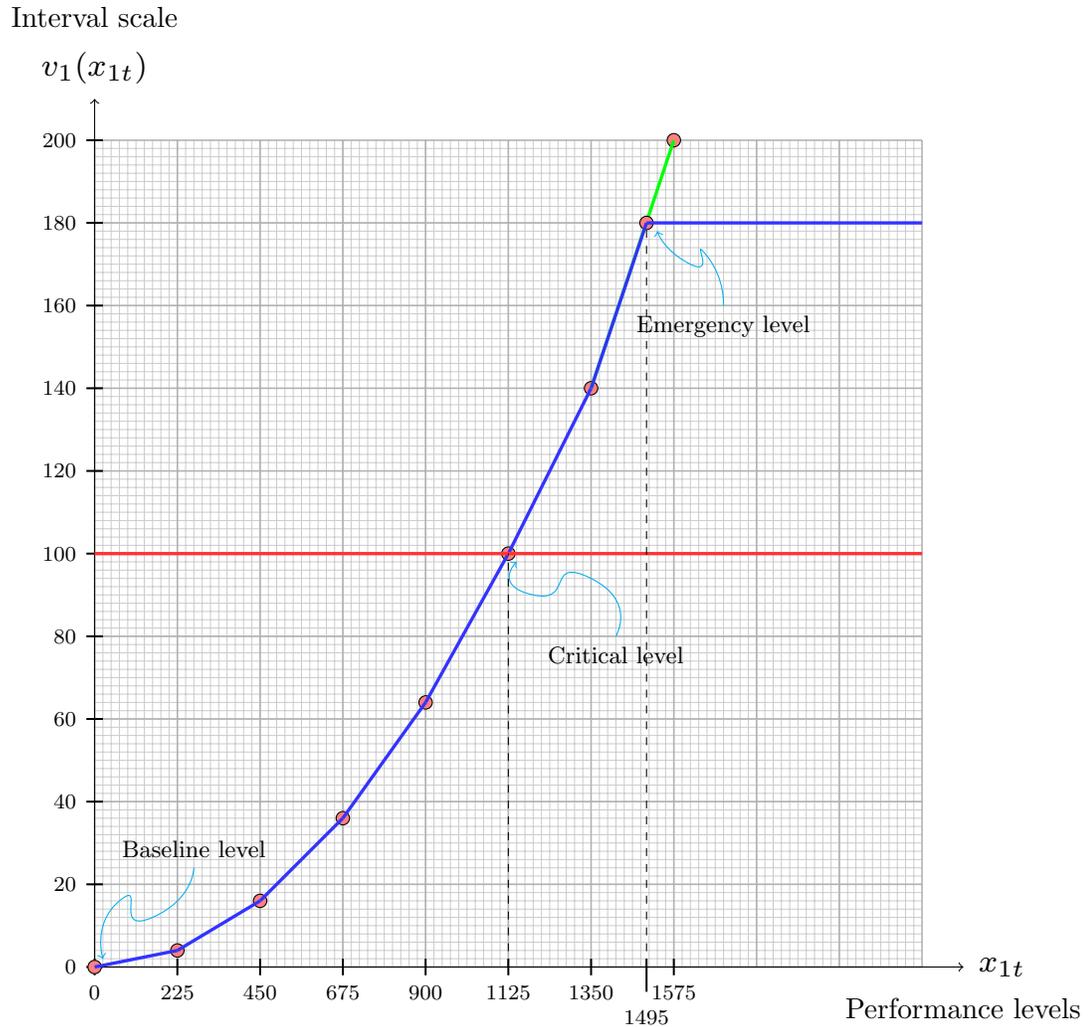
\begin{figure}[htbp]
    \centering
    \begin{tikzpicture}[scale=1.1]
        \draw[help lines,gray!40,step=.1] (0,0) grid (10,10);
        \draw[help lines,gray!60,line width=.6pt,step=1] (0,0) grid (10,10);
        \tikzstyle{axes}=[]
        \tikzstyle{important line}=[very thick]
        \tikzstyle{information text}=[rounded corners,inner sep=1ex]
        \startscope[style=axes]
            \draw[->] (0,0) -- (10.5,0) node[scale=1.5,right]{\scriptsize $x_{1t}$} coordinate(x axis);
            \node[below] at (10.5,-0.25){Performance levels};
            \draw[->] (0,0) -- (0,10.5) node[scale=1.5,above]{\scriptsize $v_1(x_{1t})$ } coordinate(y axis);
            \node[above] at (0,11.25){Interval scale};
             \draw[-,line width=0.75pt,black](0.1,0) -- (-0.1,0);
            \node[left] at (-0.1,0) {\scriptsize {$0$}};
             \draw[-,line width=0.75pt,black](0.1,1) -- (-0.1,1);
            \node[left] at (-0.1,1) {\scriptsize {$20$}};
             \draw[-,line width=0.75pt,black](0.1,2) -- (-0.1,2);
            \node[left] at (-0.1,2) {\scriptsize {$40$}};
             \draw[-,line width=0.75pt,black](0.1,3) -- (-0.1,3);
            \node[left] at (-0.1,3) {\scriptsize {$60$}};
             \draw[-,line width=0.75pt,black](0.1,4) -- (-0.1,4);
            \node[left] at (-0.1,4) {\scriptsize {$80$}};
             \draw[-,line width=0.75pt,black](0.1,5) -- (-0.1,5);
            \node[left] at (-0.1,5) {\scriptsize {$100$}};
             \draw[-,line width=0.75pt,black](0.1,6) -- (-0.1,6);
            \node[left] at (-0.1,6) {\scriptsize {$120$}};
             \draw[-,line width=0.75pt,black](0.1,7) -- (-0.1,7);
            \node[left] at (-0.1,7) {\scriptsize {$140$}};
             \draw[-,line width=0.75pt,black](0.1,8) -- (-0.1,8);
            \node[left] at (-0.1,8) {\scriptsize {$160$}};
             \draw[-,line width=0.75pt,black](0.1,9) -- (-0.1,9);
            \node[left] at (-0.1,9) {\scriptsize {$180$}};
             \draw[-,line width=0.75pt,black](0.1,10) -- (-0.1,10);
            \node[left] at (-0.1,10) {\scriptsize {$200$}};
            \draw[-,line width=0.75pt,black](0,0.1) -- (0,-0.1);
            \node[below] at (0,-0.1) {\scriptsize {$0$}};
            \draw[-,line width=0.75pt,black](1,0.1) -- (1,-0.1);
            \node[below] at (1,-0.1) {\scriptsize {$225$}};
            \draw[-,line width=0.75pt,black](2,0.1) -- (2,-0.1);
            \node[below] at (2,-0.1) {\scriptsize {$450$}};
            \draw[-,line width=0.75pt,black](3,0.1) -- (3,-0.1);
            \node[below] at (3,-0.1) {\scriptsize {$675$}};
            \draw[-,line width=0.75pt,black](4,0.1) -- (4,-0.1);
            \node[below] at (4,-0.1) {\scriptsize {$900$}};
            \draw[-,line width=0.75pt,black](5,0.1) -- (5,-0.1);
            \node[below] at (5,-0.1) {\scriptsize {$1125$}};
            \draw[-,line width=0.75pt,black](6,0.1) -- (6,-0.1);
            \node[below] at (6,-0.1) {\scriptsize {$1350$}};
            \draw[-,line width=0.75pt,black](7,0.1) -- (7,-0.1);
            \node[below] at (7,-0.1) {\scriptsize {$1575$}};
            \node[draw,scale=0.4,circle,fill=red!50,text=blue] at (0,0) {$\;$};
            \node[draw,scale=0.4,circle,fill=red!50,text=blue] at (1,0.2) {$\;$};
            \node[draw,scale=0.4,circle,fill=red!50,text=blue] at (2,0.8) {$\;$};
            \node[draw,scale=0.4,circle,fill=red!50,text=blue] at (3,1.8) {$\;$};
            \node[draw,scale=0.4,circle,fill=red!50,text=blue] at (4,3.2) {$\;$};
            \node[draw,scale=0.4,circle,fill=red!50,text=blue] at (5,5) {$\;$};
            \node[draw,scale=0.4,circle,fill=red!50,text=blue] at (6,7) {$\;$};
            \node[draw,scale=0.4,circle,fill=red!50,text=blue] at (7,10) {$\;$};
            \draw[-,line width=1.25pt,blue!80](0,0) -- (1,0.2);
            \draw[-,line width=1.25pt,blue!80](1,0.2) -- (2,0.8);
            \draw[-,line width=1.25pt,blue!80](2,0.8) -- (3,1.8);
            \draw[-,line width=1.25pt,blue!80](3,1.8) -- (4,3.2);
            \draw[-,line width=1.25pt,blue!80](4,3.2) -- (5,5);
            \draw[-,line width=1.25pt,blue!80](5,5) -- (6,7);
            \draw[-,line width=1.25pt,green](6,7) -- (7,10);
            \draw[-,line width=0.75pt,black](6.67,0.1) -- (6.67,-0.3);
            \node[below] at (6.67,-0.4) {\scriptsize {$1495$}};
            \node[draw,scale=0.4,circle,fill=red!50,text=blue] at (6.67,9) {$\;$};
            \draw[-,line width=1.25pt,blue!80](6,7) -- (6.67,9);
            \draw[-,line width=1.25pt,blue!80](6.67,9) -- (10,9);
            \draw[-,line width=1.25pt,red!80](0,5) -- (10,5);
            \draw[dashed](5,0) -- (5,5);
            \draw[dashed](6.67,0) -- (6.67,9);
            \draw [->, cyan] plot [smooth, tension=2] coordinates {(1.2,1.2) (0.7,0.6) (0.3,0.8) (0.1,0.1)};
            \node[above] at (1.2,1.2) {\footnotesize Baseline level};
            \draw [->, cyan] plot [smooth, tension=2] coordinates {(6.3,4) (6,4.7) (5.3,4.5)(5.1,4.9)};
            \node[below] at (6.3,4) {\footnotesize Critical level};
            \draw [->, cyan] plot [smooth, tension=2] coordinates {(7.6,8) (7.4,8.6) (7.2,8.5) (6.8,8.9)};
            \node[below] at (7.6,8) {\footnotesize Emergency level};
        \stopscope
    \end{tikzpicture}
    \caption{Shape of the value function for incidence}
    \label{fig:function_v1}
    \end{figure}
    \end{center}
     A marginal change in the lower part of the scale has less impact on the pandemic than the same marginal change in the upper part of the scale, as depicted in Figure \ref{fig:function_v1}. Before reaching the emergency level, this is a convex value function, and we can observe the marginal increase in the value of the impact. Moving from $0$ to $225$ implies an increase in the indicator value, from $0$ to $4$, while moving from $1125$ to $1350$ implies an increase in the impact value from $100$ to $140$. The same number of additional units $225$ in the upper part of the scale produces a much higher impact ($40$) than in the lowest part of the scale (only $4$). It increases closer to the critical level than to the baseline level. This was a strong requirement established by the experts.
    \item \textit{Remark}. Please note that the values of the break/saturation levels were calculated considering Portugal's capacity and, from a certain limit, the normal and extraordinary capacity was exceeded, and a rupture state was reached. All these values can, and should be, adjusted according to the characteristics of each country and can be revised according to the adaptations of each country. When we performed the PACI, the values for Portugal were the ones used and there was no need for further readjustment. There are feedback mechanisms from the ground to adjust the levels and the PACI itself, by monitoring the situation. It contributes to the allocation of resources and adjustment of the limits. In Portugal there was no such need, till this moment.
  \item \textit{The output of the model and possible approximations}. One of the outputs of the model is a piecewise linear function, whose mathematical expression of which can be stated as follows:
  \begin{equation}\label{eq:function_x1}
    v_1(x_{1t}) = \left\{
    \begin{array}{lcl}
      4x_{1t}/225      & \mbox{if} & x_{1t} \in [0,\,225[ \\
      4x_{1t}/75 -8    & \mbox{if} & x_{1t} \in [225,\,450[ \\
      4x_{1t}/45 -24   & \mbox{if} & x_{1t} \in [450,\,675[ \\
      28x_{1t}/225 -48 & \mbox{if} & x_{1t} \in [675,\,900[ \\
      4x_{1t}/25 - 80  & \mbox{if} & x_{1t} \in [900,\,1125[ \\
      44x_{1t}/225-120 & \mbox{if} & x_{1t} \in [1125,\,1350[ \\
      52x_{1t}/225-168 & \mbox{if} & x_{1t} \in [1350,\,1495[ \\
               180  & \mbox{if} & x_{1t} \in [1495,\,+\infty[ \\
    \end{array}
    \right.
    \end{equation}
    This particular function could be approximated by a quadratic function without losing much information, but such an approximation needs to be validated by the experts:
    \[
    \widetilde{v_{1}}(x_{1t}) =
        \left\{
        \begin{array}{ccc}
            100(x/1125)^{2} & \text{if} & x_{1t} \in [0,\,1125] \\
            & & \\
            180 & \text{if} & x_{1t} \in [1125,\,+\infty[.
        \end{array}
        \right.
    \]
    The mean of the Euclidean distance between the functions $\widetilde{g_{1}}$ and $g_{1}$ is given by the following expression:
    \[
        \frac{\sqrt{\int\limits_{0}^{+\infty }\left( g_{1}\left( \xi \right) -
        \widetilde{g_{1}}\left( \xi \right) \right) ^{2}d\xi }}{\sqrt{%
        \int\limits_{0}^{+\infty }\left( g_{1}\left( \xi \right) \right) ^{2}d\xi }}
        =6.40389\times 10^{-4},
    \]
    which is almost negligible and shows that the approximation does not lead to the loss of much information.
  \item \textit{Missing, imprecise, and inconsistent judgments.} The method described in \cite{Corrente2021} also allows us to deal with missing, imprecise, and inconsistent judgements. The inconsistency analysis is performed by using linear programming, similarly to that performed in other MCDA tools, as for example, in \cite{MousseauEtAl2002}. The team members’ experience and the way the sociotechnical interaction was conducted largely facilitated the information gathering process, not requiring the use of more sophisticated functionalities of the PaCo-DCM tool. Time pressures meant it was not possible to use all of method's functionalities, including the fact that posing more complex questions leads to possible inconsistent judgements, but since the experts validated the results (value function and weights, in this case), there was no real need to render the dialogue more complex. This is not a question of corrupting the process of the method application. In a sociotechnical co-constructive process, when constructing the value functions and the weights, there are no true values for such parameters because there is no true reality, i.e., there are no true value functions and true weights. We used the questions that were the most adequate and accepted by the experts given the time constraints. With no such a constraint, we could pose more questions, and render the process more complex, but since this process had been validated, there were no need to introduce additional questions. In the future we could introduce more complexity, but again it does not mean corrupting a process. This is a sociotechnical approach with some limitations inherent to all sociotechnical processes because our system is a mental construct and our model is another mental construct too.
  \item \textit{The break level}. After running the model for the whole set of days during the pandemic, the experts realised they could set a maximum of $180$ points for this function since all the situations beyond such a point would be equally bad and out of control. The function was thus truncated at the level $1495$, which is the first level with a value of $180$. After this performance level the situation collapses and all the performance levels are felt as serious as the break level.
\end{enumerate}

The piecewise functions for the other four criteria, as well as the number of blank cards in between consecutive levels, are provided in the Appendix.

\begin{remark}{(Fragility Point 2)} \textit{Subjectivity in building the value functions}. There is some obvious subjectivity in the construction of the value functions since the experts are not precise instruments like high tech thermometers. In addition, for example, there is no true value function for modelling the incidence; this function is a construct, which can be more or less adequate to the situation. This is another kind of fragility point in our model, which justifies the use of sensitivity analyses as we will present in Section \ref{sec:robustness}.
\end{remark}

\subsubsection{Weighting coefficients (ratio scales)}\label{sec:ratioscale}
\noindent The assignment of a value for each criterion weights was also performed through PaCo-DCM, but the interaction protocol with the experts and the nature of the judgements were presented in different way. The weights are interpreted here as scaling factors or substitution rates. The dialogue with the experts was conducted as follows:

\begin{enumerate}
  \item \textit{Constructing dummy situations}. A set of five dummy situations (or statuses) one \textit{per} criterion, representing the swings between the baseline level and the critical level were built as follows (also see \citealt{DinisEtAl2021}).
      \begin{itemize}[label={--}]
        \item $p_1 = (1125, \, 0, \, 0, \, 0, \, 0) \equiv (100, \, 0, \, 0, \, 0, \, 0)$. This situation represents the impact on the pandemic of the swing (regarding the first criterion) from the baseline level to the critical level, maintaining the remaining criteria at their baseline levels.
        \item $p_2 = (0, \, 1, \, 0, \, 0, \, 0) \equiv (0, \, 100, \, 0, \, 0, \, 0)$. The meaning of this situation is similar to the one provided in the first situations.  A transmission rate equal to one is considered adequate by the experts to represent the critical level.
        \item $p_3 = (0, \, 0, \, 3.6, \, 0, \, 0) \equiv (0, \, 0, \, 100, \, 0, \, 0)$.  The meaning of this situation is similar to the one provided for the first situation. As a definition of this situation, the experts considered that half of the maximum value of $g_3(t)$ along the pandemic in Portugal corresponds to $100$ points. Thus, $\max \left\{g_3(t)\right\} = 7.19148$. Consequently, $100$ points correspond to the value of lethality of $3.59574 \approx 3.6$.
        \item $p_4 = (0, \, 0, \, 0, \, 2500, \, 0) \equiv (0, \, 0, \, 0, \, 100, \, 0)$.  The meaning of this situation is similar to the one provided in the first situation. The $2500$ represent $15\%$ of the total number of beds, which is  an adequate number for defining the critical level.
        \item $p_5 = (0, \, 0, \, 0, \, 0, \, 200) \equiv (0, \, 0, \, 0, \, 0, \, 200)$.  The meaning of this situation is similar to the one provided for the first situation. The $200$ beds represent $80\%$ of the difference between the current number of beds and the number of beds in existence before the pandemic, this number was defined by the experts as adequate to represent the critical level.
      \end{itemize}
    The concept of swings is in line with the swing weighting technique by \cite{VonWinterfeldtEd1986} and the use of two reference levels with the concepts of ``neutral'' and ``good'' by \cite{BanaEtAl2016}.
  \item \textit{Ranking the dummy situations with possible ties}. The experts received five cards, one with each one of the previous situations and the analyst team asked them to provide a ranking of these five cards, with possible ties, according to the impact that the swings have on the pandemic. The situation(s) leading to the highest impact was (were) placed in first position, the one(s) with the second greatest impact on the second, and so on. The following ranking was proposed by the experts.
      \[
        \{p_1\}\;\, \{p_3,\,p_4,\,p_5\}\;\,\{p_2\}\;
      \]
      The analysts explained to the experts that the situation in the first situation will receive the highest weight, the ones in the second position the second highest weight, and the situation in the last position the lowest weight.
  \item \textit{Inserting blank cards}. The experts were invited to insert blank cards in between consecutive positions to differentiate the role each weight (swing) would have on the impact of pandemic, after telling them the meaning of swings and substitution rates. The following set of blank cards (in between brackets) was provided by the experts.
      \[
        \{p_1\}\;[2]\; \{p_3,\,p_4,\,p_5\}\;[3]\;\{p_2\}\;
      \]
      Similar to the value functions, a more sophisticated PaCo-DCM procedure could be used for such a purpose, but the experts felt comfortable with the information they provided.
  \item \textit{Assessing the value of the substitution rates}. This was the most difficult question for the experts. We need to establish a relation between the weight of the criterion in the first position of the ranking (incidence) and the weight of the criterion in the last position of the ranking (transmission). In PaCo-DCM, this is called the $z-$ratio, used to build a ratio scale. After a long discussion and several attempts, the experts provided the following relation between the two weights: $z=\hat{w}_1/\hat{w}_2=2$. We are using  $\hat{w_j}$, for the non-normalised weights of criterion $g_j$, for $j=1,\ldots,5$.
  \item \textit{Calculations}. The computations are similar to the ones performed for the value functions:
      \begin{itemize}[label={--}]
      \item The values of the non-normalised weights of the situations in the first and last positions of the ranking, i.e., $\hat{w}_1=2$ and $\hat{w}_1=1$.
      \item The number of units in between them, i.e., $h = (2+1)+(3+1) = 7$.
      \item The value of the unit, $\alpha = \big(\hat{w}_1 - \hat{w}_2\big)/h = (2-1)/7 = 0.14286$.
      \item The non-normalised weights: $\hat{w}_2 = 1$, $\hat{w}_3 = \hat{w}_4 = \hat{w}_5 = 1.42858$, and $\hat{w}_1 = 2$.
      \item The normalised weights: $w_2 = 1/7.28574=0.13725$, $w_3 = w_4 = w_5 = 1.42858/7.28574 = 0.19608$, and $w_1 = 2/7.28574 = 0.27451$.
    \end{itemize}
  \item \textit{Final adjustments}. After adjusting the model results to the real pandemic data and some discussions with the experts, the following weights were proposed for this model: $w_2=0.141$, $w_3 = w_4 = w_5 = 0.193$, and $w_1 = 0.280$.
\end{enumerate}

\begin{remark}{(Fragility Point 3)} \textit{Subjectivity in building the weights of criteria}.
    The justification is in line with the one provided in Remark 2, which also requires the use of sensitivity analyses (see Section \ref{sec:robustness}).
\end{remark}

\subsection{The possible existence of dependence between criteria}\label{sec:dependence}
\noindent The possible dependence between some criteria could be a major concern and was subject to reflection before taking the decision to propose the current criteria model (with the five criteria within two pillars) as well as the aggregation preference model; an additive model with measurable multi-criteria value functions, as proposed in \cite{DyerSa1979}. This topic is on slippery terrain since the path or the theory can lead to nonsense or useless models. We needed to ensure we did not make serious mistakes and reached a good compromise. We prepared some questions that we thought were legitimate: What is the real meaning of the dependence between or among criteria? Are the dependencies clear enough to be understood by the actors involved in the process and more importantly by the general population? If they really exist, how could we highlight and assess them? We will see that there are many forms and definitions of dependence between criteria. Also, intelligibility is not always easy to guarantee, and some of the identification and assessment procedures are quite time consuming in a context with multiple experts and many other actors, who may have different opinions about these dependencies

There are two major categories of dependencies. One is related to the (direct or indirect) factors we consider in the definition of the criteria, or that influence the output of the criteria in some way (i.e., the performance levels); this is called \textit{structural dependence} \citep[see,][]{Roy1996}. There is also possible technical dependence between the performance levels of the criteria – \textit{statistical dependence} – which can only be identified when sufficient data are available. Structural and statistical dependencies can be present simultaneously. The most important feature of these dependencies is that they do not depend on the experts’ impact judgements; they are ‘objective’ in the sense they are consensual for all the actors involved in the process. A different category of dependence requires the intervention of the experts’ judgements; it is subjective, and also related to the aggregation model. Here we refer to this type of dependence as \textit{subjective dependence}. In what follows we will provide some details about the two categories of dependence in the context of our tool.

\begin{enumerate}
    \item \textit{Factor structural and statistical dependencies.} These are related to the links between the factors that contribute to the definition of the concept of each criterion and/or the relation between the data on these factors. In our case these relations can exist in the formation of the criteria levels from the raw data.
        \begin{enumerate}
            \item \textit{Structural dependencies}. In activity pillar, the criteria transmission and incidence both share the same common factors contributing to their evolution: the transmission rate ($g_1$) depends on the active cases (infectious people), contacts (with people and/or contaminated spaces, surfaces, or objects), and the characteristics of the virus variants, while incidence ($g_2$) also depends on active cases, contacts, and on the transmission rate itself. This shows that there is a clear structural link between transmission and incidence. A structural link can also be observed in the formula of the criteria model of $g_2$: the numerator inside the product is equal to $7$ in $g_1$. Does it mean that there is also a statistical link? Is it enough to conclude about double or triple counting? Were these two criteria reduced to a single one, would it be beneficial and intelligible for the actors and population? The answer of the statistical link will be provided in the next paragraph, while in the comments below, we will try to answer the remaining three questions raised.
            \item \textit{Statistical dependencies}. For this dependence we compute the correlation between all pairs of criteria as presented in the tables below at two distinct moments: several days before the press conference and last December. \\
                \begin{figure}[h!]
                \caption{Correlation between criteria (10/7/2021)}
                \begin{center}
                \begin{tabular}{|c|c|c|c|c|}\hline
                        & $g_2$ & $g_3$ & $g_4$ & $g_5$ \\ \hline
                  $g_1$ &-0.034 & 0.178 & \textbf{0.919} &\textbf{0.836}  \\ \hline
                  $g_2$ &       &-0.082 &-0.214 &-0.280 \\ \hline
                  $g_3$ &       &       & 0.214 &0.163  \\ \hline
                  $g_4$ &       &       &       &\textbf{0.967}  \\\hline
                \end{tabular}
                \end{center}
                \end{figure}
                \begin{figure}[h!]
                \caption{Correlation between criteria (21/12/2021)}
                \begin{center}
                \begin{tabular}{|c|c|c|c|c|}\hline
                        & $g_2$ & $g_3$ & $g_4$ & $g_5$ \\ \hline
                  $g_1$ &-0.017 & 0.147 & \textbf{0.877} & \textbf{0.801} \\ \hline
                  $g_2$ &       &-0.043 &-0.181 &-0.249 \\ \hline
                  $g_3$ &       &       &-0.249 & 0.233 \\ \hline
                  $g_4$ &       &       &       & \textbf{0.967} \\\hline
                \end{tabular}
                \end{center}
                \end{figure}
            \item \textit{Comments}. The double counting argument cannot be used to justify replacing two criteria with one with different units and a dummy meaning which is difficult for people to understand \citep[see, for example,][for a more detailed explanation on this aspect]{Roy1996}. There are some technical aspects regarding the dynamical systems the question of our concrete models, which can be highlighted:
                \begin{itemize}[label={--}]
                    \item We consider the daily rate of change, regarding the system as a discrete dynamical system, or its derivative, which regards the system as a continuous dynamical system, or the basic reproductive number, regarding the system from the perspective of epidemiology. All these criteria are topologically conjugated, i.e., basically represent the same observation with different approaches. We preferred to use the daily rate of change, since it is very simple to compute and understand. In dynamical systems theory, the momentum and the position are not dependent criteria. The daily change drives the incidence, but depends on the behaviour of the population and on the rate of virus transmission; incidence depends on the rate of change and on the incidence of day minus one. The first and second criteria are not even dependent in the mathematical sense. They are structurally related, naturally, but not dependent in a subjective sense (see point 2).
                    \item Lethality depends on the severity of the virus, the average immune response of the population, the rate of vaccination, the efficiency of that vaccination and the efficiency of the therapeutics at home, on the number of patients admitted to wards and, finally, on the number of patients admitted to ICU. It varies with time and has changed dramatically with different variants of the virus and with the evolution of therapeutic measures. If we consider the variables in short periods of time, we believe that some dependence  is misleading.
                    \item There was also a decoupling also between patients in number of patients admitted to wards and number of patients admitted to ICU with with over time. If we want to measure the evolution of the the severity of the epidemic’s severity, we must consider the two criteria mut be considered. The PACI is a measure of the severity combined with transmissibility and incidence. The time evolution of the pandemic in Portugal had proven that all the five criteria were equality relevant.
                    \item The question is quite interesting if we were to consider the evaluation of an epidemic that has no changes in therapeutics, no mutations, and constant immunity. In that case, an indicator could be built using only incidence and transmissibility.
                    \item To conclude a statistical link does not mean an intelligible and accepted dependence. This is thus not enough to conclude about a double. triple, or multiple counting.
                \end{itemize}
        \end{enumerate}
    \item \textit{Subjective dependencies}. The two conditions below are generally imposed to guarantee the existence of an additive measurable multi-criteria value function (for three or more criteria).
        \begin{enumerate}
            \item \textit{Mutual preference (impact) independence}. The criteria are mutually preference (or impact) independent when all the proper sets of the set of criteria are preference (or impact), independent of their complements. If we were to consider two alternatives (or, time periods) with the same performance levels on the criteria belonging to the complement set; there is independence between the two alternatives if the preference between both does not depend on the criteria in the complement set. This condition guarantees a function that provides an ordinal order, the alternatives or time periods (in our settings).
            \item \textit{Mutual difference independence}.  The criteria are mutually difference independent when all the proper sets of the set of criteria, are difference independent of their complements, where a difference independence means that a difference between two alternatives characterised on several criteria and differing only on one, does not depend on the performance levels of the other criteria \citep{Dyer2016}. This condition guarantees the additivity of the aggregation model and captures the strength of preference or impact (in our settings).
            \item \textit{Comments}. It should be noted that the construction of measurable value functions requires the check not only for the two previous items, but also a condition called difference consistency, as well as some more technical assumptions \citep[see][]{Dyer2016, DyerSa1979}).
            \item \textit{Mutual preference (impact) independence test}. Given the time constraints, we did not conduct any independence test with the experts. Even if we were able to make it possible, there was no guarantee that all the experts would agree on the existence of the two previous types of dependence. However, and since one of the analysts was also an expert, it was legitimate to carry out some tests with him to check mutual preference (impact) independence, between the criteria of the activity pillar against the severity pillar. The indifference relation (denoted by the symbol ``$\sim$'') was used for doing such tests. We proceeded as follows with some pairs of time periods: considering a time period performances (\textbf{0.962}, 3, 500, 100) and a second time period performances (\textbf{1.000}, 3, 500, 100), we asked the analyst if these two time periods were indifferent in terms of the impact, i.e., if they produce the same impact. The answer was ``Yes'', so we made a change to the performance levels of the three criteria of pillar two and asked the same question again. The answered was again positive, i.e.,  (\textbf{0.962}, 6, 1000, 200) $\sim$ (\textbf{1.000}, 6, 1000, 200). Please note that a negative answer would imply that the two first criteria would not be independent of the last three). We proceeded the same way with different pairs of criteria and there was no hesitation, i.e., any doubt about the violation of mutual preference (impact) independence.
            \item \textit{Comments}. This test and the assumption of independence seemed to be an acceptable working hypothesis to build our model. The results of our model also confirmed our choice was adequate.
        \end{enumerate}
\end{enumerate}

The set of criteria was conceived in such a way it can be used to construct, transform, make, evolve and justify the impact judgements of the experts based on concrete elements. It means that the five criteria are related to concrete elements that the experts understand as well as the general population. They are thus related to particular significance axes. If some of these criteria were replaced by a dummy or abstract one, there would be a strong risk of loss of the meaning of the initial criteria \citep{Roy1996}. This is the reason why we kept this model as simple as possible, since we validated that it was an adequate representation of the problem.

\subsection{Illustrative example}\label{sec:example}
\noindent This is an illustrative example with five actions, i.e., five different time points in the pandemic. Moment $t=0$ is four days before the press conference with the \textit{media} (14 July 2021) at the PMA (in Lisbon). The other moments, $t$, were set with respect to the number of the days before $t=0$ and corresponds to the first lock down in Portugal (20 March 2020), one of the lowest activity and severity periods (31 July 2020), Christmas (24 December 2020), and some days after the second lock down (24 January 2021). Table \ref{tab:performance} presents the activity as severity performance levels for the five considered criteria, according to the two pillars.

\begin{center}
\begin{table}[htbp]
\centering
\begin{tabular}{clrrrrr}
\toprule
& & \multicolumn{2}{c}{Pillar I ({\sc{ACT}})} & \multicolumn{3}{c}{Pillar II ({\sc{SEV}})} \\
\cmidrule(r){3-4}
\cmidrule(r){5-7}
  $t$   & Date & $x_{1t}$ ({\sc{incid}}) & $x_{2t}$ ({\sc{trans}}) & $x_{3t}$ ({\sc{letha}}) & $x_{4t}$ ({\sc{wards}}) & $x_{5t}$ ({\sc{icu}})\\
\midrule
-474 & 2020-03-20  & 194   & 1.301 & 4.160 & 128  & 41   \\
-343 & 2020-07-31  & 197   & 0.978 & 1.140 & 340  & 41   \\
-197 & 2020-12-24  & 3574  & 0.987 & 2.180 & 2348 & 505  \\
-166 & 2021-01-24  & 12341 & 1.039 & 3.460 & 5375 & 742  \\
 0   & 2021-07-10  & 3658  & 1.042 & 0.382 & 488  & 144  \\
\bottomrule
\end{tabular}
\caption{Performances levels for five moments of the pandemic}
\label{tab:performance}
\end{table}
\end{center}

From the data in Table \ref{tab:performance} and by applying the previously constructed piecewise linear value functions, we obtained the results shown in Table \ref{tab:values}. The last column of this table provides the overall value of each moment of the pandemic, after applying the Model (\ref{eq:add_model}) formula.

\begin{center}
\begin{table}[htbp]
\centering
\begin{tabular}{cccccccccc}
\toprule
& & \multicolumn{2}{c}{Pillar I ({\sc{ACT}})} & \multicolumn{3}{c}{Pillar II ({\sc{SEV}})} & \multicolumn{1}{c}{Value} \\
\cmidrule(r){3-4}
\cmidrule(r){5-7}
\cmidrule(r){8-8}
     &  & $v_1(x_{1t})$ & $v_2(x_{2t})$ & $v_3(x_{3t})$ & $v_4(x_{4t})$  & $v_5(x_{5t})$  & \\
   $t$  & Date    & $w_1= 0.280$ & $w_2=0.141$  & $w_3=0.193$  & $w_4=0.193$   & $w_5=0.193$   & $v(t)$ \\
\midrule
-471 & 2020-03-20  & 3.441 & 180.00 & 115.571 & 1.0240 & 4.300 & 49.6800 \\
-342 & 2020-07-31  & 3.503 & 60.900 & 31.6990 & 2.7200 & 4.300 & 17.0400 \\
-197 & 2020-12-24  & 180.0 & 76.702 & 60.5640 & 89.056 & 180.0 & 124.832 \\
-166 & 2021-01-24  & 180.0 & 180.00 & 96.1120 & 180.00 & 180.0 & 163.810 \\
 0   & 2021-07-10  & 180.0 & 180.00 & 10.6240 & 3.9040 & 52.80 & 88.7700 \\
\bottomrule
\end{tabular}
\caption{Value functions scores for the five moments of the pandemic}
\label{tab:values}
\end{table}
\end{center}

Our pandemic indicator, PACI, reached its highest value in January 2021 and the lowest in July 2020. The four first moments of this example are displayed in Figure \ref{fig: pieces} (Appendix) and were used to test the experts and some anonymous people about the validity of the indicator.

\subsection{Chromatic classification system (ordinal scales)}\label{sec:ordinalscale}
\noindent The chromatic classification system is a tool that makes use of colours for better visualising the ordinal scale built with the experts. The colours selected for our model were inspired by the ones used in the RM of the Portuguese health authorities since the Portuguese population was already familiar with them.

Five fundamental states were defined with the experts: residual, alert, alarm, critical, and break. In addition, two more states were considered at the extremes, a baseline (the very lowest one) and a saturation or emergency state (the highest). All these states are zones defined in between two consecutive levels or cut-off lines:

    \begin{itemize}[label={--}]
      \item \textit{Baseline level} (cut-off line value $= 0$). The five performance baseline levels are presented in the following list: $[0,0,0,0,0]$. Each baseline level has a precise meaning for the experts. In this case, the first two mean there was  no pandemic activity recorded over the last seven days, which does not mean, of course, the pandemic was extinct, but simply that we did not register activity over the past seven days. The other three values, mean there were no deaths over the last seven days and there were no hospitalised {\sc{Covid-19}} patients.
      \item \textit{Residual level} (cut-off line value $= 10$). The five performance residual levels are presented in the following list: $[338,0.93, 0.36, 750,60]$. The experts were given this list, the values of each level on each value function, and its adequacy to represent an overall value of $10$. These elements were validated by the experts
      \item \textit{Alert level} (cut-off line value $= 40$). The five performance alert levels were presented in the following list: $[707,0.963,1.43,1571,126]$. The discussion with the experts was performed as in the previous case.
      \item \textit{Alarm level} (cut-off line value $= 80$). The five performance alarm levels were presented in the following list: $[1000,0.989,2.89,2222,178]$. The discussion with the experts was performed as in the residual and alert cases.
      \item \textit{Critical level} (cut-off line value $= 100$). The five performance critical levels were presented in the following list: $[1125,1,3.6,2500,2727,200]$. As the baseline levels, these levels are reference levels for the experts and have a particular meaning.
      \item \textit{Break level} (cut-off line value $= 120$). The five performance break levels were presented in the following list: $[1227,1.009,4.31,2727,218]$. The discussion with the experts was performed as in the residual, alert, and alarm cases.
      \item \textit{Emergency level} (cut-off line value $= 180$). The five performance saturation levels are resented in the following list: $[1506,1.034,6.47,3346,268]$. In this case, the experts agreed that any performance higher than the ones presented in the list will be considered as serious as the ones in the list. This corresponds to what we considered a saturation level
    \end{itemize}

The five fundamental states can be represented as in Figure 4. The transition between colours or states is not necessarily abrupt.

    \begin{figure}[H]
    \centering
    \includegraphics[width=17.5cm,height=3cm]{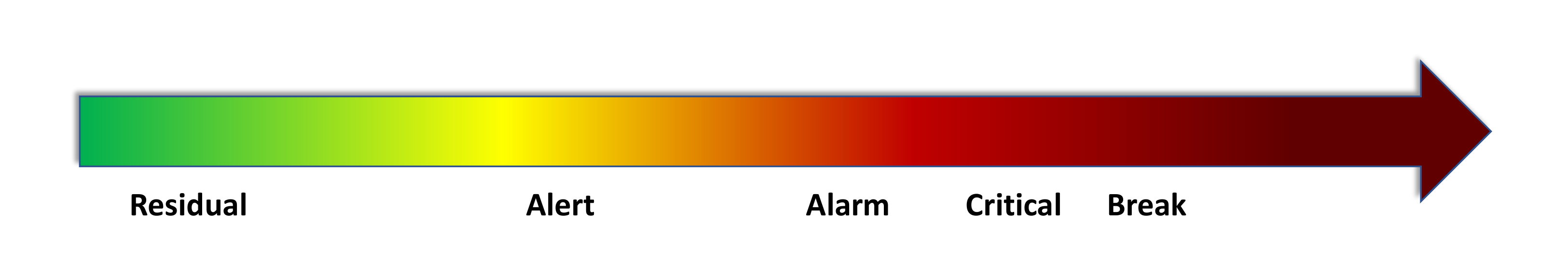}
    \caption{Chromatic classification system}\label{fig: chromatic}
    \end{figure}

A smooth transition can be considered since the policy/decision-makers cannot necessarily make the decisions automatically, after moving to a different state. It is important to see the evolution of the pandemic over subsequent days, after definitely moving from the current to a new state and implementing the measures/recommendations for this new state. A lower and an upper threshold for each cut-off line could be considered instead to make a smooth transition possible.

As can also be seen in the figure, and given the way the value functions were built, whenever we move to the next state, there is less room to make the decisions, i.e., it moves quickly, for example, from the alarm state to the critical state, than from the residual state to the alert state. This feature was a strong requirement of the experts. The main states can be briefly defined as follows, with respect to their impact on the health system:

\begin{itemize}[label={--}]
  \item \textit{Residual}: Absent or minimal pandemic activity without any impact on health structures (i.e., at the normal operating level) and without compromising the system tolerance.
  \item \textit{Alert}: Mild pandemic activity, still without impact on the normal activity of health structures, but reaching the usual flexibility, adaptability and safety tolerance threshold (e.g., increase in the emergency room visits and/or in the occupancy rate of hospital admissions).
  \item \textit{Alarm}: Moderate pandemic activity, already impacting the normal activity of health structures, with reallocation of technical and human resources and commitment to other health needs, reaching the functional reserve threshold.
  \item \textit{Critical}: Strong pandemic activity, having already exceeded the system's reserve threshold, conditioning effort and disruption in the activity of health structures allocated almost exclusively to the pandemic.
  \item \textit{Emergency}: Very strong pandemic activity and imminent collapse of health structures.
\end{itemize}

\begin{remark}{(Fragility Point 4)} \textit{Cut-off lines subjectivity}. This is in line with the previous two fragility points. The definition of the cut-off lines is subjective since they result from a co-constructive interactive process with the experts. However, defining thresholds for modelling a smooth transition between successive states can mitigate the subjectivity behind the definition of these cut-off lines.
\end{remark}

\subsection{Graphical model for visualization and communication}\label{sec:agregation}
\noindent One of the main features of our model is the visualisation functionalities to enable easy communication with the general public. Apart from other minor graphical functionalities, four types of graphical tools were developed:

\begin{enumerate}
  \item A graphic which displays the evolution of the indicator behaviour with coloured states and cut-off lines to separate each state (see Figure \ref{fig: Evolution}).
  \item An animation graphical tool with the cumulative contribution of each criterion to the pandemic (see Figure \ref{fig:cumulative}).
  \item A graphical representation of the (positive) impact of the vaccination plan to mitigate the progression of the pandemic in the country (see Figure \ref{fig:vaccination}).
  \item A state chromatic line, as in Figure \ref{fig: chromatic}.
\end{enumerate}

More details about these graphical tools will be provided in the next section.


\section{Results, sensitivity analyses, and simulations}\label{sec:results_overall}
\noindent This section is devoted to the implementation issues and verifications tests, results presentation, their validation, sensitivity analyses, and some final comments.

\subsection{Implementation issues and verification tests}\label{sec:results}
\noindent Our application was coded in the software Wolfram Mathematica, version 12.0\footnote{https://writings.stephenwolfram.com/2019/04/version-12-launches-today-big-jump-for-wolfram-language-and-mathematica/ }. All the functionalities of the Mathematica code were verified checked in several small examples with particularly extreme and pathological cases. This step includes verifications in thechecks on the debugging, input of criteria performance levels parameters, calculation of the criteria performance levels, calculation of the value functions and weights parameters, calculation of the comprehensive values for each time unit, all the graphical models outputs, sensitivity analyses, as well as the checking whetherr all the logical structure of the models  was correctly represented on the computer. The entire application has been designed to translate convert all the three models (criteria model, MAVT aggregation model, and graphical visualiszation and communication model) successfully in their entirety, as well as some additional functionalities for validation, simulation, and other sensitivity analyses purposes.

A Microsoft Excel version of the model was also implemented with fewer functionalities. The computations of the Excel PACI model are available for public consultation on an IST website (indicadorcovidl9.técnico.ulisboa.pt/) and on the PMA’s website (ordemdosmedicos.pt/iap/). This software automatically computes the daily changes in all the five criteria performance levels and the actual transmission rate, which is computed using the Robert Koch Institute formula (see \citealt{Koch2020}).

\subsection{Results}\label{sec:results}
\noindent The results provided information on three main aspects: the pandemic evolution, the cumulative contribution of each criterion to the evolution, and the impact of the vaccination plan. Please note that PACI can be used for forecasting if there is a good prediction techniques for the raw data used in the five criteria formulas.

\subsubsection{Pandemic evolution}\label{sec:evolution}
\noindent Before running the MAVT-based model with our set of criteria, we tested it with the two RM criteria, by setting the baseline, critical, and emergency levels as we did for our PACI tool, and considering linear value functions (since moving from an $R(t) = 0.1$ to an $R(t) = 0.2$ has the same impact as moving from an $R(t) = 0.9$ to an $R(t) = 1.0$, and this is true for any location along the $R(t)$ scale; the same reasoning applies to the incidence criterion). We also considered equal weights for each criterion. For other cut-off lines we set them as in our model. The evolution line can be observed in Figure \ref{fig: RM_Evolution} (Appendix). It is clear that this MAVT-based model does not adequately represent how the Portuguese people felt about the impact of the pandemic. The impact is always rather high, above the cut-off line in the alarm state. For the experts and the general public, the low impact at some points of the pandemic, as for example, in August 2020, cannot be seen clearly on this evolution line.

\begin{figure}[H]
    \centering
    \includegraphics[width=14cm,height=8cm]{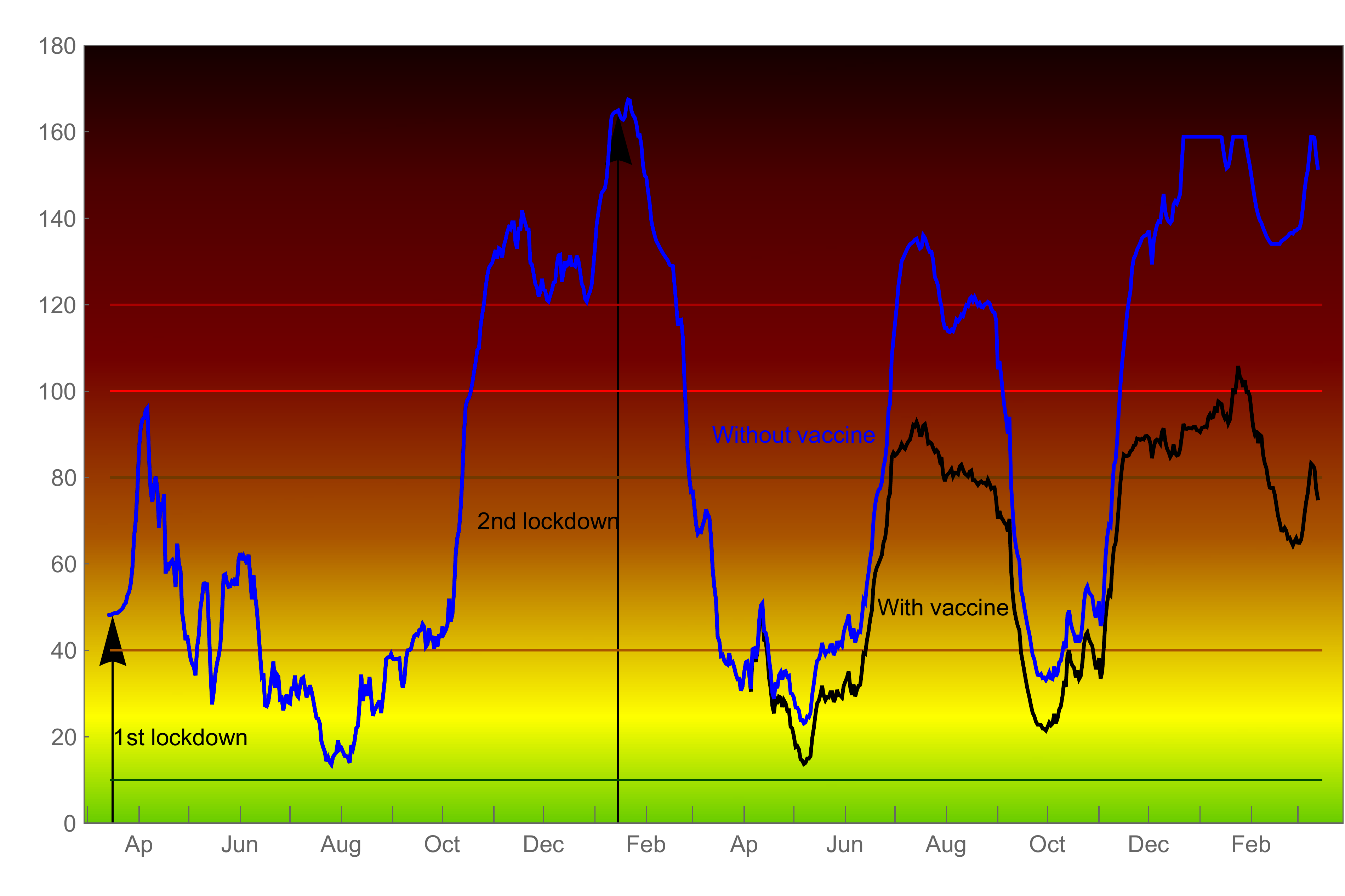}
    \caption{Evolution of the PACI till 2022-03-13}\label{fig: Evolution}
    \end{figure}

In Figure \ref{fig: Evolution} shows the global evolution of the PACI values in Portugal, along with the cut-off lines for separating the chromatic states.

The indicator gives an intuitive human perception of the impact of the pandemic in Portugal, particularly in the four main waves. The lockdown occurred exactly during the first wave, when the PACI value was close to $50$ points. The inertia of the system brought the PACI to a maximum on 6 April 2020, with $95$ points. The summer of 2020 was relatively mild with minimal values of the PACI value near to 14 between 26 July and 6 August 2020. We reached intolerable values between November 2020 and January 2021. In terms of pure pandemic impact, we can see that  the autumn of 2020 and the winter of 2021 corresponded to the same pandemic wave. The day of maximum  impact of the pandemic in Portugal, with $167.5$ points, was precisely 21 January 2021, which preceded the reinforcement of the lockdown measures/recommendations on 22 January 2021. The sharp drop from 25 January 2021 was caused by the lockdown of 15 January 2021. This decrease was reinforced substantially from 6 February 2021, when the population began to feel the effects of closing the schools  on the incidence and transmission rate values. After April 2021 and due to the vaccination effect, there was a clear detachment in the curves taking or not the new setting into account. The combined results of vaccination and non-pharmacological measures was maximal on 7 May, 2021, with 13.4 points, corresponding to the absolute minimum of the indicator since the appearance of {\sc{Covid-19}} in Portugal. After that date, we noticed an increase in the indicator, related to the surge in the Delta variant in the country, reaching a local maximum of $92.3$ on July 9, 2021. In the absence of lockdown, after that date, there was a slow trend in the decay of the impact, which can be related to the positive evolution of the vaccination in Portugal. With the surge in the Omicron variant, there was a strong increase in the PACI, which was related principally to the incidence and transmissibility in late November 2021. There was a fourth pandemic wave which differed from the previous ones in terms of severity. The maximum level of the last wave was attained between 21 and 28 January 2022, when values exceeded 100 points and reached a maximum of $106.4$ on 24 January.

\subsubsection{Cumulative contribution of each criterion}\label{sec:comulative}
\begin{figure}[H]
    \centering
    \includegraphics[width=14cm,height=8cm]{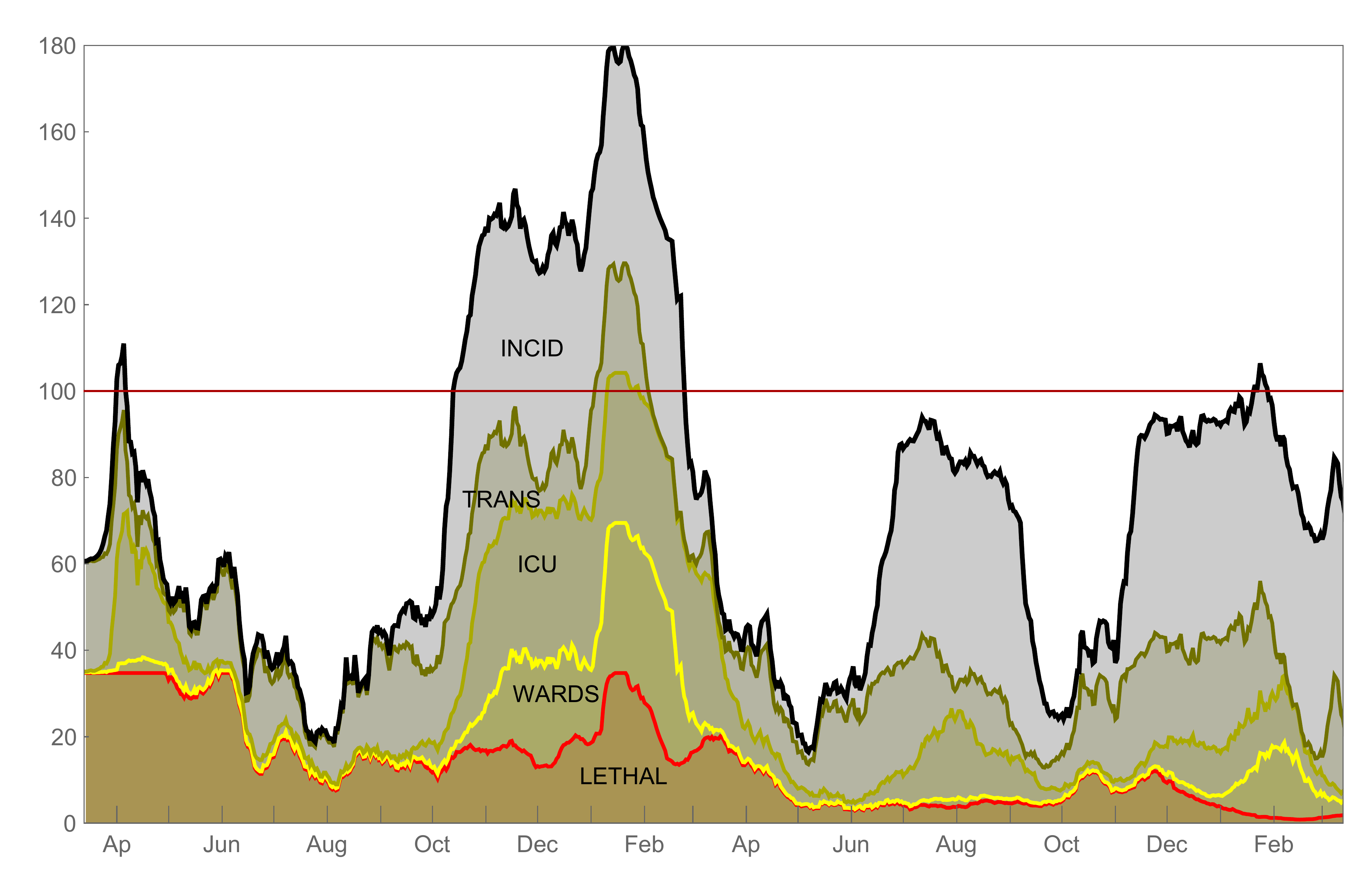}
    \caption{Cumulative contribution for the PACI till 2022-03-13}\label{fig:cumulative}
    \end{figure}

\noindent  In Figure \ref{fig:cumulative} we can observe the individual contributions of each criterion to the PACI value in Portugal.

The long term and gradual decrease in the contribution of lethality (case fatality rate, or {\sc{lethal}}  line) to the indicator is noteworthy. The contribution of this criterion increased dramatically during the January 2021 crisis. The role of the occupancy of the general number of beds for patients admitted to wards ({\sc{wards}} line) was particularly significant during the second pandemic wave in Portugal (between October 2020 and February 2021). This partial contribution reacts strongly when the spread of the disease is out of control. It is interesting to note that the patients admitted to intensive care beds ({\sc{icu}} line) criterion contribution to the PACI is significant in each pandemic wave. ICU bed occupancy grew after the increase in incidence ({\sc{incid}} line) with a delay of $10–12$ days.

Nevertheless, the relationship between the two criteria (incidence and ICU bed occupancy) is very clear and appeared in every pandemic wave, including the last wave related to the Omicron variant. Naturally, with the increment of vaccination, the relation between ICU occupancy and incidence dropped with time. In contrast to the ICU bed criterion contribution, the growth rate contribution ({\sc{trans}} line) appeared at an early stage of each wave and was the first alarm sign of a future increase in incidence, which was natural and expected. For instance, in the first wave of 2020, the PACI was mainly due to the growth rate and the case fatality rate for the first days of the {\sc{Covid-19}} infection in Portugal. The same effect is clear in the second wave, in October, in the last wave before June 2021 and, finally, in the last wave.

Finally, the incidence contribution to the PACI was severe in the months between October 2020 and February 2021 and significant with the last wave, due to the Omicron variant. The softening of governmental control measures/recommendations in Portugal at Christmas 2020 took place when the contribution of the incidence was high. The softening of measures in the last wave, in contrast, was appropriate, due to the protective effect of previous natural immunity and vaccination amongst the Portuguese population.

The introduction and the effects of the Delta variant are visible in the contribution of the incidence to the PACI in June 2021 and the effects of the Omicron variant are evident in the last pandemic wave in November 2021. Fortunately, this increase in the incidence contribution was balanced by the drop in the case fatality rate, overall number of patients admitted to wards occupancy, and number of beds for patients admitted to ICU occupancy relative to the values before the introduction of the Delta and Omicron variants and generalisation of the vaccination in Portugal.

\subsubsection{The impact of vaccination}\label{sec:vaccination}
\begin{figure}[H]
    \centering
    \includegraphics[width=14cm,height=8cm]{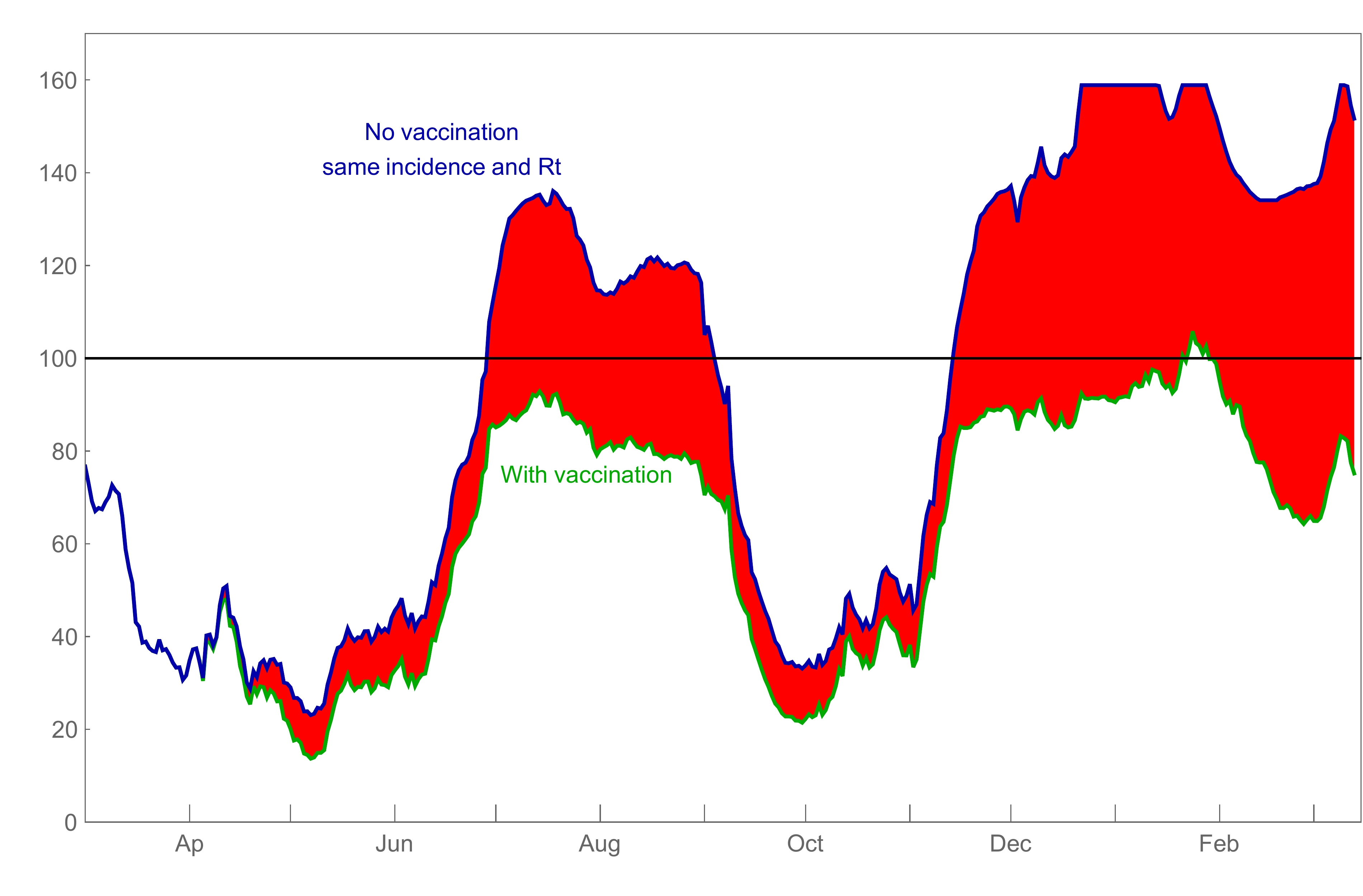}
    \caption{The vaccination impact till 2022-03-13}\label{fig:vaccination}
    \end{figure}

\noindent  Figure \ref{fig:vaccination} shows a comparison of the PACI values for actual data/parameters (lower curve) and an estimation of the indicator computed without the introduction of vaccination (upper curve). The upper curve was computed with the same observed incidence and growth rate of the actual PACI (lower curve), i.e., with the two criteria of Pillar I, but with no reduction in the severity of the disease (see the methodology described below), i.e., without the three criteria of Pillar II. The upper curve is, naturally, a lower bound estimate of the indicator without vaccination, since the immunisation process related to vaccination also the incidence and growth rate.

The effect of vaccination began to be measurable after $390$ days of PACI computation, i.e., after 4 April. The methodology to obtain this lower bound estimate of the indicator without vaccination, is to consider the averages of the proportion of ICU beds and general number of patients admitted to wards relative to the incidence before 4 April, and fix that rate, instead of using the actual data. Finally, in the case of fatality, we used the average of the same first $390$ days of the pandemic in Portugal. At this point, on day $390$, there was no difference between the values of the actual PACI indicator and the values of the lower estimate PACI without vaccination. After that date the values began to differ. We noticed a dramatic increase in the difference between the curves with time after 4 April. Naturally, since vaccination also reduced incidence and growth, after some time, more precisely on 15 July 2021, the indicator PACI began  to decrease, followed immediately by the lower bound estimate of PACI without vaccination.  The decoupling of severity and incidence in the last pandemic wave is noteworthy. Without vaccination, the last wave, due to the Omicron variant, could have had dramatic effects as we can see in the lower bound estimate of PACI without vaccination in Figure \ref{fig:vaccination}.

\subsection{Sensitivity analyses}\label{sec:robustness}
\noindent Sensitivity analyses conducted in this study are in line with the definition proposed in \cite{Roy2010}, where a change in the results was observed after a (simultaneous) change in all the parameters that are affected by some fragility aspects (see Remarks 1, 2, and 3). There are two major types of sensitivity analyses, one based on simulation (also called pseudo-sensitivity analysis) and another based on an exact characterisation of the effects of the changes in parameters. We detail these two types of sensitivity analyses in the following paragraphs.

\subsubsection{Simulation on the weights}\label{sec:robweigths}
\noindent The fragility point of Remark 3, which is related to subjectivity when assessing the weights of criteria, is one of the most critical fragility points in practice. A simulation analysis (also called pseudo-sensitivity) is dominated by an exact sensitivity analysis, as explained in the next subsection, but it has the advantage of being able to produce a large set of lines. Their shape conveys an idea of the evolution of our indicator (just for an illustrative purpose, see the last figure in Appendix, i.e., Figure \ref{fig: rob_add}, with a $\pm20\%$ change in the weights used in the application of PACI to the Portuguese pandemic).

In order to study the sensitivity of the results, we performed a strong change in the weights, thus enabling a variation in the range of each one, from $0$ to $1$ (the sum of all weights being equal to $1$). A Monte Carlo simulation made this study possible. Figure \ref{fig: simul_1} (Appendix) displays $400$ lines, among the $10\;000$ simulations performed (the representation of more lines is time consuming and led to a software crash due to a limited memory capacity). In this figure, the shape of all the lines is more or less the same. There is a drastic variation in the weights, since the simulations go from an extremely unbalanced situation, where only one criterion counts for the impact on the pandemic (the one with weight equal to one), to a complete balanced situation, where all the criteria contribute equally to the impact (equal weights for all criteria). A more realistic simulation, with a perturbation of $5\%$ below the weights values elicited with PaCo-DMC (for each day) and the same $5\%$ above the elicited weights, produces the curves (also $400$) represented in Figure \ref{fig: simul_2}. All the lines almost coincide with the shape of the curve of our PACI model. The results are quite robust for realistic variations in the weights, which shows the adequate behaviour of the PACI tool.

\subsubsection{Exact sensitivity analysis}\label{sec:all}
\noindent The exact sensitivity analysis was performed by taking into account all the critical fragility points of Remarks 1, 2, and 3. We made changes in the data provided by the incidence, transmission, and lethality criteria (the other two were not affected by strong imprecision, i.e., the number of beds in number of patients admitted to wards and number of patients admitted to ICU is relatively precise), in the five value functions, given the subjectivity when building them, and in the weights of criteria, for the same reason as the previous one.

A first strong (direct) perturbation on all the data/parameters related to the first three criteria and the five value functions revealed: $10\%$ below and above their daily performance levels/values. As for all the criteria weights, a variation of $10\%$ enabled us to build a polyhedron in a five-dimensional space. In this polyhedron, for each day, there was a maximum and a minimum value for the indicator by using linear programming techniques. For the maximum value, we considered a $+10\%$ change in the performance levels of the first three criteria and the values of the five value functions and computed the maximum of the indicator formula over the polyhedron of the weights; this is done for each day (see the upper envelope curve in Figure \ref{fig: rob_1}, Appendix). For computing the minimum for each day, by considering a $-10\%$ change in the performances levels of the first three criteria and values of the five value functions (see the lower envelope curve in Figure \ref{fig: rob_2}, Appendix). The difference of the upper and lower curves is in an average of $47.12$ points with a standard deviation of $11.0848$ points.

Then, a more realistic sensitivity analysis with a $\pm5\%$ change in the performance levels of the first three criteria and the values of the five value functions, and a similar construction of the polyhedron for the weights was performed. The results are presented in Figure \ref{fig: rob_3} and Figure \ref{fig: rob_4} (Appendix). The average of the difference is now $28.2366$ impact points and the standard deviation is $6.9615$ points.

The outline of the exact sensitivity procedure can be presented as follows:
\begin{enumerate}
  \item Take a certain time unit, $t$ (e.g., a given day during the pandemic).
  \item Consider the first three criteria: $g_1$, $g_2$, and $g_3$. Proceed as follows (\textit{Fragility Point 1}):
    \begin{enumerate}
        \item Consider their performance levels at time $t$: $x_{1t}$, $x_{2t}$, $x_{3t}$.
        \item Make a decrease of these levels. Let $x^{-}_{1t}$, $x^{-}_{2t}$, and $x^{-}_{3t}$, denote the new performance levels.
    \end{enumerate}
  \item Consider all value functions: $v_1$, $v_2$, $v_3$, $v_4$, and $v_5$. Proceed as follows (\textit{Fragility Point 2}):
    \begin{enumerate}
        \item Consider their values at time $t$, tanking into account the modified performance levels of the first three criteria (as we did in the previous step) and the  performance levels of the last two criteria: $v_{1}(x^{-}_{1t})$, $v_{2}(x^{-}_{2t})$, $v_{3}(x^{-}_{3t})$, $v_{4}(x_{4t})$, and $v_{5}(x_{5t})$.
        \item Decrease these values. Let $v^{-}_{1t}$, $v^{-}_{2t}$, $v^{-}_{3t}$, $v^{-}_{4t}$, and $v^{-}_{t}$, denote the new values.
        \item (for the sake of simplicity of notation) Consider the vector, $v^{(t)-} = (v^{-}_{1t}, v^{-}_{2t}, v^{-}_{3t}, v^{-}_{4t}, v^{-}_{5t})$.
    \end{enumerate}
  \item Consider all weights: $w_1$, $w_2$, $w_1$, $w_2$, and $w_3$. Proceed as follows (\textit{Fragility Point 3}):
    \begin{enumerate}
        \item Decrease the weights. Let $w^{-}_1$, $w^{-}_2$, $w^{-}_3$, $w^{-}_4$, and $w^{-}_5$, denote the lower bound values for the weights.
        \item Increase of the weights. Let $w^{+}_1$, $w^{+}_2$, $w^{+}_3$, $w^{+}_4$, and $w^{+}_5$, denote the upper bound values for the weights.
        \item Remark: These changes in the weights are not indexed to the time period; they are thus valid for any $t$.
        \item Construct the polyhedron of the weights, denoted by $W$, as the result of the intersection on the following constraints:
            \begin{itemize}[label={--}]
              \item Bounding constraints: $w^{-}_l \leqslant w_j \leqslant w^{+}_j$, for $j=1,\ldots,5$.
              \item Normalization constraint: $w_1 + w_2 + w_3 + w_4 + w_5 = 1$.
              \item Consistence constraints: $0 \leqslant w_j \leqslant 1$, for $j=1,\ldots,5$ (these constraints avoid to have negative weights of weights with values strictly greater than one).
            \end{itemize}
        \item (for the sake of simplicity of notation) Consider the vector, $w = (w_1, w_2, w_3, w_4, w_5)$. Any feasible $w$ is an element of $W$, i.e., $w \in W$.
    \end{enumerate}
  \item Solve the following linear programming problem:
        \[
            v^{-}(t) = \min\left\{w^{\top}v^{(t)-} \;:\; w \in W \right\},
        \]
        where $v^{-}(t)$ is the lowest (which, in our case, corresponds to the best) value of the PACI model.
  \item Proceed in a similar way to obtain $v^{+}(t)$, i.e., the worst value of the PACI model for time unit $t$.
\end{enumerate}

\subsection{Validation, adjustments, and comments}\label{sec:validation}
\noindent In our case, the validation consisted of presenting the actors (mainly the experts), but also some anonymous people, with the shape of the pandemic indicator and obtaining some comments to validate or make adjustments to our models. All the tests were conducted before and after performing the sensitivity analysis, but they have more credibility after such an analysis is performed. More precisely, we performed the following tests:

\begin{enumerate}
  \item In a first step, we built a figure displaying different moments of the pandemic in the country, see Figure  \ref{fig: pieces} (Appendix). The moments were not chronologically ordered. We asked the experts to look at the figure and tell us if they were able to identify such moments and relate them to a state of the pandemic. We selected the following moments: the beginning of the pandemic, July 2020, January 2021, and Christmas 2020. All the experts were able to easily identify all the moments. Only a slight hesitation was recorded for one expert with respect to the moment related to the start of the pandemic. The team members also asked the same question to some anonymous people. We performed tests with 30 individual, mostly university administrative staff ($10$ individuals), students ($10$ individuals), and random people from the streets of Lisbon ($10$ individuals). We gave some additional explanations about the fact that the impact was represented on a scale from $0$ to $180$ points as well as the minimum and the maximum PACI values. The number of people asked to perform the validation testes was rather low. A more systematic and complete study would be important to obtain more input for validating our PACI model. Those presented with the four pieces of our graphic of Figure \ref{fig: pieces}, only the start of the pandemic lead to some hesitations. The lowest and the highest impact of the pandemic were easily identified, and Christmas 2020 was almost always recognised.
  \item In a second step, we showed experts the whole evolution curve since the beginning of the pandemic. They commented and justified all the moments of the pandemic and the different critical situations, i.e., the waves that occurred during the disease evolution in Portugal, namely the initial growth, the autumn and winter crisis and, finally, the surge of the Delta variant. They were also able to identify the calm situation of the spring/summer of 2020 and the relaxation of April 2021. The same exercise was also performed with the anonymous people. Most of the tests were positive, with $28$ people being able to comment the initial stage, the autumn/winter crisis, and the Delta variant wave.
  \item In a third step, we asked the experts to comment about the reasons that lead to some moments of the curve and to explain the reasons leading to such behaviour in the PACI values. We also chose some particular points on the curve and asked the experts to comment and justify them. This is a different approach from the previous one, since the moments were not chosen by the experts, but by the analysts.
  \item Finally, we asked the experts to provide some raw data characteristic of each state, run the model, and show the results. For example, the following list of performance levels $[1250, \, 1.02, \, 2.8, \, 2235,$ $195]$ is a profile that should be considered in the critical state. After running the model, we get the value $104.2$, which is within the critical state of the chromatic model. Almost all the results led to the state provided by the experts. This test led to a very slight adjustment of the weights (see point 6 of subsection 3.3.2.).
\end{enumerate}


\section{Lessons learned}\label{sec:recommendations}
\noindent  This section is devoted to understanding the lessons learned from the use of the PACI since the beginning of its application for modelling the impact of COVID-19 in Portugal. We present sequentially, in the next subsections, the impact of PACI on \textit{media} and on the Portuguese general population, the success features of the tool as well as its failures and, finally, the aspects we can consider for improving the PACI in the future.

\subsection{Impact on \textit{media} and on the general public}\label{sec:ImpactPACI}
\noindent  The PACI was presented at a press conference of the PMA and the IST on 14 July 2021 and was mentioned by all the media. For example, the two largest daily Portuguese newspapers, Diário de Notícias\footnote{https://www.dn.pt/sociedade/especialistas-apresentam-novo-indicador-para-a-matriz-de-risco-13937495.html}   and Público\footnote{https://www.publico.pt/2021/07/14/sociedade/noticia/ordem-medicos-propoe-nova-matriz-risco-funciona-1970375}, the largest national weekly, Expresso\footnote{Example: https://expresso.pt/coronavirus/covid-19-portugal-ultrapassa-critical-level-of-severity-of-the-pandemic/}, and all Portuguese televisions channels, including for example, RTP\footnote{ https://www.rtp.pt/noticias/pais/especialistas-apresentam-novo-indicador-que-esperam-venha-a-ser-matriz\_n1335331}, reported on the press conference and interviewed several members of the team. During the presentation, several questions were asked about the PACI and the reception of the indicator by the \textit{media} was very positive and highly praised. Although it has not replaced the official RM, the PACI has been regularly used by the \textit{media} for a better characterisation of pandemic activity at a national level  and mainly to monitor the different phases of the pandemic, justify the allocation of resources and, above all, explain the situation to the population in order to guarantee their involvement in prevention and control measures for the transmission. The use of the PACI has been frequent and regular and has justified equally praiseworthy references to the PACI by regular commentators from different media. About one year after its presentation, the PACI has become a simple, reliable, and useful indicator for assessing the impact of the pandemic that is regularly used by the \textit{media} in the news for the public. The PACI has proved to be an important tool at a national level for the assessment and control of the pandemic.

\subsection{Learning from Pandemic Assessment Composite Indicator successes}\label{sec:SuccessesPACI}
\noindent The successes of the PACI indicator can be considered at several levels, as follows.
\begin{enumerate}
    \item \textit{Criteria model}. The criteria model comprises the most relevant criteria built around the fundamental concepts followed two principles: familiarity and negligibility, as explained in subsection 3.2. Integrating more criteria would provide a more complex composite indicator, which would not add much more value as since the beginning of the pandemic none of the five criteria have been questioned; this reveals their acceptance and success.
    \item \textit{Aggregation model}. The additive aggregation model was also considered a success for the same reason; it was well accepted and not questioned even by the scientific community. There is, of course, room for improvement, such as considering different value functions and weights, but this was largely tested with sensitivity analysis and simulations.
    \item \textit{Socio-technical process}. The co-constructive interactive process between the group of experts and the group of analysts was also rather a successfully socio-technical process but given the time constraints this would benefit from more interaction to gain more confidence in some preference judgements. Subsection 5.4 will provide some pointers for future improvement.
    \item \textit{Communication model}. This communication model was probably one of the highest successes of our toll. Composed of three main graphical tools, the evolution graphic with the cut-off linear, the cumulative evolution graphic, and the graphic on the impact of the vaccination plan, make the information very informative for the general population. We regret, however, not being able to include the forecasting capabilities of our tool for a one or two week forecasting process, which would be extremely important for the short term planning of measures and recommendations.
    \item \textit{Sensitivity analysis and simulation tools}. These two tools were of utmost important for guaranteeing confidence in the indicator. Performing drastic and realistic changes in the model parameters, we observed the effects of the indicator and have a good idea of how strong it was.
    \item \textit{Validation tests}. The validation process with experts and some individuals (a rather small number) was important and also a success, but it was rather, given the time constraints. This is an aspect we would like to improve in the future and will be outlined in subsection 5.4.
\end{enumerate}

\subsection{Learn from Pandemic Assessment Composite Indicator failures}\label{sec:FailluresPACI}
\noindent  The failures of the PACI indicator can be considered at several levels, as follows.
\begin{enumerate}
    \item \textit{Criteria model}. A possible drawback of our indicator is that despite it being an early warning indicator of future changes, due to de presence of the incidence growth rate, it does not incorporate short term forecasting of the future evolution of the pandemic. This could be included in the indicator but has some drawbacks: computation is heavy; there are computational costs; it becomes more difficult to communicate to the public; and there is a margin of error in the predictions.
    \item \textit{Aggregation model}. As it is a more or less a quadratic non-linear model, the indicator reacts very quickly to an increase in the two pillars. This appears to be a drawback, i.e., when the pandemic approaches a hazardous status, the indicator increases very quickly giving little time for the authorities to react. But, in fact, this is a consequence of reality, rather than a failure of the indicator. When there is a near exponential explosion of cases, severity and transmissibility grow very sharply and the indicator reflects this behaviour. Whenever there is an increase in the indicator at moderate values, i.e., close to 40, there is a strong need for counter measures.
    \item \textit{Socio-technical process}. Despite being well accepted by the Portuguese society, including the Portuguese President of the Republic and media, the local Health Authorities did not officially change the evaluation of the pandemic in Portugal using the RM concept, as we explain in the conclusion to this article. The co-constructive interactive process between the group of experts and the group of analysts was also rather a successfully socio-technical process but give the time limitations this would have benefitted from greater interaction to gain more confidence in some preference judgements. Subsection 5.4 will provide some pointers for future improvement.
    \item \textit{Communication model}. The communication model, despite its success, would benefit from a more adequate web-based platform, with additional graphical tools and the possibility of graphical forecasting.
    \item \textit{Sensitivity analysis and simulation tools}. These tools were the adequate given the type of indicator designed. We cannot see any failure with respect to this point.
    \item \textit{Validation tests}. The validation testes are needed for a wider application in the general public given that the number of tests performed was rather low, as explained in the validation subsection of this paper. Also, here a web-based platform would be of utmost importance to conduct a more in depth study related to the validation tests.
\end{enumerate}

\subsection{Future improvements in the Pandemic Assessment Composite Indicator }\label{sec:ImprovementsPACI}
\noindent  Future improvements to the PACI can be summarised as follows:
\begin{enumerate}
    \item \textit{Criteria model}. The PACI model makes use of incidence criteria (new cases and transmission) and severity (number of total admissions and ICU admissions and lethality) in order to streamline a limited number of easily accessible parameters, and that enabled the monitoring of the social and health impact of the disease and the pandemic. The use of other indicators, such as the vaccination coverage rate, the number of tests and the percentage of positive testing, have an impact on both incidence and severity. They can be used in the future as additional criteria, but we considered that their use, in addition or replacement, did not add value to the indicator. A possible improvement would involve adding further criteria and checking their effect with respect to the current version of the PACI.
    \item \textit{Aggregation model}. The value function and the weights of each indicator was established in collaboration by the two teams and according to analyses of the impact of the pandemic in the period before the development of the indicator (period of more than 12 months). Although it was not a failure, we could add more value functions and reconstruct the set of the weights by considering \ref{table:criterion_g1} the possibility of adding more criteria, as mentioned in the previous point.
    \item \textit{Socio-technical process}. With respect to the point, and as was mentioned in the subsection failures, the same point, we would like to perform more interaction with experts and the population, and even text the deck of card method with the possibility of inconsistent judgements to model the impact judgements of the experts and participants in a more accurate manner.
    \item \textit{Communication model}. The improvements to the communication model are the ones mentioned in the subsection of the model failures. An interactive web-based platform would benefit the general population substantially and render the information clearer and make it more widely available.
    \item \textit{Sensitivity analysis and simulation tools}. This indicator was also tested in the flu epidemics prior to 2020, which confirmed the sensitivity and quality of the indicator. Due to the lack of available data, it was not possible to test the indicator in a pandemic prior to the current one, for example, in the 2009 influenza A (H1N1) pandemic. After its release in July 2021, the continuous use of the indicator lead to the assessment of its reliability. The need to adjust either the weights or the inclusion or substitution of new criteria has not been detected so far. However, we admit that the Portuguese reality may not be replicable in other countries, so it is essential to validate the PACI in other countries and, eventually, establish different weights for different realities. This is indisputably one of the objectives of the dissemination of the PACI internationally.
    \item \textit{Validation tests}. The tests should also follow the ideas pointed out in the failure subsection to mitigate the less effectiveness of the current PACI tool with respect to the validation tests.
\end{enumerate}


\section{Conclusions}\label{sec:context}
\noindent In this paper we proposed an innovative application of a MAVT additive model for building a PACI and a chromatic ordinal classification system to assist in the management of the {\sc{Covid-19}} pandemic in Portugal. This indicator was built by following a sociotechnical co-constructive interactive process between the CCIST and the
GCMO teams, and to the best of our knowledge, it is the first MAVT model proposed to analyse the evolution the pandemic and mitigate the impacts of {\sc{Covid-19}} in the world. It was designed with the particular purpose of answering several questions posed by the Portuguese population: How is the pandemic evolving in the country? In which pandemic state are we currently? What is the impact of the vaccination plan established by the government? All these questions have been answered, and our indicator had strong acceptation in Portugal. It still continues to be followed and frequently mentioned in the \textit{media}, even though the Portuguese health authorities did not adopt it officially as another indicator for effective policy/decision-making.

In Portugal, pandemic-related decisions are centralised at the Ministry of Health in conjunction with the other members of the Government and with the full involvement of the Prime Minister. The Ministry of Health is responsible for two official bodies of a technical nature, the Directorate- General for Health, and the National Institute of Public Health Dr Ricardo Jorge. They monitor the pandemic activity and design the rules and guidelines to fight the pandemic and minimise its impact on health systems and on the country’s social and economic life. At the beginning of the pandemic, the National Institute of Public Health Dr Ricardo Jorge developed an RM  based on the number of cases \textit{per} $100\;000$ inhabitants at 14 days (national and regional incidence) and the transmission, $R(t)$. Before the massive vaccination campaign of the Portuguese population, this matrix was a useful tool to support measures at both a national and regional level. With vaccination, the relationship between incidence and transmission and severity was lost. For this reason, the PMA and the IST developed the PACI to monitor activity and severity better and also support decisions with coherence and greater involvement of the general public. On several occasions, PMA and IST were in touch with the Ministry of Health to provide the use of the PACI. Despite much praise for the new indicator, official authorities chose to avoid public changes to the pandemic assessment model and keep the RM already in place. The PACI was reserved for more restricted use by official bodies in monitoring pandemic activity. This decision was also influenced by the fact that the National Institute of Public Health Dr Ricardo Jorge is an official body under the Ministry of Health and that the replacement of the RM in force could harm collaboration between the different official bodies.

Despite the fragility points related to the data and the construction of the aggregation model itself, it has several advantages that are widely acknowledged by academics, opinion-makers, media, and Portuguese general public. Even the technicalities of the method, including computations, can be reproducible for any reader with basic mathematical knowledge. The parameters of the aggregation model can be adjusted, if justified, during the pandemic evolution. This comprises the shape of the value functions, the weights, the cut-off lines, and the reference levels, in particular the critical level (if the number of beds in intensive care increases it is normal for the critical level related to the fifth criterion to change accordingly). Also, the formulas of the criteria model can be adjusted or replaced by more suitable ones (this also implies changes in the aggregation model).

The flexibility of the PACI mode opens up several avenues for possible future research:

\begin{itemize}[label={--}]
  \item This is an open model, in the sense that it can easily accommodate the inclusion of more criteria and even more pillars to account for other points of view, as for example the economic impact of the pandemic. Other aspects, for example, the interaction effect between/among criteria are also possible, but they are more sophisticated and require adaptations in the judgement assessment techniques, such as PaCo-DCM in the case of some additional multiplicative terms (see \citealt{KeeneyRa1993}). In the case of the Choquet aggregation model, PaCo-DCM could be more or less easy to adapt (see \citealt{BotteroEtAl2018}), but as all the value functions are between $0$ and $1$, it would require a re-scaling.
  \item This model can be applied to all territorial units (country, regions, counties, sets of counties, etc.) with available data and possibly with some readjustments of the critical levels of some criteria.
  \item Scalability to other countries is also a possibility, but all the criteria would have to be reconsidered, as well as all the levels, in particular, the baseline  critical levels, and the cut-off lines. Comparison with other countries would be of great importance in analysing the impact of different measures taken by other countries.
  \item The model can be applied to other diseases, and other health problems, and even in different sectors where the building of composite indicators is important.
  \item The model also has forecasting capabilities; it only needs to have good estimates of the raw data  $N(\cdot)$, $O(\cdot)$, $H(\cdot)$, and $U(\cdot)$.
  \item One of the most interesting avenues for future research is that of making use of constructive preference-learning techniques as some adaptations of the GRIP method (see \citealt{FigueiraEtAl2009})) for building the composite indicator. This is a kind of machine learning based tool that infers representative value functions from examples provided by the experts. After this first study, the experts now have a much better understanding of the entire issue and can easily provide a ``good'' set of examples for helping in the construction of the model parameters, value functions, weights, and even the cut-off lines.
\end{itemize}

Proposing the PACI was possible thanks to collaboration between the CCIST and the GCOM teams, which continue carry out the research proposed in the previous listed avenues for future work. In conclusion, it is important to note that this new indicator is not necessarily a competitor for the RM as both can be used at the same time to inform the Portuguese health authorities better when making decisions.


\section*{Acknowledgements}
\addcontentsline{toc}{section}{\numberline{}Acknowledgements}
\noindent  This paper is dedicated to the memory of our friend and colleague Carlos Alves (Full Professor of Instituto Superior Técnico).  The authors would like to acknowledge all the members of Crisis Office of Portuguese Medical Association (GCOM) for their valuable contributions, Professor Fernando Mira da Silva (IST) for providing the conditions to create the PACI webpage at the informatics platform of IST, Romana Borja-Santos, Ana Rodrigues, and Beatriz Santiago, for their appreciated and priceless help during all this process, and the anonymous people who helped us in validating the model. Jos\'e Rui Figueira acknowledges the project PTDC/EGE-OGE/30546/2017 (hSNS: Portuguese public hospital performance assessment using a multi-criteria decision analysis framework) supported by national funds through Funda\c{c}{\~a}o para a Ciência e a Tecnologia (FCT),  and Ana Paula Serro the project LISBOA-01-0145-FEDER-072536 (CAPTURE - Use of functionalized particles for enrichment and efficient detection of SARS-CoV-2 in clinical and environmental samples from Programa Operacional Regional de Lisboa). Henrique M. Oliveira was partially supported through FCT funding for CEMAT projects with reference UID/MAT/04459/2020. Finally, the authors also would like to acknowledge the anonymous reviewers for their valuable comments, insights, and suggestions, which proved to be vital for the improvement of the quality of this work.

\addcontentsline{toc}{section}{\numberline{}References}
\bibliographystyle{model2-names}
\bibliography{API_Indicator_1ST_R}

\begin{thebibliography}{51}
\expandafter\ifx\csname natexlab\endcsname\relax\def\natexlab#1{#1}\fi
\expandafter\ifx\csname url\endcsname\relax
  \def\url#1{\texttt{#1}}\fi
\expandafter\ifx\csname urlprefix\endcsname\relax\def\urlprefix{URL }\fi
\providecommand{\eprint}[2][]{\url{#2}}
\providecommand{\bibinfo}[2]{#2}
\ifx\xfnm\relax \def\xfnm[#1]{\unskip,\space#1}\fi
\bibitem[{Azevedo et~al.(2020)Azevedo, Pereira, Ribeiro and
  Soares}]{AzevedoEtAl2020}
\bibinfo{author}{Azevedo, L.}, \bibinfo{author}{Pereira, M.},
  \bibinfo{author}{Ribeiro, M.}, \bibinfo{author}{Soares, A.},
  \bibinfo{year}{2020}.
\newblock \bibinfo{title}{Geostatistical {COVID-19} infection risk maps for
  {P}ortugal}.
\newblock \bibinfo{journal}{International Journal of Health Geographics}
  \bibinfo{volume}{19}, \bibinfo{pages}{25}.
\bibitem[{Belton and Stewart(2002)}]{BeltonSt2002}
\bibinfo{author}{Belton, V.}, \bibinfo{author}{Stewart, {\relax T.J.}.},
  \bibinfo{year}{2002}.
\newblock \bibinfo{title}{Multiple Criteria Decision Analysis: An Integrated
  Approach}.
\newblock \bibinfo{publisher}{Kluwer Academic Publishers},
  \bibinfo{address}{Dordrecht, The Netherlands}.
\bibitem[{Bottero et~al.(2018)Bottero, Ferretti, Figueira, Greco and
  Roy}]{BotteroEtAl2018}
\bibinfo{author}{Bottero, M.}, \bibinfo{author}{Ferretti, V.},
  \bibinfo{author}{Figueira, {\relax J.R.}.}, \bibinfo{author}{Greco, S.},
  \bibinfo{author}{Roy, B.}, \bibinfo{year}{2018}.
\newblock \bibinfo{title}{On the {C}hoquet multiple criteria preference
  aggregation model: {T}heoretical and practical insights from a real-world
  application}.
\newblock \bibinfo{journal}{European Journal of Operational Research}
  \bibinfo{volume}{271}, \bibinfo{pages}{120--140}.
\bibitem[{Bouyssou and Pirlot(2016)}]{BouyssouPi2016}
\bibinfo{author}{Bouyssou, D.}, \bibinfo{author}{Pirlot, M.},
  \bibinfo{year}{2016}.
\newblock \bibinfo{title}{Conjoint measurement tools for {MCDM}: {A} brief
  introduction}, in: \bibinfo{editor}{Greco, S.}, \bibinfo{editor}{Ehrgott,
  M.}, \bibinfo{editor}{Figueira, J.} (Eds.), \bibinfo{booktitle}{Multiple
  Criteria Decision Analysis}. \bibinfo{publisher}{Springer Science+Business
  Media New York}, \bibinfo{address}{New York, NY, USA}, pp.
  \bibinfo{pages}{97--151}.
\bibitem[{Brans and De~Smet(2016)}]{BransDe2016}
\bibinfo{author}{Brans, {\relax J.-P.}.}, \bibinfo{author}{De~Smet, Y.},
  \bibinfo{year}{2016}.
\newblock \bibinfo{title}{{\sc{Promethee}} methods}, in:
  \bibinfo{editor}{Greco, S.}, \bibinfo{editor}{Ehrgott, M.},
  \bibinfo{editor}{Figueira, J.} (Eds.), \bibinfo{booktitle}{Multiple Criteria
  Decision Analysis}. \bibinfo{publisher}{Springer Science+Business Media New
  York}, \bibinfo{address}{New York, NY, USA}, pp. \bibinfo{pages}{187--219}.
\bibitem[{Català et~al.(2021)Català, Marchena, Conesa, Palacios, Urdiales,
  Alonso, Alvarez-Lacalle, Lopez, Cardona and Prats}]{CatalaEtAl2021}
\bibinfo{author}{Català, M.}, \bibinfo{author}{Marchena, M.},
  \bibinfo{author}{Conesa, D.}, \bibinfo{author}{Palacios, P.},
  \bibinfo{author}{Urdiales, T.}, \bibinfo{author}{Alonso, S.},
  \bibinfo{author}{Alvarez-Lacalle, E.}, \bibinfo{author}{Lopez, D.},
  \bibinfo{author}{Cardona, P.J.}, \bibinfo{author}{Prats, C.},
  \bibinfo{year}{2021}.
\newblock \bibinfo{title}{Monitoring and analysis of {COVID-19} pandemic: {T}he
  need for an empirical approach}.
\newblock \bibinfo{journal}{Frontiers in Public Health} \bibinfo{volume}{9},
  \bibinfo{pages}{633123}.
\bibitem[{Corrente et~al.(2021)Corrente, Figueira and Greco}]{Corrente2021}
\bibinfo{author}{Corrente, S.}, \bibinfo{author}{Figueira, {\relax J.R.}.},
  \bibinfo{author}{Greco, S.}, \bibinfo{year}{2021}.
\newblock \bibinfo{title}{Pairwise comparison tables within the deck of cards
  method in multiple criteria decision aiding}.
\newblock \bibinfo{journal}{European Journal of Operational Research}
  \bibinfo{volume}{291}, \bibinfo{pages}{738--756}.
\bibitem[{Bana~e Costa et~al.(2016)Bana~e Costa, De~Corte and
  Vansnick}]{BanaEtAl2016}
\bibinfo{author}{Bana~e Costa, C.}, \bibinfo{author}{De~Corte, {\relax
  J.-M.}.}, \bibinfo{author}{Vansnick, {\relax J.-C.}.}, \bibinfo{year}{2016}.
\newblock \bibinfo{title}{Mathematical foundations of {MACBETH}}, in:
  \bibinfo{editor}{Greco, S.}, \bibinfo{editor}{Ehrgott, M.},
  \bibinfo{editor}{Figueira, J.} (Eds.), \bibinfo{booktitle}{Multiple Criteria
  Decision Analysis}. \bibinfo{publisher}{Springer Science+Business Media New
  York}, \bibinfo{address}{New York, NY, USA}, pp. \bibinfo{pages}{421--463}.
\bibitem[{Cox(2008)}]{Cox2008}
\bibinfo{author}{Cox, {\relax A.Jr.}.}, \bibinfo{year}{2008}.
\newblock \bibinfo{title}{What's wrong with risk matrices?}
\newblock \bibinfo{journal}{Risk Analysis} \bibinfo{volume}{28},
  \bibinfo{pages}{497--512}.
\bibitem[{Dinis et~al.(2021)Dinis, Figueira and Teixeira}]{DinisEtAl2021}
\bibinfo{author}{Dinis, D.}, \bibinfo{author}{Figueira, {\relax J.R.}.},
  \bibinfo{author}{Teixeira, {\relax A.P.}.}, \bibinfo{year}{2021}.
\newblock \bibinfo{title}{A multiple criteria approach for ship risk
  classissication: {A}n alternative to the {P}aris {M}o{U} {S}hip {R}isk
  {P}rofile}.
\newblock \bibinfo{journal}{CoRR} \bibinfo{volume}{abs/2107.07581}.
\newblock \eprint{2107.07581}.
\bibitem[{Dlamini et~al.(2020)Dlamini, Dlamini, Mabaso and
  Simelane}]{DlaminiEtAl2020}
\bibinfo{author}{Dlamini, W.}, \bibinfo{author}{Dlamini, S.},
  \bibinfo{author}{Mabaso, S.}, \bibinfo{author}{Simelane, S.},
  \bibinfo{year}{2020}.
\newblock \bibinfo{title}{Spatial risk assessment of an emerging pandemic under
  data scarcity: {A} case of {COVID-19} in {E}swatini}.
\newblock \bibinfo{journal}{Applied Geography} \bibinfo{volume}{125},
  \bibinfo{pages}{102358}.
\bibitem[{Doumpos and Zopounidis(2002)}]{DoumposZo2002}
\bibinfo{author}{Doumpos, M.}, \bibinfo{author}{Zopounidis, C.},
  \bibinfo{year}{2002}.
\newblock \bibinfo{title}{Multicriteria Decision Aid Classification Methods}.
\newblock \bibinfo{publisher}{Kluwer Academic Publishers},
  \bibinfo{address}{Dordrecht, The Netherlands}.
\bibitem[{Dubois and Perny(2016)}]{DuboisPe2016}
\bibinfo{author}{Dubois, D.}, \bibinfo{author}{Perny, P.},
  \bibinfo{year}{2016}.
\newblock \bibinfo{title}{A review of fuzzy sets in decision sciences:
  {A}chievements, limitations and perspectives}, in: \bibinfo{editor}{Greco,
  S.}, \bibinfo{editor}{Ehrgott, M.}, \bibinfo{editor}{Figueira, J.} (Eds.),
  \bibinfo{booktitle}{Multiple Criteria Decision Analysis}.
  \bibinfo{publisher}{Springer Science+Business Media New York},
  \bibinfo{address}{New York, NY, USA}, pp. \bibinfo{pages}{637--691}.
\bibitem[{Dyer(2016)}]{Dyer2016}
\bibinfo{author}{Dyer, J.}, \bibinfo{year}{2016}.
\newblock \bibinfo{title}{Multiattribute utility theory (maut)}, in:
  \bibinfo{editor}{Greco, S.}, \bibinfo{editor}{Ehrgott, M.},
  \bibinfo{editor}{Figueira, J.} (Eds.), \bibinfo{booktitle}{Multiple Criteria
  Decision Analysis}. \bibinfo{publisher}{Springer Science+Business Media New
  York}, \bibinfo{address}{New York, NY, USA}, pp. \bibinfo{pages}{285--314}.
\bibitem[{Dyer and Sarin(1979)}]{DyerSa1979}
\bibinfo{author}{Dyer, J.}, \bibinfo{author}{Sarin, R.}, \bibinfo{year}{1979}.
\newblock \bibinfo{title}{Measurable multiattribute value functions}.
\newblock \bibinfo{journal}{Operations Research} \bibinfo{volume}{27},
  \bibinfo{pages}{810--822}.
\bibitem[{El~Gibari et~al.(2019)El~Gibari, G\'omez and Ruiz}]{ElGibariEtAl2019}
\bibinfo{author}{El~Gibari, S.}, \bibinfo{author}{G\'omez, T.},
  \bibinfo{author}{Ruiz, F.}, \bibinfo{year}{2019}.
\newblock \bibinfo{title}{Building composite indicators using multicriteria
  methods: {A} review}.
\newblock \bibinfo{journal}{Journal of Business Economics}
  \bibinfo{volume}{89}, \bibinfo{pages}{1--24}.
\bibitem[{Figueira et~al.(2021)Figueira, Greco and Roy}]{FigueiraEtAl2021}
\bibinfo{author}{Figueira, {\relax J.R.}.}, \bibinfo{author}{Greco, S.},
  \bibinfo{author}{Roy, B.}, \bibinfo{year}{2021}.
\newblock \bibinfo{title}{{\sc{Electre-Score}}: {A} first outranking based
  method for scoring actions}.
\newblock \bibinfo{journal}{European Journal of Operational Research}
  \bibinfo{note}{DOI: 10.1016/j.ejor.2021.05.017}.
\bibitem[{Figueira et~al.(2009)Figueira, Greco and
  S{\l}owi{\'n}ski}]{FigueiraEtAl2009}
\bibinfo{author}{Figueira, {\relax J.R.}.}, \bibinfo{author}{Greco, S.},
  \bibinfo{author}{S{\l}owi{\'n}ski, R.}, \bibinfo{year}{2009}.
\newblock \bibinfo{title}{Building a set of additive value functions
  representing a reference preorder and intensities of preference: {GRIP}
  method}.
\newblock \bibinfo{journal}{European Journal of Operational Research}
  \bibinfo{volume}{195}, \bibinfo{pages}{460--486}.
\bibitem[{Figueira et~al.(2016)Figueira, Mousseau and Roy}]{FigueiraEtAl2016}
\bibinfo{author}{Figueira, {\relax J.R.}.}, \bibinfo{author}{Mousseau, V.},
  \bibinfo{author}{Roy, B.}, \bibinfo{year}{2016}.
\newblock \bibinfo{title}{{\sc{Electre}} methods}, in: \bibinfo{editor}{Greco,
  S.}, \bibinfo{editor}{Ehrgott, M.}, \bibinfo{editor}{Figueira, J.} (Eds.),
  \bibinfo{booktitle}{Multiple Criteria Decision Analysis}.
  \bibinfo{publisher}{Springer Science+Business Media New York},
  \bibinfo{address}{New York, NY, USA}, pp. \bibinfo{pages}{155--185}.
\bibitem[{Figueira and Roy(2002)}]{FigueiraRo2002}
\bibinfo{author}{Figueira, {\relax J.R.}.}, \bibinfo{author}{Roy, B.},
  \bibinfo{year}{2002}.
\newblock \bibinfo{title}{Determining the weights of criteria in the
  {E}{\sc{lectre}} type methods with a revised {S}imos' procedure}.
\newblock \bibinfo{journal}{European Journal of Operational Research}
  \bibinfo{volume}{139}, \bibinfo{pages}{317--326}.
\bibitem[{Ghimire et~al.(2021)Ghimire, Parajuli, Khatiwada, Poudel, Sharma and
  Mishra}]{GhimireEtAl2021}
\bibinfo{author}{Ghimire, B.}, \bibinfo{author}{Parajuli, R.},
  \bibinfo{author}{Khatiwada, B.}, \bibinfo{author}{Poudel, S.},
  \bibinfo{author}{Sharma, K.}, \bibinfo{author}{Mishra, B.},
  \bibinfo{year}{2021}.
\newblock \bibinfo{title}{Covira: {A} {COVID-19} risk assessment, visualization
  and communication tool}.
\newblock \bibinfo{journal}{SoftwareX} \bibinfo{volume}{16}.
\bibitem[{Grabisch and Labreuche(2016)}]{GrabischLa2016}
\bibinfo{author}{Grabisch, M.}, \bibinfo{author}{Labreuche, {\relax Ch.}.},
  \bibinfo{year}{2016}.
\newblock \bibinfo{title}{Fuzzy measures and integrals in {MCDA}}, in:
  \bibinfo{editor}{Greco, S.}, \bibinfo{editor}{Ehrgott, M.},
  \bibinfo{editor}{Figueira, J.} (Eds.), \bibinfo{booktitle}{Multiple Criteria
  Decision Analysis}. \bibinfo{publisher}{Springer Science+Business Media New
  York}, \bibinfo{address}{New York, NY, USA}, pp. \bibinfo{pages}{553--603}.
\bibitem[{Greco et~al.(2016)Greco, Matarazzo and
  S{\l}owi{\'{n}}ski}]{GrecoEtAl2016}
\bibinfo{author}{Greco, S.}, \bibinfo{author}{Matarazzo, B.},
  \bibinfo{author}{S{\l}owi{\'{n}}ski, R.}, \bibinfo{year}{2016}.
\newblock \bibinfo{title}{Decision rule approach}, in: \bibinfo{editor}{Greco,
  S.}, \bibinfo{editor}{Ehrgott, M.}, \bibinfo{editor}{Figueira, J.} (Eds.),
  \bibinfo{booktitle}{Multiple Criteria Decision Analysis}.
  \bibinfo{publisher}{Springer Science+Business Media New York},
  \bibinfo{address}{New York, NY, USA}, pp. \bibinfo{pages}{497--552}.
\bibitem[{Haghighat(2021)}]{Haghighat2021}
\bibinfo{author}{Haghighat, F.}, \bibinfo{year}{2021}.
\newblock \bibinfo{title}{Predicting the trend of indicators related to
  {Covid-19} using the combined {MLP-MC} model}.
\newblock \bibinfo{journal}{Chaos, Solitons and Fractals}
  \bibinfo{volume}{152}, \bibinfo{pages}{111399}.
\bibitem[{Hale et~al.(2021)Hale, Angrist, Goldszmidt, Kira, Petherick,
  Phillips, Webster, Cameron-Blake, Hallas, Majumdar and Tatlow}]{HaleEtAl2021}
\bibinfo{author}{Hale, T.}, \bibinfo{author}{Angrist, N.},
  \bibinfo{author}{Goldszmidt, R.}, \bibinfo{author}{Kira, B.},
  \bibinfo{author}{Petherick, A.}, \bibinfo{author}{Phillips, T.},
  \bibinfo{author}{Webster, S.}, \bibinfo{author}{Cameron-Blake, E.},
  \bibinfo{author}{Hallas, L.}, \bibinfo{author}{Majumdar, S.},
  \bibinfo{author}{Tatlow, H.}, \bibinfo{year}{2021}.
\newblock \bibinfo{title}{A global panel database of pandemic policies
  ({O}xford {COVID-19} government response tracker)}.
\newblock \bibinfo{journal}{Nature Human Behaviour} \bibinfo{volume}{5},
  \bibinfo{pages}{529--538}.
\bibitem[{Keeney and Raiffa(1993)}]{KeeneyRa1993}
\bibinfo{author}{Keeney, R.L.}, \bibinfo{author}{Raiffa, H.},
  \bibinfo{year}{1993}.
\newblock \bibinfo{title}{Decisions with Multiple Objectives}.
\newblock \bibinfo{publisher}{Cambridge University Press},
  \bibinfo{address}{New York, NY, USA}.
\bibitem[{Keeney(1992)}]{Keeney1992}
\bibinfo{author}{Keeney, {\relax R.L.}.}, \bibinfo{year}{1992}.
\newblock \bibinfo{title}{Value-{F}ocused {T}hinking: {A} {Path} to {C}reative
  {Decisionmaking}}.
\newblock \bibinfo{publisher}{Harvard University Press},
  \bibinfo{address}{Cambridge, MA, USA}.
\bibitem[{Koch(2020)}]{Koch2020}
\bibinfo{author}{Koch, R.}, \bibinfo{year}{2020}.
\newblock \bibinfo{title}{Erl\''{a}uterung der sch\''{a} tzung der zeitlich
  variierenden reproduktionszahl r}.
\newblock \bibinfo{journal}{Memoir} \bibinfo{volume}{R/7-Tages-R (15.5.2020)}.
\bibitem[{Li et~al.(2021)Li, Rong and Zhang}]{LiEtAl2021}
\bibinfo{author}{Li, T.}, \bibinfo{author}{Rong, L.}, \bibinfo{author}{Zhang,
  A.}, \bibinfo{year}{2021}.
\newblock \bibinfo{title}{Assessing regional risk of {COVID-19} infection from
  wuhan via high-speed rail}.
\newblock \bibinfo{journal}{Transport Policy} \bibinfo{volume}{106},
  \bibinfo{pages}{226--238}.
\bibitem[{Martcheva(2015)}]{Martcheva2015}
\bibinfo{author}{Martcheva, M.}, \bibinfo{year}{2015}.
\newblock \bibinfo{title}{An Introduction to Mathematical Epidemiology}.
\newblock \bibinfo{publisher}{Springer Science+Business Media LLC New York},
  \bibinfo{address}{New York, NY, USA}.
\bibitem[{Martel and Matarazzo(2016)}]{MartelBe2016}
\bibinfo{author}{Martel, {\relax J.-M.}.}, \bibinfo{author}{Matarazzo, B.},
  \bibinfo{year}{2016}.
\newblock \bibinfo{title}{Other outranking approaches}, in:
  \bibinfo{editor}{Greco, S.}, \bibinfo{editor}{Ehrgott, M.},
  \bibinfo{editor}{Figueira, J.R.} (Eds.), \bibinfo{booktitle}{Multiple
  Criteria Decision Analysis}. \bibinfo{publisher}{Springer Science+Business
  Media New York}, \bibinfo{address}{New York, NY, USA}, pp.
  \bibinfo{pages}{221--282}.
\bibitem[{Moshkovich et~al.(2016)Moshkovich, Mechitov and
  Olson}]{MoshkovichEtAl2016}
\bibinfo{author}{Moshkovich, H.}, \bibinfo{author}{Mechitov, A.},
  \bibinfo{author}{Olson, D.}, \bibinfo{year}{2016}.
\newblock \bibinfo{title}{Verbal decision analysis}, in:
  \bibinfo{editor}{Greco, S.}, \bibinfo{editor}{Ehrgott, M.},
  \bibinfo{editor}{Figueira, J.} (Eds.), \bibinfo{booktitle}{Multiple Criteria
  Decision Analysis}. \bibinfo{publisher}{Springer Science+Business Media New
  York}, \bibinfo{address}{New York, NY, USA}, pp. \bibinfo{pages}{605--636}.
\bibitem[{Mousseau et~al.(2003)Mousseau, Figueira, Dias, Gomes~da Silva and
  {\relax J.C.N.}}]{MousseauEtAl2002}
\bibinfo{author}{Mousseau, V.}, \bibinfo{author}{Figueira, {\relax J.R.}.},
  \bibinfo{author}{Dias, {\relax L.C.}.}, \bibinfo{author}{Gomes~da Silva, C.},
  \bibinfo{author}{{\relax J.C.N.}, C.}, \bibinfo{year}{2003}.
\newblock \bibinfo{title}{Resolving inconsistencies among constraints on the
  parameters of an {MCDA} model}.
\newblock \bibinfo{journal}{European Journal of Operational Research}
  \bibinfo{volume}{147}, \bibinfo{pages}{72--93}.
\bibitem[{Nelken et~al.(2020)Nelken, Siems, Infantino, Genicot, Amariles and
  Harrington}]{NelkenEtAl2020}
\bibinfo{author}{Nelken, D.}, \bibinfo{author}{Siems, M.},
  \bibinfo{author}{Infantino, M.}, \bibinfo{author}{Genicot, N.},
  \bibinfo{author}{Amariles, D.}, \bibinfo{author}{Harrington, J.},
  \bibinfo{year}{2020}.
\newblock \bibinfo{title}{{COVID-19} and the social role of indicators: {A}
  preliminary assessment}.
\newblock \bibinfo{journal}{Workin Paper} \bibinfo{volume}{EUI Working Paper
  LAW 2020/17}.
\bibitem[{Neyens et~al.(2020)Neyens, Faes, Vranckx, Pepermans, Hens, Van~Damme,
  Molenberghs, Aerts and Beutels}]{NeyensEtAl2020}
\bibinfo{author}{Neyens, T.}, \bibinfo{author}{Faes, C.},
  \bibinfo{author}{Vranckx, M.}, \bibinfo{author}{Pepermans, K.},
  \bibinfo{author}{Hens, N.}, \bibinfo{author}{Van~Damme, P.},
  \bibinfo{author}{Molenberghs, G.}, \bibinfo{author}{Aerts, J.},
  \bibinfo{author}{Beutels, P.}, \bibinfo{year}{2020}.
\newblock \bibinfo{title}{Can covid-19 symptoms as reported in a large-scale
  online survey be used to optimise spatial predictions of covid-19 incidence
  risk in belgium?}
\newblock \bibinfo{journal}{Spatial and Spatio-temporal Epidemiology}
  \bibinfo{volume}{35}, \bibinfo{pages}{100379}.
\bibitem[{Pang et~al.(2021)Pang, Hu and Wen}]{PangEtAl2021}
\bibinfo{author}{Pang, S.}, \bibinfo{author}{Hu, X.}, \bibinfo{author}{Wen,
  Z.}, \bibinfo{year}{2021}.
\newblock \bibinfo{title}{Environmental risk assessment and comprehensive index
  model of disaster loss for {COVID-19} transmission}.
\newblock \bibinfo{journal}{Environmental Technology and Innovation}
  \bibinfo{volume}{23}, \bibinfo{pages}{101597}.
\bibitem[{Pictet and Bollinger(2008)}]{PictetBo2008}
\bibinfo{author}{Pictet, J.}, \bibinfo{author}{Bollinger, D.},
  \bibinfo{year}{2008}.
\newblock \bibinfo{title}{Extended use of the cards procedure as a simple
  elicitation technique for {MAVT.} application to public procurement in
  {S}witzerland}.
\newblock \bibinfo{journal}{European Journal of Operational Research}
  \bibinfo{volume}{185}, \bibinfo{pages}{1300--1307}.
\bibitem[{Roy(1996)}]{Roy1996}
\bibinfo{author}{Roy, B.}, \bibinfo{year}{1996}.
\newblock \bibinfo{title}{Multicriteria {M}ethodology for {D}ecision {A}iding}.
\newblock \bibinfo{publisher}{Kluwer Academic Publishers},
  \bibinfo{address}{Dordrecht, The Netherlands}.
\bibitem[{Roy(2010)}]{Roy2010}
\bibinfo{author}{Roy, B.}, \bibinfo{year}{2010}.
\newblock \bibinfo{title}{Robustness in operational research and decision
  aiding: {A} multi-faceted issue}.
\newblock \bibinfo{journal}{European Journal of Operational Research}
  \bibinfo{volume}{200}, \bibinfo{pages}{629--638}.
\bibitem[{Roy et~al.(2014)Roy, Figueira and Almeida-Dias}]{RoyEtAl2014}
\bibinfo{author}{Roy, B.}, \bibinfo{author}{Figueira, {\relax J.R.}.},
  \bibinfo{author}{Almeida-Dias, J.}, \bibinfo{year}{2014}.
\newblock \bibinfo{title}{Discriminating thresholds as a tool to cope with
  imperfect knowledge in multiple criteria decision aiding: {T}heoretical
  results and practical issues}.
\newblock \bibinfo{journal}{Omega - The International Journal of Management
  Science} \bibinfo{volume}{43}, \bibinfo{pages}{9--20}.
\bibitem[{Saaty(2016)}]{Saaty2016}
\bibinfo{author}{Saaty, {\relax T.L.}.}, \bibinfo{year}{2016}.
\newblock \bibinfo{title}{The {A}nalytic {H}ierarchy and {A}nalytic {N}etwork
  {P}rocesses for the measurement of intangible criteria and for
  decision-making}, in: \bibinfo{editor}{Greco, S.}, \bibinfo{editor}{Ehrgott,
  M.}, \bibinfo{editor}{Figueira, J.} (Eds.), \bibinfo{booktitle}{Multiple
  Criteria Decision Analysis}. \bibinfo{publisher}{Springer Science+Business
  Media New York}, \bibinfo{address}{New York, NY, USA}, pp.
  \bibinfo{pages}{363--419}.
\bibitem[{Sangiorgio and Parisi(2020)}]{SangiorgioPa2020}
\bibinfo{author}{Sangiorgio, V.}, \bibinfo{author}{Parisi, F.},
  \bibinfo{year}{2020}.
\newblock \bibinfo{title}{A multicriteria approach for risk assessment of
  {Covid-19} in urban district lockdown}.
\newblock \bibinfo{journal}{Safety Science} \bibinfo{volume}{130},
  \bibinfo{pages}{104862}.
\bibitem[{Sarkar and Chouhan(2021)}]{SarkarCh2021}
\bibinfo{author}{Sarkar, A.}, \bibinfo{author}{Chouhan, P.},
  \bibinfo{year}{2021}.
\newblock \bibinfo{title}{{COVID-19}: {D}istrict level vulnerability assessment
  in {I}ndia}.
\newblock \bibinfo{journal}{Clinical Epidemiology and Global Health}
  \bibinfo{volume}{9}, \bibinfo{pages}{204--215}.
\bibitem[{Shadeed and Alawna(2021)}]{ShadeedAl2021}
\bibinfo{author}{Shadeed, S.}, \bibinfo{author}{Alawna, S.},
  \bibinfo{year}{2021}.
\newblock \bibinfo{title}{{GIS}-based {COVID-19} vulnerability mapping in the
  {W}est {B}ank, {P}alestine}.
\newblock \bibinfo{journal}{International Journal of Disaster Risk Reduction}
  \bibinfo{volume}{64}, \bibinfo{pages}{102483}.
\bibitem[{Simos(1989)}]{Simos1989}
\bibinfo{author}{Simos, J.}, \bibinfo{year}{1989}.
\newblock \bibinfo{title}{L'\'Evaluation Environnementale: Un Processus
  Cognitif N\'egoci\'e}.
\newblock Ph.D. thesis. \'Ecole Polytechnique F\'ed\'erale de Lausanne (EPFL),
  Suisse.
\bibitem[{Siskos and Tsotsolas(2015)}]{SiskosTo2015}
\bibinfo{author}{Siskos, E.}, \bibinfo{author}{Tsotsolas, N.},
  \bibinfo{year}{2015}.
\newblock \bibinfo{title}{Elicitation of criteria importance weights through
  the {S}imos method: {A} robustness concern}.
\newblock \bibinfo{journal}{European Journal of Operational Research}
  \bibinfo{volume}{246}, \bibinfo{pages}{543--553}.
\bibitem[{Siskos et~al.(2016)Siskos, Grigoroudis and
  Matsatsinis}]{SiskosEtAl2016}
\bibinfo{author}{Siskos, Y.}, \bibinfo{author}{Grigoroudis, E.},
  \bibinfo{author}{Matsatsinis, {\relax N.F.}.}, \bibinfo{year}{2016}.
\newblock \bibinfo{title}{{UTA methods}}, in: \bibinfo{editor}{Greco, S.},
  \bibinfo{editor}{Ehrgott, M.}, \bibinfo{editor}{Figueira, J.} (Eds.),
  \bibinfo{booktitle}{Multiple Criteria Decision Analysis}.
  \bibinfo{publisher}{Springer Science+Business Media New York},
  \bibinfo{address}{New York, NY, USA}, pp. \bibinfo{pages}{315--361}.
\bibitem[{Tsotsolas et~al.(2019)Tsotsolas, Spyridakos, Siskos and
  Salmon}]{TsotsolasEtAl2019}
\bibinfo{author}{Tsotsolas, N.}, \bibinfo{author}{Spyridakos, A.},
  \bibinfo{author}{Siskos, E.}, \bibinfo{author}{Salmon, I.},
  \bibinfo{year}{2019}.
\newblock \bibinfo{title}{Criteria weights assessment through prioritizations
  {(WAP)} using linear programming techniques and visualizations}.
\newblock \bibinfo{journal}{Operational Research: An International Journal}
  \bibinfo{volume}{19}, \bibinfo{pages}{135--150}.
\newblock \bibinfo{note}{DOI: https://doi.org/10.1007/s12351-016-0280-7}.
\bibitem[{Tsoukiàs and Figueira(2006)}]{TsoukiasFi2006}
\bibinfo{author}{Tsoukiàs, A.}, \bibinfo{author}{Figueira, {\relax J.R.}.},
  \bibinfo{year}{2006}.
\newblock \bibinfo{title}{{RAND} {M}emorandum – 5868 by {H}oward {R}aiffa
  ({S}pecial {I}ssue), \textit{Journal of Multi-Criteria Decision Analysis},
  14(113)}.
\bibitem[{von Winterfeldt and Edwards(1986)}]{VonWinterfeldtEd1986}
\bibinfo{author}{von Winterfeldt, D.}, \bibinfo{author}{Edwards, W.},
  \bibinfo{year}{1986}.
\newblock \bibinfo{title}{Decision Analysis and Behavioral Research}.
\newblock \bibinfo{publisher}{Cambridge University Press},
  \bibinfo{address}{New York, NYk, USA}.
\bibitem[{Zopounidis and Doumpos(2002)}]{ZopounidisDo2002}
\bibinfo{author}{Zopounidis, C.}, \bibinfo{author}{Doumpos, M.},
  \bibinfo{year}{2002}.
\newblock \bibinfo{title}{Multicriteria classification and sorting methods: {A}
  literature review}.
\newblock \bibinfo{journal}{European Journal of Operational Research}
  \bibinfo{volume}{138}, \bibinfo{pages}{229--246}.

\end{thebibliography}

\vfill\newpage

\section*{Appendix}\label{sec:appendix}
\addcontentsline{toc}{section}{\numberline{}Appendix}
\noindent In the first part of this appendix, we provide some elements regarding the value functions for criteria, $g_2$ to $g_4$.

\begin{enumerate}
  \item Criterion $g_2$ (Transmission - \c{TRANS}). The performance levels, after discretising the scale of g2, and the blank cards inserted in between consecutive levels, are as follows:
    \[
    \{[0,0.92]\} [0] \{.0.96\} [2] \{0.98\} [4] \{1.00\} [6] \{1.02\} [8] \{1.040\}
    \]

    As in criterion $g_1$, we also used further levels to understand the evolution of the number of blank cards inserted in between consecutive levels. This was very similar to the value function for criterion $g_1$. The piecewise linear value function obtained is presented as follows:

    \begin{equation}\label{eq:function_x2}
    v_2(x_{2t}) = \left\{
    \begin{array}{lcl}
      4.34783x_2    & \mbox{if} & x_{2t} \in [0.000,\,0.920[ \\
      600x_{2t}-548    & \mbox{if} & x_{2t} \in [0.920,\,0.940[ \\
      1000x_{2t}-924   & \mbox{if} & x_{2t} \in [0.940,\,0.960[ \\
      1400x_{2t}-1308  & \mbox{if} & x_{2t} \in [0.940,\,0.980[ \\
      1800x_{2t}-1700  & \mbox{if} & x_{2t} \in [0.980,\,1.000[ \\
      2200x_{2t}-2100  & \mbox{if} & x_{2t} \in [1.000,\,1.020[ \\
      2600x_{2t}-2508  & \mbox{if} & x_{2t} \in [1.020,\,1.034[ \\
         180        & \mbox{if} & x_{2t} \in [1.034,\,+\infty[ \\
    \end{array}
    \right.
    \end{equation}
    As in $g_1$, it also can be approximated by a quadratic function.
  \item Criterion $g_3$ (Lethality - \c{LETHA}). This  is  a  different  type  of  value  function  from  the  ones  constructed for criteria $g_1$ and $g_2$. When discretizing the scale range of criterion $g_3$ and asking the experts to add blank cards in between consecutive levels, they always considered the same number of blank cards. It means that this function is a linear function and the reason is obvious, as explained by the experts: one death is always very serious and does not depend on the place we are on the scale of this criterion, i.e., moving from one to two deaths has the same impact as moving from $49$ to $50$. The function can thus be presented as follows:
   \begin{equation}\label{eq:function_x3}
    v_3(x_{3t}) = \left\{
    \begin{array}{lcl}
      {\displaystyle \frac{200x_{3t}}{7.2}}    & \mbox{if} & x_{3t} \in [0,\,6.48[ \\
      & & \\
      180                                   & \mbox{if} & x_{3t} \in [6.48,\,+\infty[ \\
    \end{array}
    \right.
    \end{equation}
    The saturation level at 180 is used for making an upper level and limit the values of the indicator, and not because the number of deaths after a certain level has the same impact as the number of deaths leading to the saturation/emergency level. The increase in the number of deaths always has a strong impact in terms of the severity of the pandemic.
  \item Criterion $g_4$ (Number of patients admitted to wards - \sc{WARDS}). The scale of criterion $g_4$ is a discrete scale, which leads also to a discrete value function. The levels selected by the experts  and blank cards inserted in between consecutive levels are presented below:
   \[
    \{0\} [0] \{500\} [2] \{1000\} [4] \{1500\} [6] \{2000\} [8] \{2500\} [10] \{3000\} [12] \{3500\}
   \]
    The value function can be states as follows:
    \begin{equation}\label{eq:function_x4}
    v_4(x_{4t}) = \left\{
    \begin{array}{lcl}
        0.008x_{4t}      & \mbox{if} & x_{4t} \in \{0,1,\ldots,498,499\} \\
        0.024x_{4t}-8    & \mbox{if} & x_{4t} \in \{500,501,\ldots,999,999\} \\
        0.04x_{4t}-24    & \mbox{if} & x_{4t} \in \{1000,1001,\ldots,1498,1499\} \\
        0.056x_{4t}-48   & \mbox{if} & x_{4t} \in \{1500,1501,\ldots,1999,1999\} \\
        0.072x_{4t}-80   & \mbox{if} & x_{4t} \in \{2000,2001,\ldots,2498,2499\} \\
        0.088x_{4t}-120  & \mbox{if} & x_{4t} \in \{2500,2501,\ldots,2998,2999\} \\
        0.104x_{4t}-168  & \mbox{if} & x_{4t} \in \{3000,3001,\ldots,3344,3345\} \\
        180           & \mbox{if} & x_{4t} \in \{3346,3347,\ldots\} \\
    \end{array}
    \right.
    \end{equation}
    It is a discrete function, but it has the same kind of “shape” and behaviour as the functions for criteria $g_1$ and $g_2$.
  \item Criterion $g_5$ (Number of patients admitted to ICU - \sc{ICU}). As in the previous case, this scale is also a discrete one, which also leads to a discrete value function. The levels selected by the experts and the number of blank cards inserted in between consecutive levels are presented below:
    \[
    \{0\} [0] \{40\} [2] \{80\} [4] \{120\} [6] \{160\} [8] \{200\} [10] \{240\} [12] \{280\}
    \]
    From the previous information and useing PaCo-DCM we can derive the following function.
   \begin{equation}\label{eq:function_x5}
    v_5(x_{5t}) = \left\{
    \begin{array}{lcl}
       0.1x_{5t}   & \mbox{if} & x_{5t} \in \{0,1,\ldots,39,40\} \\
       0.3x_{5t}-8    & \mbox{if} & x_{5t} \in \{40,41,\ldots,78,79\} \\
       0.5x_{5t}-24   & \mbox{if} & x_{5t} \in \{80,81,\ldots,118,119\} \\
       0.7x_{5t}-48   & \mbox{if} & x_{5t} \in \{120,121,\ldots,158,159\} \\
       0.9x_{5t}-80   & \mbox{if} & x_{5t} \in \{160,161,\ldots,198,199\} \\
       1.1x_{5t}-120  & \mbox{if} & x_{5t} \in \{200,201,\ldots,238,239\} \\
       1.3x_{5t}-168  & \mbox{if} & x_{5t} \in \{240,241,\ldots,266,267\} \\
         180       & \mbox{if} & x_{5t} \in \{268,269,\ldots\} \\
    \end{array}
    \right.
    \end{equation}
    Its ``shape'' and behaviour are similar to the previous value function.
\end{enumerate}



    \begin{figure}[htbp]
        \centering
        \begin{subfigure}[b]{0.7\textwidth}
            \centering
            \includegraphics[width=14cm, height=8.5cm]{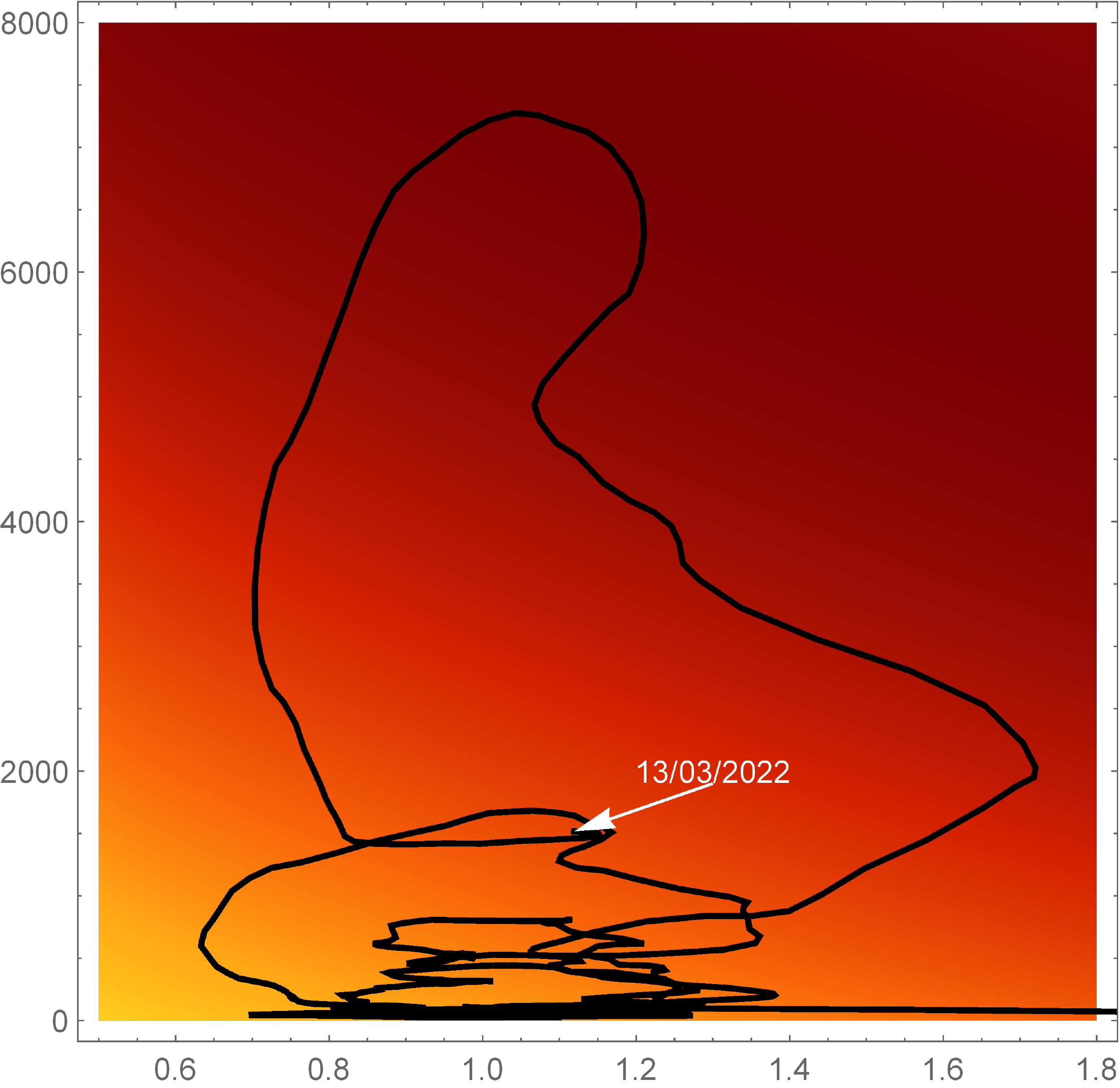}
            \caption[]{RM of Portuguese Health Authorities}
            \label{fig: OM_HA_Matrix}
        \end{subfigure}
        \\ \vspace{0.5cm}
        \begin{subfigure}[b]{0.7\textwidth}
            \centering
            \includegraphics[width=14cm, height=8.5cm]{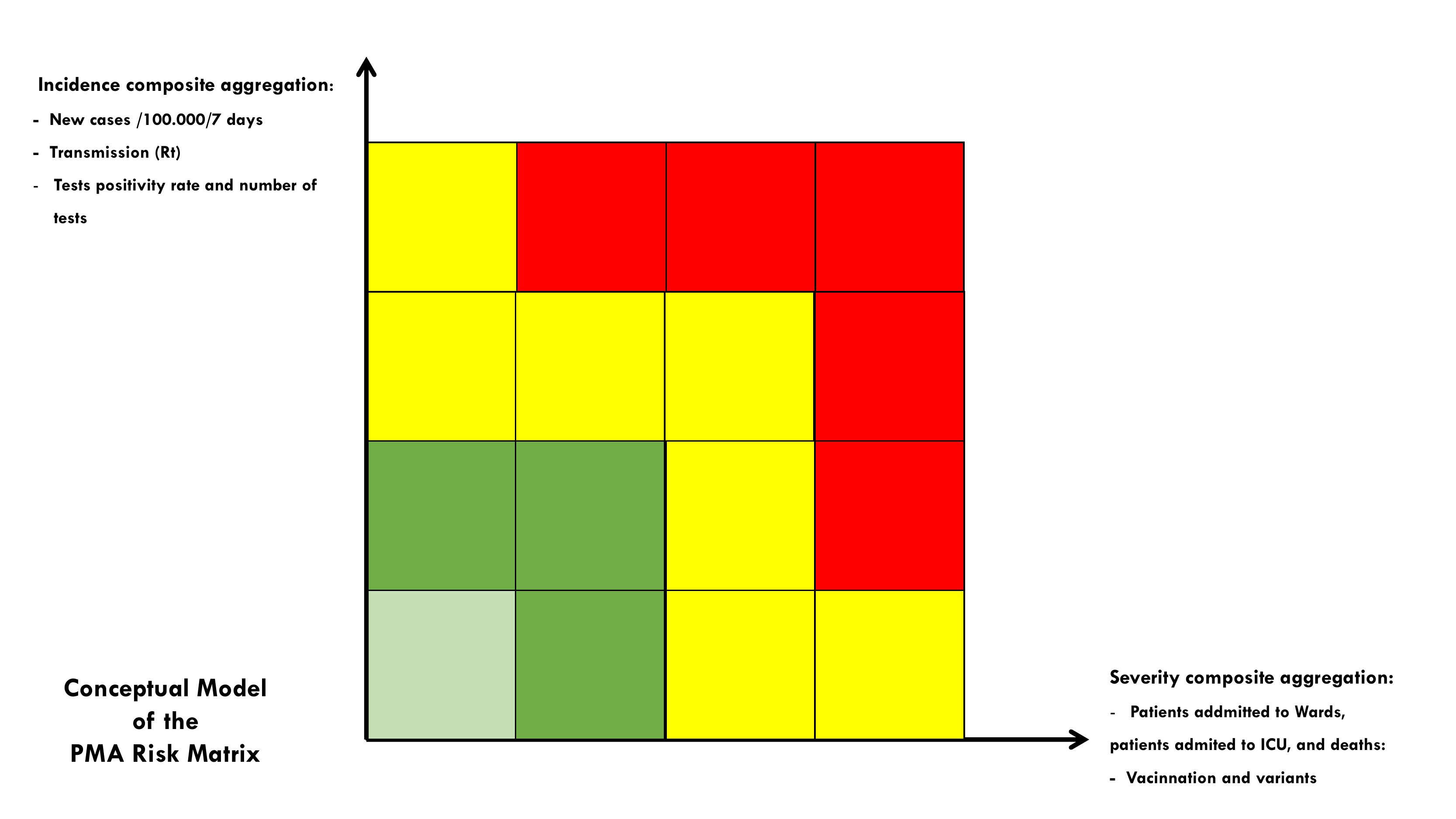}
            \caption[]{RM of the Portuguese Medical Association}
            \label{fig: RM_OM}
        \end{subfigure}
        \caption[]{Previous RM Tools}
        \label{fig: RM_Previews_}
    \end{figure}


\vfill\newpage

~~~~\\

\vspace{1.5cm}

    \begin{figure}[htbp]
        \centering
        \begin{subfigure}[b]{0.45\textwidth}
            \centering
            \includegraphics[width=9cm, height=6.5cm]{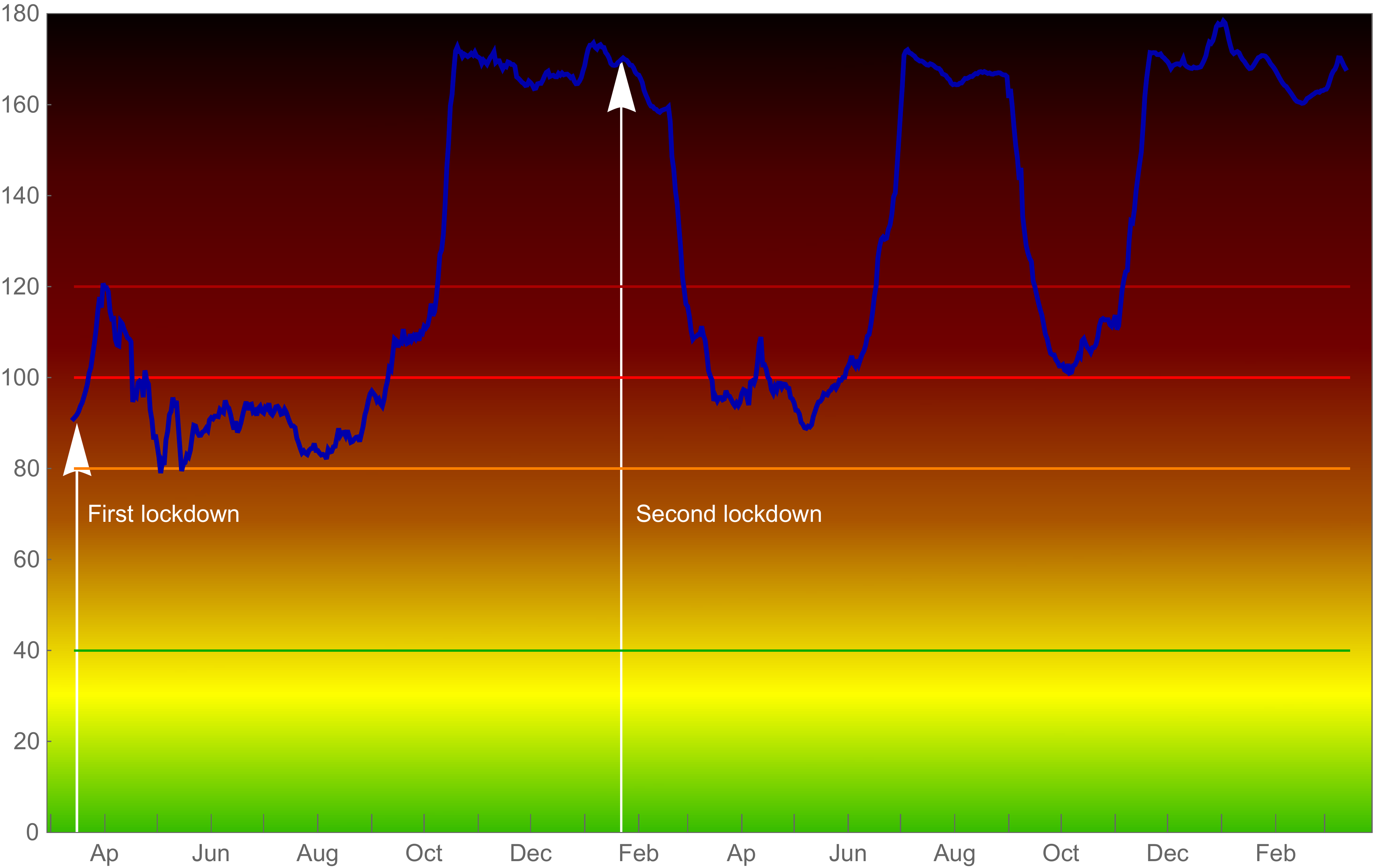}
            \caption[]{RM indicator evolution (2022-03-13)}
            \label{fig: RM_Evolution}
        \end{subfigure}
        \hfill
        \begin{subfigure}[b]{0.45\textwidth}
            \centering
            \includegraphics[width=9cm, height=6.5cm]{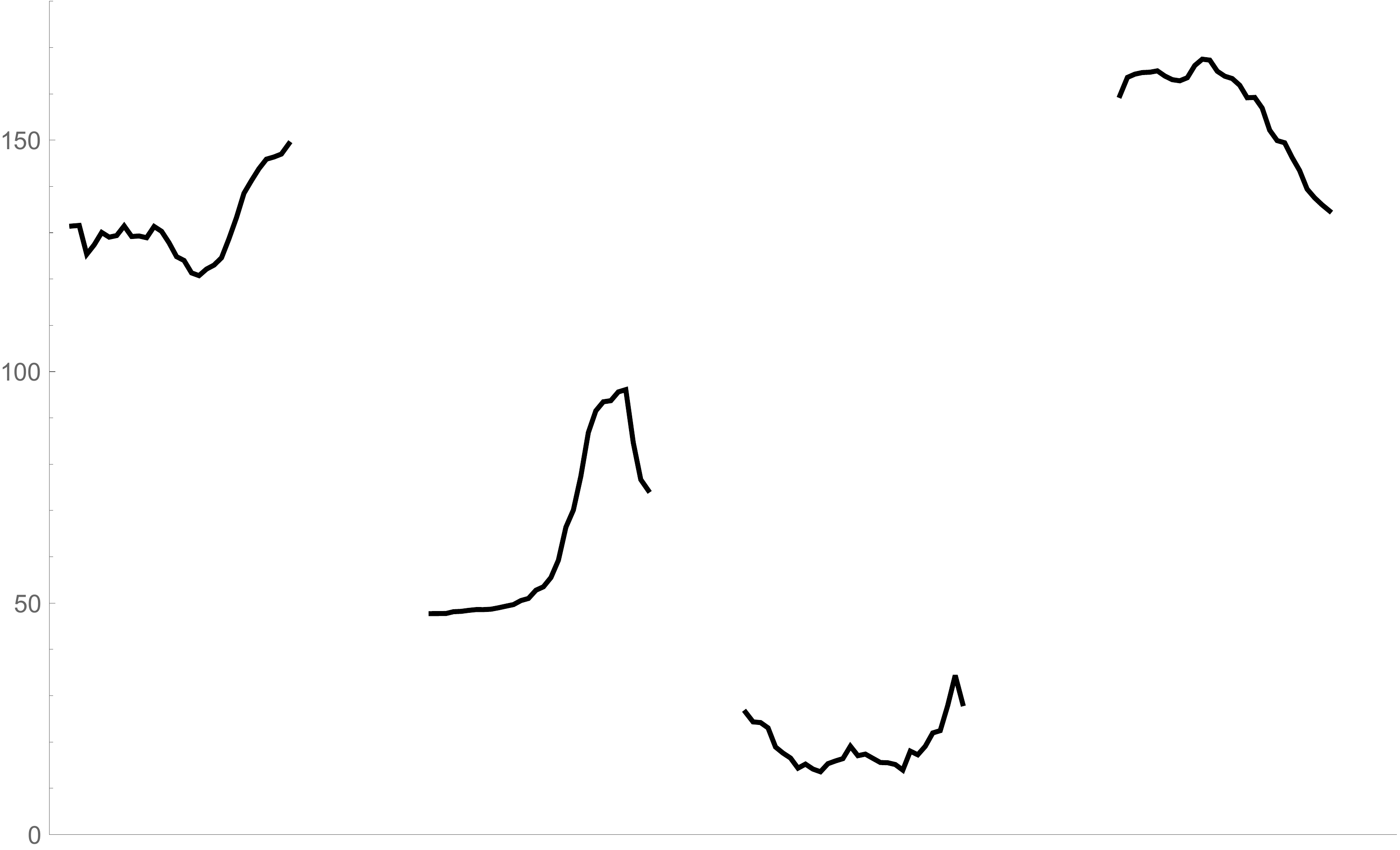}
            \caption[]{Four unordered moments of the evolution line}
            \label{fig: pieces}
        \end{subfigure}
        \\ \vspace{1.5cm}
        \begin{subfigure}[b]{0.45\textwidth}
            \centering
            \includegraphics[width=9cm, height=6.5cm]{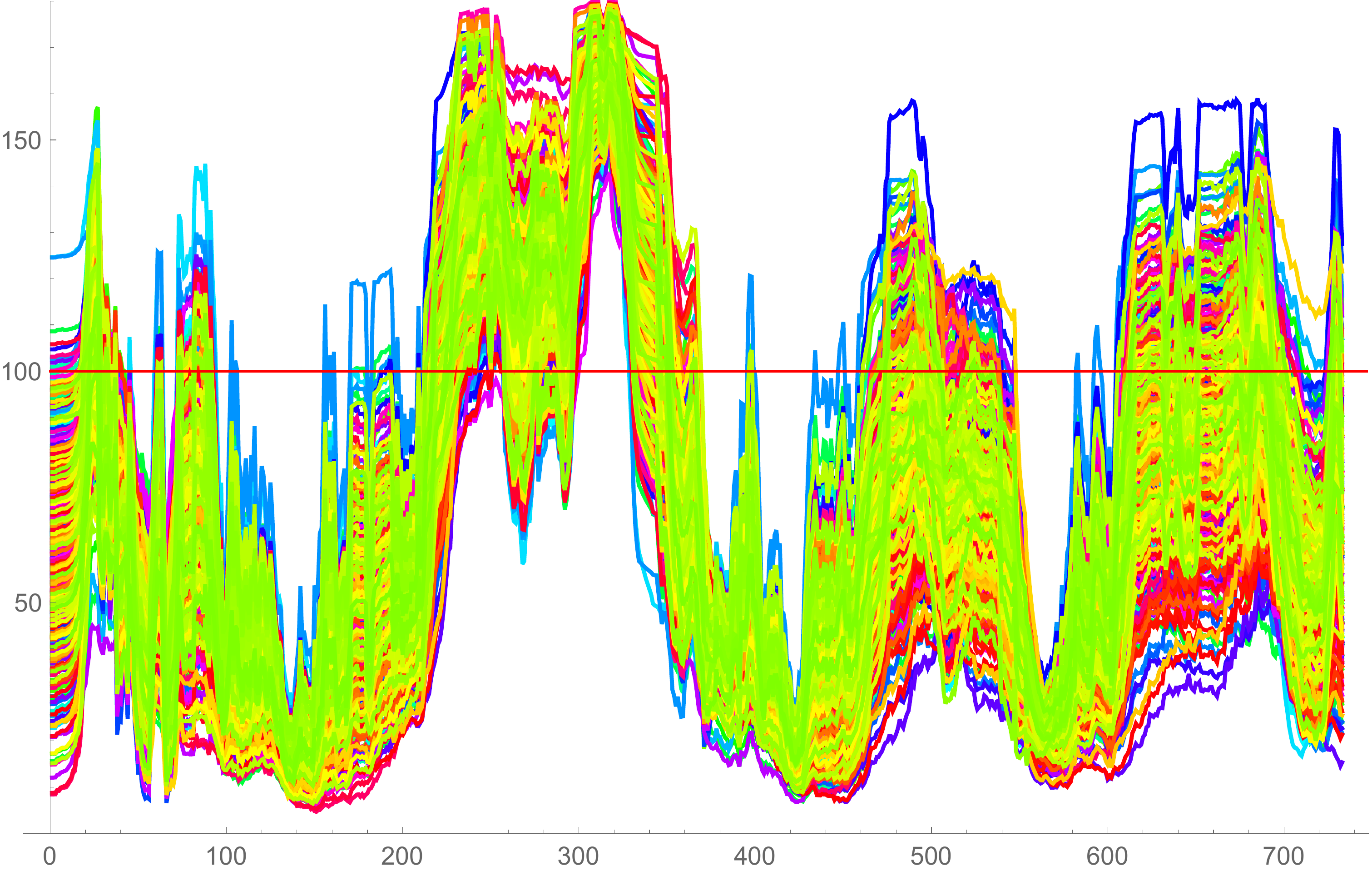}
            \caption[]{Simulating a drastic change in weights (2022-03-13)}
            \label{fig: simul_1}
        \end{subfigure}
        \hfill
        \begin{subfigure}[b]{0.45\textwidth}
            \centering
            \includegraphics[width=9cm, height=6.5cm]{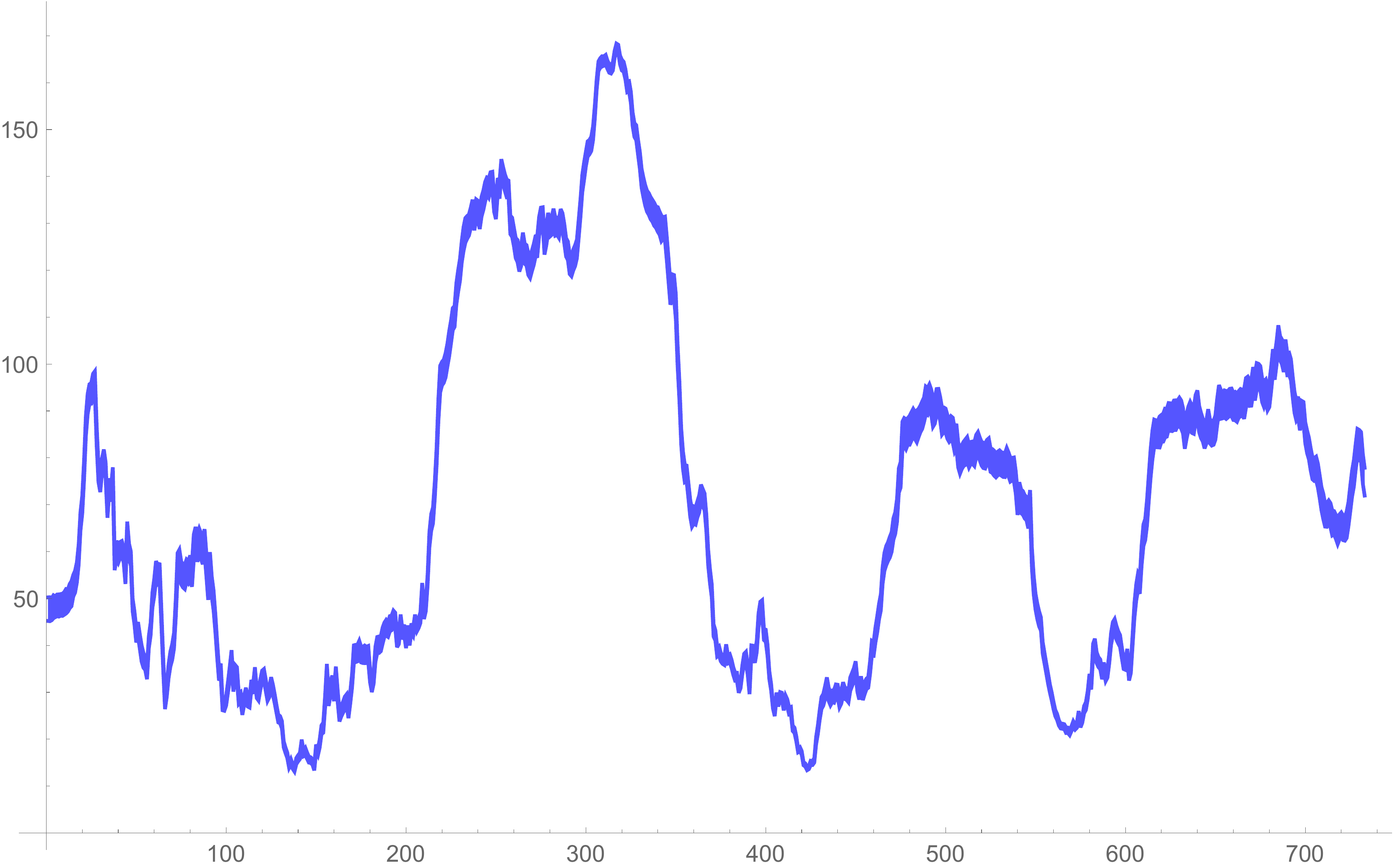}
            \caption[]{Simulating a realistic change in weights (2022-03-13)}
            \label{fig: simul_2}
        \end{subfigure}
        \caption[]{Validation and simulation analyses}
        \label{fig: valida_simul}
    \end{figure}

\vfill\newpage

~~~~\\

~~~~\\

\vspace{1.5cm}

    \begin{figure}[htbp]
        \centering
        \begin{subfigure}[b]{0.45\textwidth}
            \centering
            \includegraphics[width=9cm, height=6.5cm]{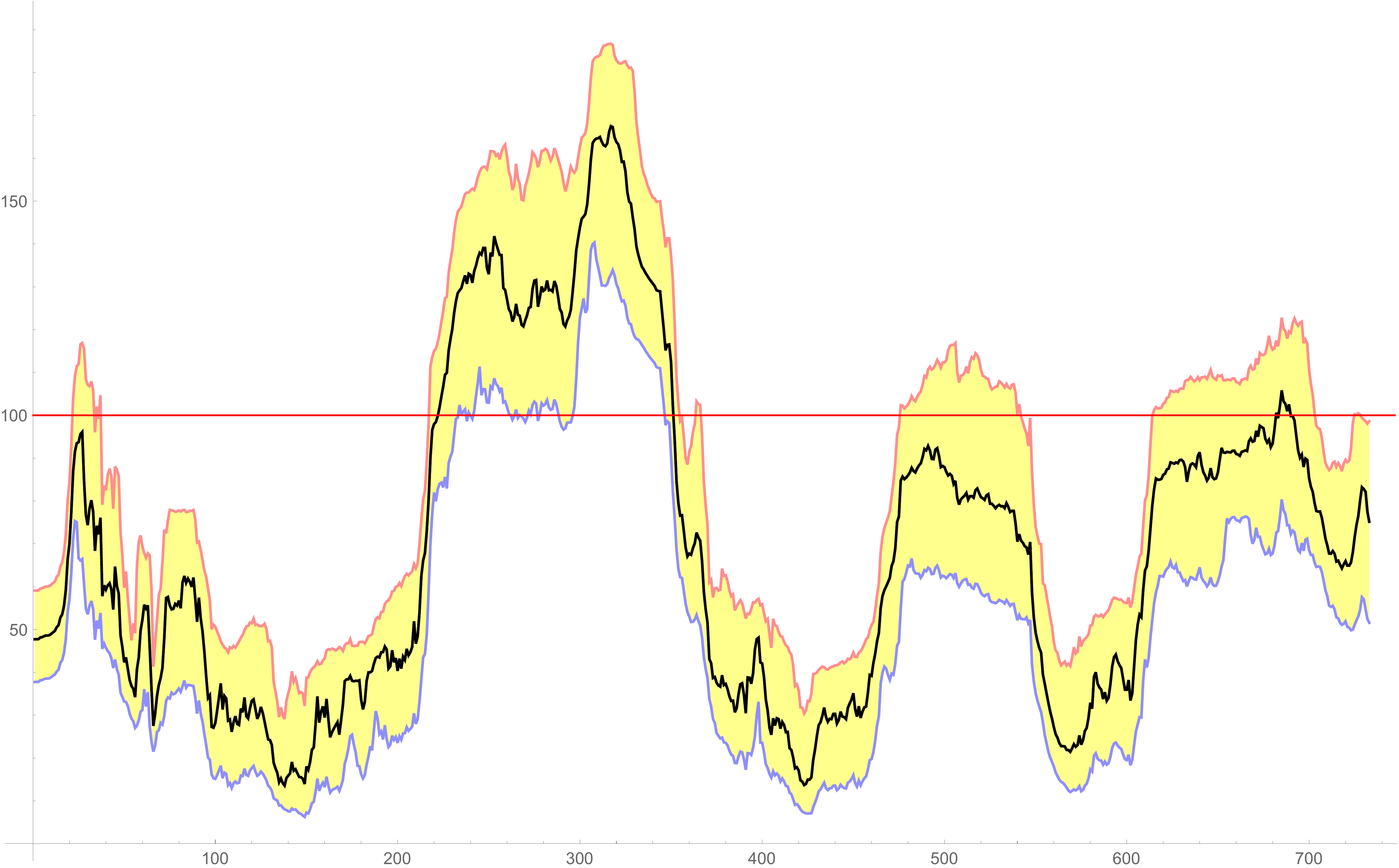}
            \caption[]{Sensitivity of a strong variation (2022-03-13)}
            \label{fig: rob_1}
        \end{subfigure}
        \hfill
        \begin{subfigure}[b]{0.45\textwidth}
            \centering
            \includegraphics[width=9cm, height=6.5cm]{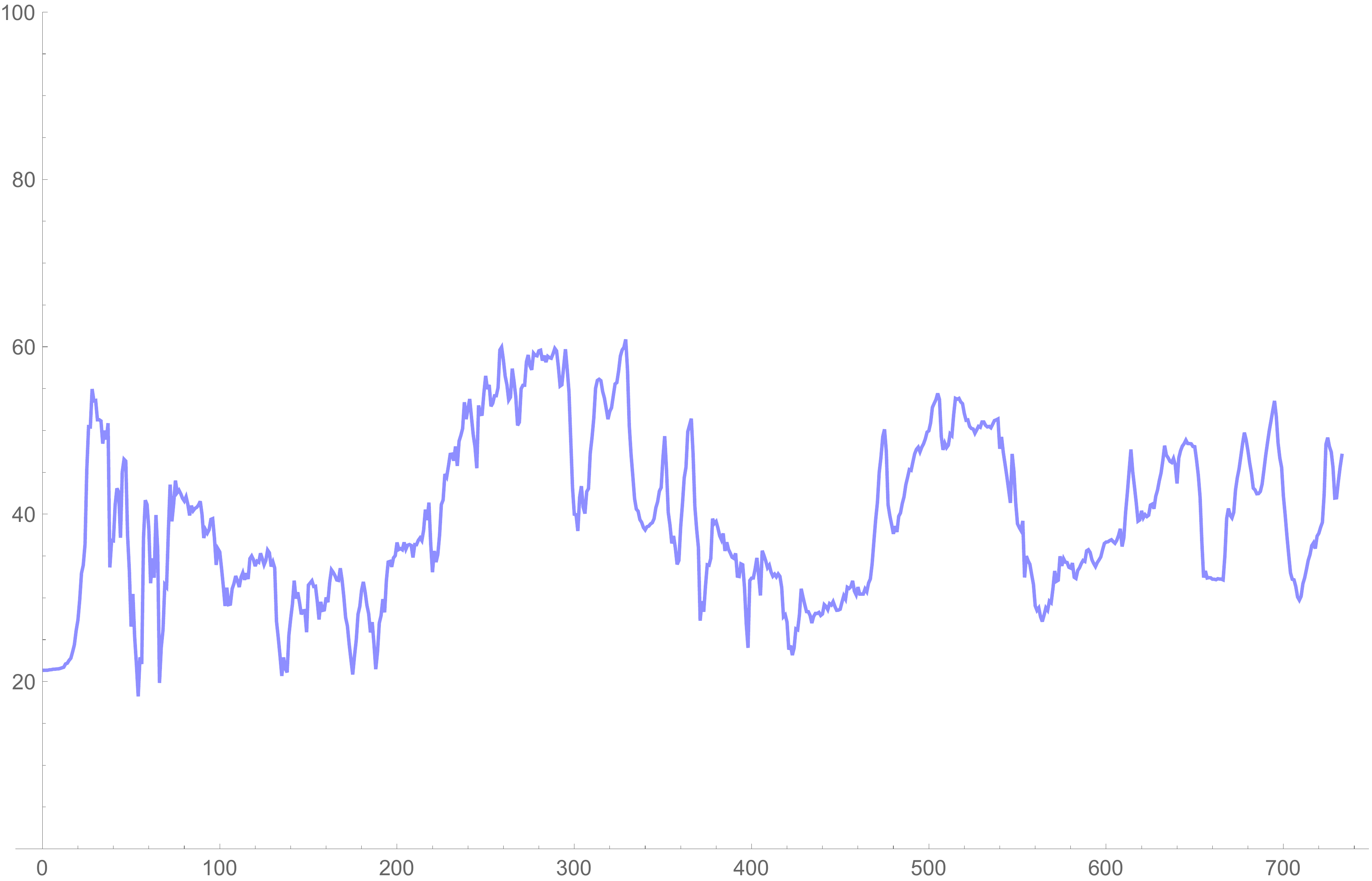}
            \caption[]{Average of the strong variation (2022-03-13)}
            \label{fig: rob_2}
        \end{subfigure}
        \\ \vspace{1.5cm}
        \begin{subfigure}[b]{0.45\textwidth}
            \centering
            \includegraphics[width=9cm, height=6.5cm]{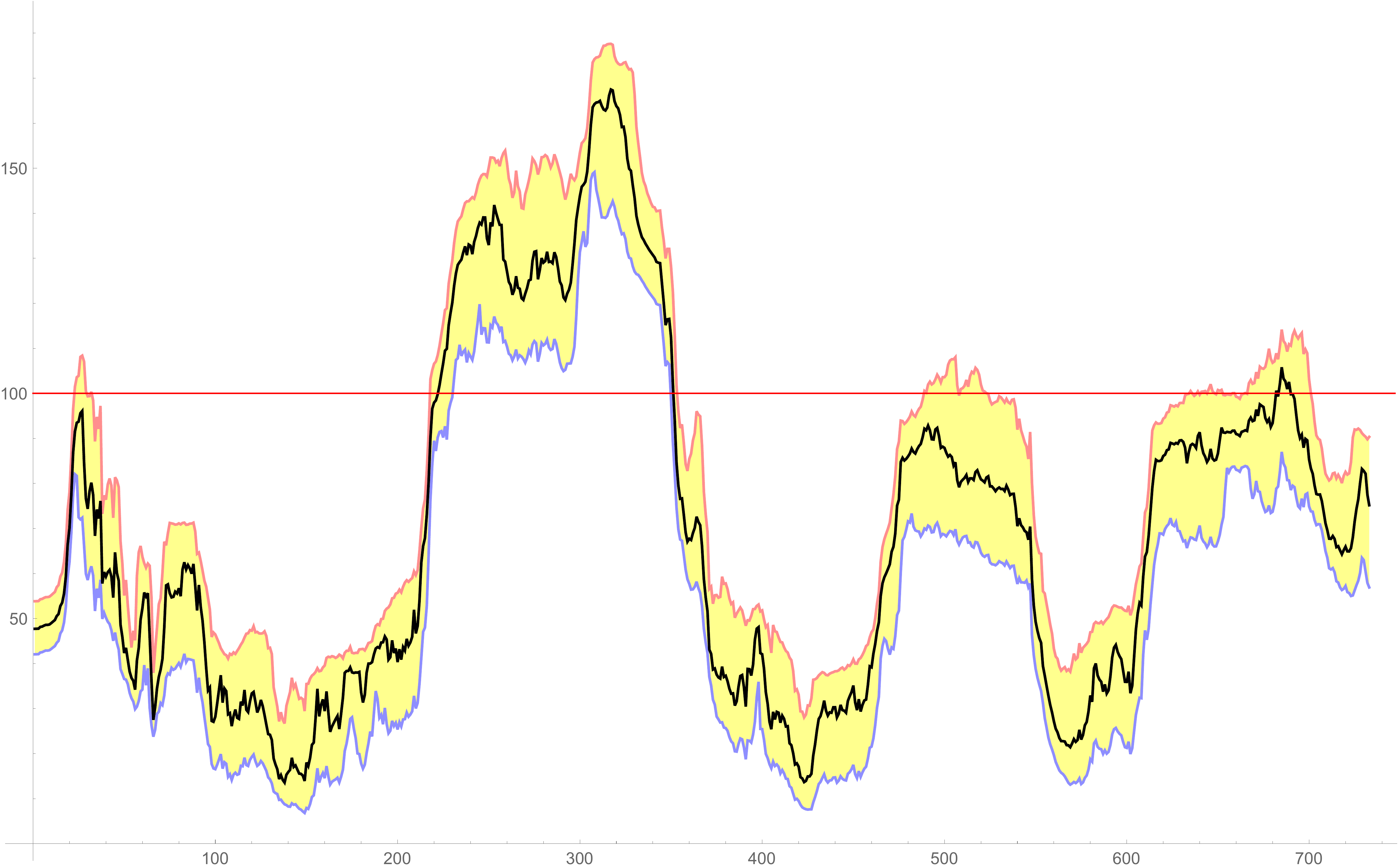}
            \caption[]{Sensitivity of a realistic variation (2022-03-13)}
            \label{fig: rob_3}
        \end{subfigure}
        \hfill
        \begin{subfigure}[b]{0.45\textwidth}
            \centering
            \includegraphics[width=9cm, height=6.5cm]{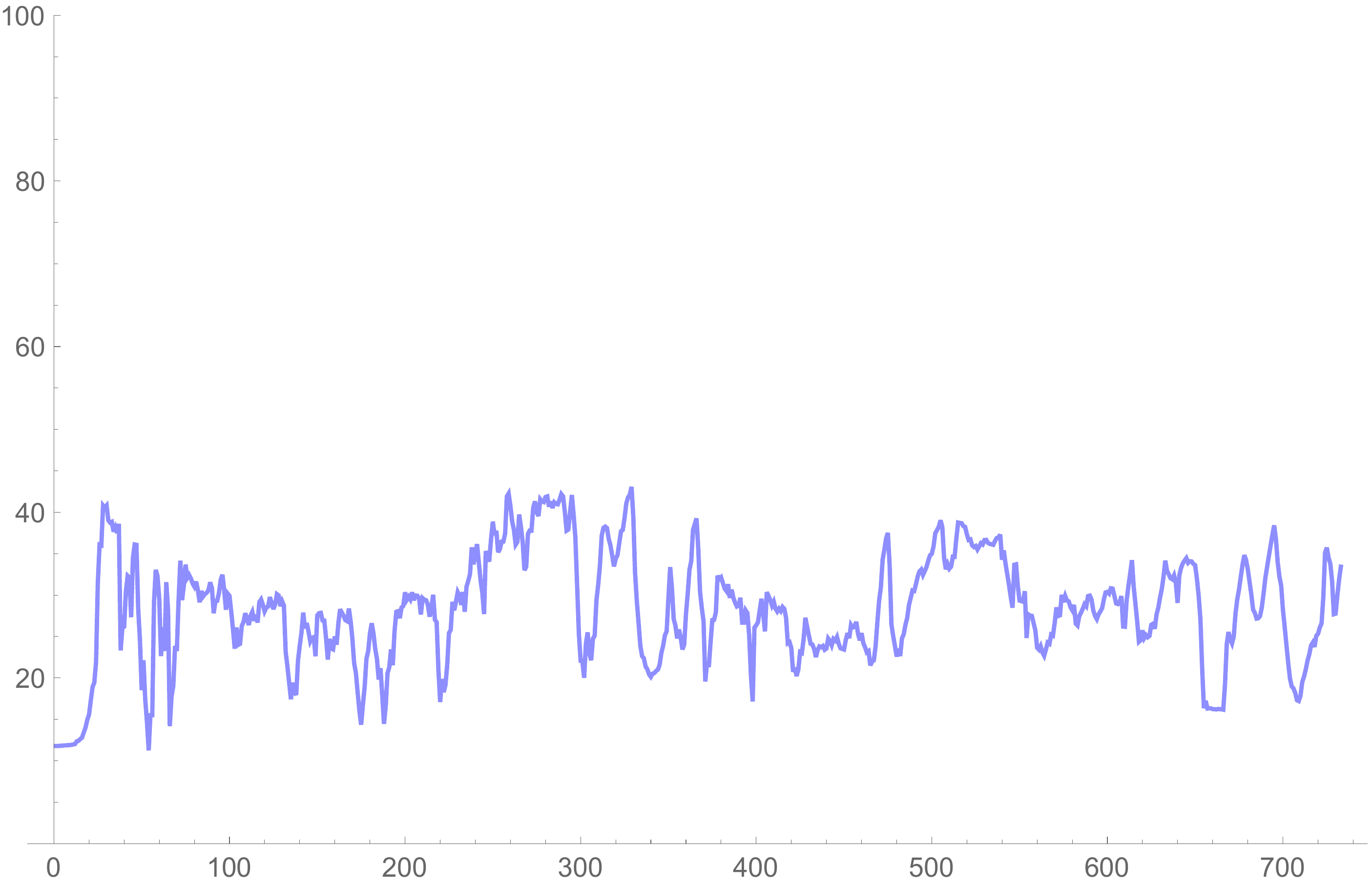}
            \caption[]{Average of the realistic variation (2022-03-13)}
            \label{fig: rob_4}
        \end{subfigure}
        \caption{sensitivity analysis}
        \label{fig: robustness}
    \end{figure}

\vspace{1.5cm}

~~~~\\

\vspace{2cm}

    \begin{figure}[htbp]
        \centering
        \includegraphics[width=9cm, height=6.5cm]{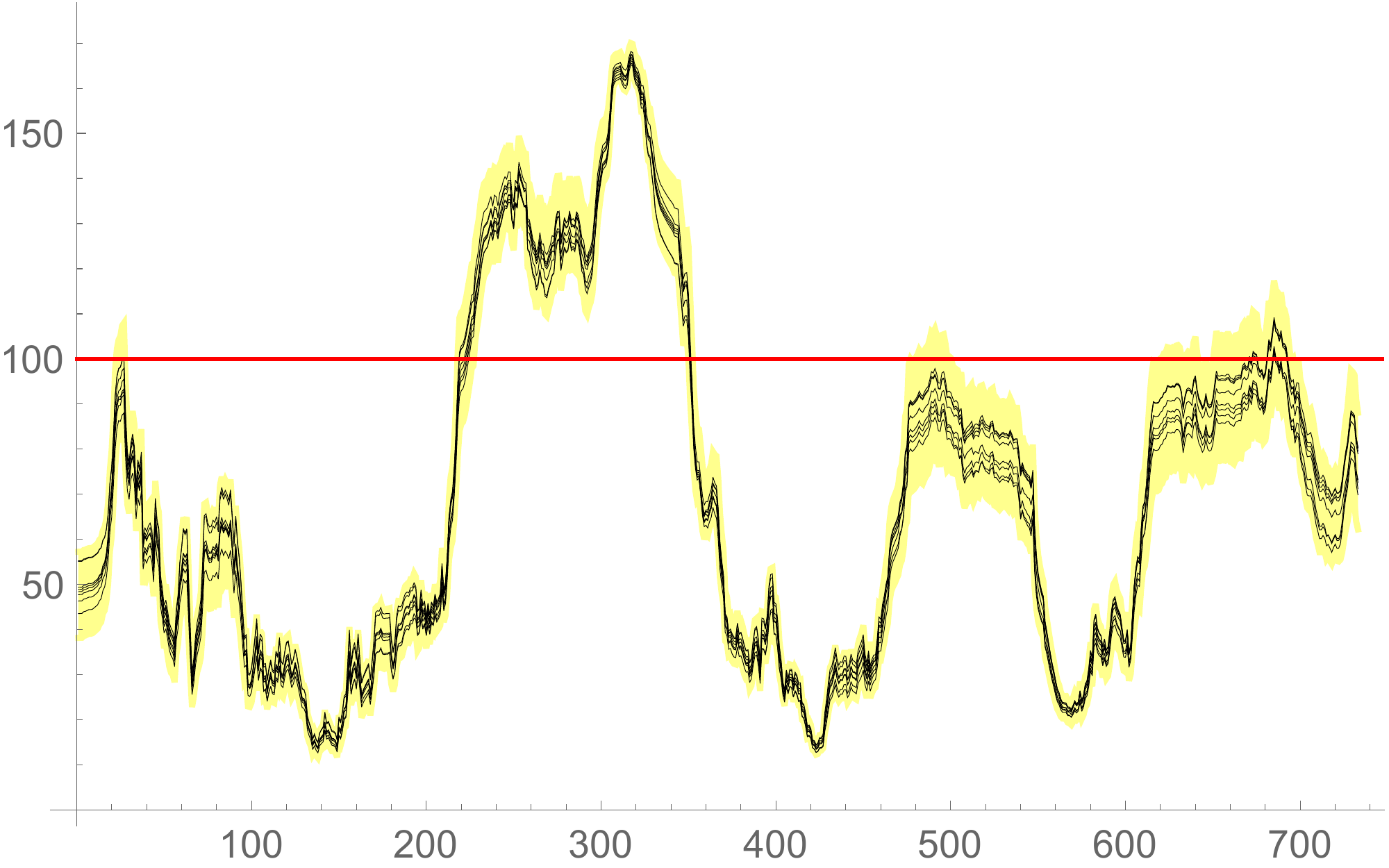}
        \caption[]{Simulation \textit{versus} Sensitivity (2022-03-13)}
        \label{fig: rob_add}
    \end{figure}

\end{document}